\documentclass[a4paper,10pt]{article}

\usepackage{import}

\usepackage[T1]{fontenc}
\usepackage[utf8]{inputenc}

\usepackage{mathrsfs}

\usepackage{framed}

\usepackage{xspace}

\usepackage[symbol]{footmisc}

\setfnsymbol{wiley}

\usepackage{etex}
\usepackage[english]{babel}
\usepackage{csquotes}

\usepackage[nocompress]{cite}

\usepackage{amsmath,amsfonts,amssymb}
\usepackage[standard,amsmath,thmmarks]{ntheorem}
\usepackage{enumitem}

\usepackage{appendix}
\usepackage{marginnote}

\usepackage{graphicx}
\usepackage[margin=0mm, font=footnotesize]{caption}

\usepackage[labelformat=simple]{subcaption}

\usepackage{tikz}
\usetikzlibrary{shapes,calc,decorations.pathreplacing,decorations.markings,arrows}
\usetikzlibrary{cd}

\usepackage{mathabx}

\usepackage[inner=30mm,outer=30mm,tmargin=25mm,bmargin=25mm]{geometry}

\usepackage{todonotes}

\usepackage{txfonts}
\usepackage{times}

\usepackage{bm}

\newcommand{\DBF}[1]{\expandafter\newcommand\csname #1\endcsname{{\mathbf{ #1 }}}}

\DBF{B}
\DBF{E}
\DBF{T}
\DBF{R}
\DBF{Q}
\DBF{N}
\DBF{Z}
\DBF{C}
\DBF{1}
\DBF{e}
\DBF{I}
\DBF{GL}
\DBF{SL}
\DBF{SO}
\DBF{n}
\DBF{F}

\newcommand{\bh}{\mathbf{h}}
\newcommand{\bg}{\mathbf{g}}
\newcommand{\bu}{\mathbf{u}}

\newcommand{\bz}{\mathbf{z}}
\newcommand{\bS}{\mathbf{S}}

\newcommand{\bPhi}{\mathbf{\Phi}}
\newcommand{\bPsi}{\mathbf{\Psi}}

\newcommand{\DCF}[1]{\expandafter\newcommand\csname cal#1\endcsname{{\mathcal{ #1 }}}}

\DCF{A}
\DCF{C}
\DCF{D}
\DCF{E}
\DCF{H}
\DCF{I}
\DCF{K}
\DCF{M}
\DCF{N}
\DCF{R}
\DCF{B}
\DCF{T}
\DCF{G}
\DCF{O}
\DCF{P}
\DCF{L}
\DCF{S}
\DCF{V}
\DCF{X}
\DCF{Y}
\DCF{F}
\DCF{W}
\DCF{Z}
\DCF{U}
\DCF{Q}

\newcommand{\DSF}[1]{\expandafter\newcommand\csname scr#1\endcsname{{{ #1 }}}}
\DSF{o}
\DSF{w}

\newcommand{\DFF}[1]{\expandafter\newcommand\csname frak#1\endcsname{{\mathfrak{ #1 }}}}

\DFF{L}
\DFF{g}
\DFF{o}
\DFF{so}
\DFF{X}

\newcommand{\DTF}[1]{\expandafter\newcommand\csname #1\endcsname {{\tt #1}}}

\newcommand {\DMO}[1]{\expandafter\DeclareMathOperator\csname #1\endcsname {#1}}

\DeclareMathOperator{\essinf}{ess~inf}
\DeclareMathOperator{\esssup}{ess~sup}

\DMO{conv}
\DMO{diag}
\DMO{supp}
\DMO{rank}
\DMO{id}
\DeclareMathOperator{\img}{im}

\DMO{diam}
\DMO{dom}
\DMO{coker}
\DMO{ind}
\DMO{Lie}
\DMO{sgn}
\DMO{crit}
\DMO{Jac}
\DMO{Hess}
\DMO{grad}
\DMO{rk}
\DMO{codim}
\DMO{SpecFlow}
\DMO{tr}
\DMO{Aut}
\DMO{Mat}
\DeclareMathOperator{\homology}{H}
\DMO{HC}
\DMO{HM}
\DMO{HTW}
\DeclareMathOperator{\HF}{HCF}

\DMO{vol}
\DMO{dist}
\DMO{card}

\DMO{Inv}
\DMO{BInv}

\DMO{div}
\DMO{period}
\DMO{SF}
\DMO{sq}

\newcommand{\sphere}{\mathbf{S}}

\mathchardef\hyphen="2D

\newcommand{\setmin}{\smallsetminus}

\newcommand{\iso}{\cong}

\newcommand\rst[2]{ \left. {#1} \right|_{#2} }

\newcommand{\eval}[3]{ \left. {#1} \right|_{#2}^{#3} }

\newcommand{\cl}[1]{ \ensuremath{\overline{#1}} }

\newcommand{\bdy}[1]{\partial {#1} }
\newcommand{\inter}[1]{\mathrm{int}\,{#1}}

\newcommand{\real}{\mathfrak{Re}\,}
\newcommand{\imag}{\mathfrak{Im}\,}

\renewcommand{\d}[1]{ \ensuremath{ \operatorname{d}\!{#1} } }

\newcommand{\set}[2]{ \left\{ #1 \; : \; #2 \right\} }

\renewcommand{\epsilon}{\varepsilon}
\renewcommand{\phi}{\varphi}

\renewcommand{\leq}{\leqslant}
\renewcommand{\geq}{\geqslant}
\renewcommand{\dots}{\ldots}

\renewcommand{\qed}{\hfill \ensuremath{\Box}}

\makeatletter
\newtheoremstyle{mynonumberplain}%
  {\item[\theorem@headerfont\hskip\labelsep ##1\theorem@separator]}%
  {\item[\theorem@headerfont\hskip \labelsep ##1\
    ##3\theorem@separator]}
\makeatother

\theoremstyle{plain}
\theoremprework{} 
\theorembodyfont{\itshape}
\theorempostwork{}
\theoremseparator{.}
\newcounter{dummy}
\numberwithin{dummy}{section}
\newtheorem{mydef}[dummy]{Definition}
\newtheorem{mythm}[dummy]{Theorem}
\renewtheorem*{mythm*}{Theorem}
\newtheorem{mylemma}[dummy]{Lemma}
\renewtheorem*{mylemma*}{Lemma}

\theoremnumbering{Alph}

\theoremnumbering{arabic}
\theoremsymbol{\ensuremath{\qed}}
\newtheorem{myremark}[dummy]{Remark}
\newtheorem{myexample}[dummy]{Example}
\theoremstyle{nonumberplain}
\theoremsymbol{}
\newtheorem{myhyp}{Hypothesis}

\theorembodyfont{\upshape}
\theoremsymbol{\ensuremath{\diamondsuit}}
\theoremseparator{.}
\theoremprework{}

\theoremstyle{mynonumberplain}
\theorembodyfont{\upshape}
\theoremsymbol{\rule{1.2ex}{1.2ex}}
\theoremseparator{.}
\newtheorem{myproof}{Proof}

\theoremstyle{nonumberplain}
\theoremheaderfont{\normalfont\bfseries}
\theoremseparator{.}
\theoremsymbol{}
\theorembodyfont{\slshape}

\makeatletter
\newtheoremstyle{mybreak}%
  {\item[\rlap{\vbox{\hbox{\hskip\labelsep \theorem@headerfont
          ##1\ ##2\theorem@separator}\hbox{\strut}}}]}%
  {\item[\rlap{\vbox{\hbox{\hskip\labelsep \theorem@headerfont
          ##1\ ##3,\ continued\theorem@separator}\hbox{\strut}}}]}
\makeatother

\theoremstyle{mybreak}

\makeatletter
\def\namedlabel#1#2{\begingroup
   \def\@currentlabel{\textrm{(#2)}}%
   \label{#1}\endgroup
}
\makeatother

\setlist[description]{leftmargin=2ex,labelindent=\parindent}

\setlist[itemize]{topsep=0pt}

\renewlist{enumerate}{enumerate}{1}
\setlist[enumerate]{label=(\arabic*),topsep=0pt}

\begin{document}


\title{Large fronts in nonlocally coupled systems using Conley--Floer homology}
\author{Bente Hilde Bakker\footnote{Mathematical Institute, Leiden University,
P.O. Box 9512, 2300 RA Leiden, The Netherlands. \href{mailto:bente.h.bakker@gmail.com}{\texttt{bente.h.bakker@gmail.com}}} 
\and Jan Bouwe van den Berg\footnote{Department of Mathematics, VU University Amsterdam, 
De Boelelaan 1081a, 1081 HV Amsterdam, The Netherlands. \href{mailto:janbouwe@few.vu.nl}{\texttt{janbouwe@few.vu.nl}}} \\[1.4 ex]
}
\date{\today}
\maketitle

\begin{abstract}
  In this paper we study travelling front solutions for nonlocal equations of the type
  \[
  \partial_t u = N * S(u) + \nabla F(u), \qquad u(t,x) \in \R^d.
  \]
  Here $N *$ denotes a convolution-type operator in the spatial variable $x \in \R$, either continuous or discrete.
  We develop a Morse-type theory, the Conley--Floer homology, which captures travelling front solutions in a topologically robust manner,
  by encoding fronts in the boundary operator of a chain complex.
  The equations describing the travelling fronts involve both forward and backward delay terms, possibly of infinite range.
  Consequently, these equations lack a natural phase space, 
  so that classic dynamical systems tools are not at our disposal.
  We therefore develop, from scratch, a general transversality theory, and a classification of bounded solutions, in the absence of a phase space.
  In various cases the resulting Conley--Floer homology can be interpreted as a homological Conley index for multivalued vector fields.
  Using the Conley--Floer homology we derive existence and multiplicity results on travelling front solutions.
\end{abstract}

\paragraph{\footnotesize 2010 Mathematics Subject Classification.}
{\footnotesize
Primary: 37L60, 45J05; \hspace{1 pt} Secondary: 53D40.
}

\paragraph{\footnotesize Acknowledgements.}
{\footnotesize
This research was supported by NWO TOP grant 613.001.351 and NWO VICI grant 639.033.109
}

\section{Introduction}
\label{sec:TWN_intro}

Evolution equations $u_t = A(u)$ arise throughout mathematical modelling.
Here the vector field $A$ may be given on finite or infinite dimensional spaces.
Vector fields on infinite-dimensional spaces are usually defined on function spaces through nonlinear functions of state variables and their spatial derivatives. 
Such vector fields are often obtained from a phenomenological point of view, where the spatial derivatives encode local interactions.
Rigorous derivation of such models as a continuum limit of microscopic laws can be daunting and requires that coupling between spatial sites $x$ and $y$ is only short-ranged.
In this paper we consider cases where this coupling may be long-ranged, so that the continuum limit retains \emph{nonlocal} interaction terms. 
More specifically, we consider here the class of evolution equations
\[
u_t = N * S(u) + \nabla F(u), \qquad u(t,x) \in \R^d,
\]
where $*$ denotes convolution in the spatial variable $x$, with either a continuous or discrete matrix-valued kernel $N$.
The convolution structure embodies the modelling assumption of translational invariance of physical space. 
In this paper we are concerned with front-type solutions in these equations.
We develop an index theory, the Conley--Floer homology, which captures front-type solutions for these nonlocal equations in a topologically robust manner, 
by \emph{encoding fronts in the boundary operator of a chain complex}.
For simplicity we consider $1$D fronts only, and assume $x \in \R$.
Conceptually, our methods generalise to higher dimensions and can deal with multidimensional fronts, comparable to \cite{bakker2017floer}, 
although compactness estimates appear to become much more delicate.

We refer to \cite{amari1977dynamics, ermentrout1993existence, pinto2001spatially, wilson1972excitatory, wilson1973mathematical} for examples of nonlocal models in neurobiology,
and \cite{bates1997traveling, bates1999integrodifferential, bates2006some, du2012analysis, du2013analysis, silling2000reformulation} for examples of nonlocal models in material sciences.
As a guiding principle, we present two examples of nonlocal equations to which our results are applicable.
\begin{myexample}[Neural field equations]
\label{ex:TWN_NFE}
Voltage-based neural field equations \cite{amari1977dynamics, ermentrout1993existence},
\begin{equation}
  \label{eq:TWN_NFE}
  \partial_t u_i = - u_i + \sum_{j=1}^d N_{ij} * \sigma_j(u_j), \qquad 1 \leq i \leq d,
\end{equation}
model interactions of $d$ thin layers of neurons through nonlocal synaptic connections.
Here $u_i$ represents the average membrane potential of neurons in layer $i$.
Typically, one thinks of the ``spatial variable'' $x$ as a feature space, grouping neurons by function and/or location.
The threshold $\sigma_i : \R \to \R$ is a smooth sigmoidal input function, which models the triggering of axonal firing of neurons in layer $i$ with membrane potential above a threshold value. 

The interaction kernels $N_{ij}(x-y)$ encode connections between neurons in layer $i$ at site $x$ and in layer $j$ at site $y$.
The sign of $N_{ij}$ encodes inhibitory versus excitatory connections.
Sign changes in $N_{ij}$ can lead to complex stationary patterns.
We make the somewhat specific assumption that the interaction kernels encode symmetric coupling, both in feature space and between layers,
i.e.,
\begin{equation}
  \label{eq:TWN_kernel_symmetry}
  N_{ij}(x) = N_{ij}(-x), \qquad \text{and} \quad N_{ij} = N_{ji}.
\end{equation}

Switching to vector notation, with the travelling wave ansatz $u(t,x) = u(\xi)$, $\xi = x - ct$ in \eqref{eq:TWN_NFE}, we obtain
\[
c u' + N*\sigma(u) - u = 0, \qquad u(\cdot) \in \R^d.
\]
Note that now both differentiation and the convolution operator act upon the spatial variable.
\end{myexample}

\begin{myexample}[Nonlocal reaction-diffusion]
\label{ex:TWN_nonlocal_diffusion}
A different class of examples are nonlocal reaction-diffusion equations of the gradient form
\[
u_t = N*u + \nabla F(u),
\]
where the matrix-valued interaction kernel $N$ may be continuous, discrete, or a combination of these two, and is assumed to be symmetric, as in \eqref{eq:TWN_kernel_symmetry}.
Equations of this type appear, for example, as models for phase transition \cite{bates1997traveling},
and as amplitude equations in neural field models \cite{bressloff2008nonlocal}.
Typically, the potential $F$ exhibits multiple local minima.
With the travelling wave ansatz $u(t,x) = u(\xi)$, $\xi = x - ct$, we obtain the integro-differential equation
\[
  c u' + N*u + \nabla F(u) = 0, \qquad u(\cdot) \in \R^d.
\]
Again, both differentiation and the convolution operator act upon the spatial variable.
\end{myexample}

For local equations, dynamical systems methods have provided powerful tools to the study of nonlinear waves.
One casts the governing equations as a dynamical system in the spatial variable $x$, 
which is then considered as a time-like variable, trying to describe the set of bounded solutions, including periodic, heteroclinic, and homoclinic solutions. 
The dynamical systems tools can then be thought of as pointwise methods, exploiting the geometry of the ``spatial phase space''.
Such ideas where already exploited by Conley, using index theory to study global bifurcations of travelling waves in local reaction-diffusion equations on the real line \cite{conley1971isolated, conley1980application}.

For nonlocal equations, much of the pointwise methods are not readily available. 
In fact,  it is often not clear how a pointwise formulation as a dynamical system in a phase space can be recovered,
since the nonlocal interactions in space typically generate forward- and backward delay terms in the time-like spatial variable. 
Much of the previous work has therefore relied on variational or perturbative techniques. 
More recently, certain dynamical systems techniques have been made available for the study of nonlocal equations. 
A first approach casts systems with finite-range interactions as ill-posed dynamical systems on an infinite-dimensional phase space, 
much like elliptic equations posed on a cylinder, 
and uses dynamical systems techniques, in particular invariant manifolds, to study the dynamics; see for instance \cite{hss,mpvl,hvl,hupkes}. 
A different avenue evokes dynamical systems techniques without ever setting up a phase space for a pointwise description, 
but rather adapting techniques from dynamical systems by reducing to basic functional analytic aspects; 
see \cite{faye2015existence} for geometric singular perturbation techniques, \cite{faye2016center}  for center manifold reductions, 
\cite{scheel2017bifurcation} for bifurcation methods, and \cite{bakker2018noether} for Hamiltonian identities.
The present work can be viewed as a continuation of this latter approach, providing a framework for topological index theory applicable to nonlocally coupled systems.

\subsection{Outline of main results}
Our theory is applicable to Examples \ref{ex:TWN_NFE} and \ref{ex:TWN_nonlocal_diffusion} with $c \neq 0$, i.e., travelling fronts.
After rescaling the spatial variable $x$, we may assume without loss of generality that $c = 1$.
These examples then fit within a more general class of equation of the type
\begin{equation}
  \label{eq:TWEnl}
   u' + \Phi(u) = 0, \qquad
    \Phi(u) = \nabla_g S(u)^T \calN[S(u)] + \nabla_g F(u).
\end{equation}
Here $u : \R \to \R^d$, and $\calN$ denotes a convolution-type operator with continuous and/or discrete kernels, i.e.,
\[
\calN[v](x) = \int_\R N(y) v(x-y) \d y + \sum_{j\in\Z} N_j v(x-\lambda_j),
\]
where the kernels are symmetric and exponentially localised. 
We emphasise that both the differential as well as the nonlocal operation $\calN$ act on the same variable $x$.
The nonlinear maps $S$ and $F$ are defined through pointwise evaluation.
The differential $\nabla_g$ denotes the gradient with respect to a Riemannian metric $g$.
Precise conditions will be formulated in Section \ref{sec:TWN_setup}.

Let us point out that to cast Example \ref{ex:TWN_NFE} into this setup, after rescaling so that $c = 1$, 
we choose $S(u) = (\sigma_1(u_1),\dots,\sigma_d(u_d))$ and define $g$ by setting, for each $u \in \R^d$,
\begin{equation*}
  g_u(v_1,v_2) := \big\langle v_1 , DS(u) v_2 \big\rangle, \qquad v_1, v_2 \in T_u \R^d = \R^d,
\end{equation*}
so that $\nabla_g S(u) = \id$.
Consequently, we choose $F$ so that $\nabla_g F(u) = -u$; such a potential $F$ exists in light of the diagonal structure of $g$.
For Example \ref{ex:TWN_nonlocal_diffusion} we simply set $S(u) = u$ and let $g$ be the flat metric on $\R^d$.

The travelling front solutions of Examples \ref{ex:TWN_NFE} and \ref{ex:TWN_nonlocal_diffusion} now correspond to solutions $u$ to \eqref{eq:TWEnl}
with $u(x) \to z_\pm$ as $x \to \pm\infty$, where $z_- \neq z_+$.
Denote by $\tau \cdot$ the shift action, defined through $\tau \cdot u(x) = u(x+\tau)$, for $\tau \in \R$.
Let $\Sigma$ be the set of all uniformly bounded solutions to \eqref{eq:TWEnl}, with the topology inherited of, say, $C^1_{\text{loc}}(\R,\R^d)$.
Since \eqref{eq:TWEnl} is left invariant by the shift action,
shifts formally define a dynamical system on $\Sigma$, where heteroclinic orbits in this system correspond to travelling front solutions.

In this paper we develop an index theory, the Conley--Floer homology, for this shift dynamical system.
The construction of the homology follows the recipe of Morse--Floer homology, but with a slight twist.
That is, \emph{we do not} obtain a variational formulation for heteroclinic solutions of \eqref{eq:TWEnl} and use these heteroclinics as the generators of the homology, 
as one familiar with Morse theory might expect.
Instead, we encode heteroclinic solutions of \eqref{eq:TWEnl} in the boundary operator of an appropriate chain complex.

\paragraph{Chain complex.}
We note that restricting $\Phi$ to the constant functions yields a vector field on $\R^d$ which may be written as $\nabla_g h$, where the potential $h$ is given by
\[
h(z) = \frac 1 2 S(z) \cdot \widetilde{\calN} S(z) + F(z), \qquad \widetilde{\calN} = \int_\R N(x) \d x + \sum_{j\in\Z} N_j.
\]
Note that to obtain this variational structure, we have used the imposed symmetry \eqref{eq:TWN_kernel_symmetry} of $\widetilde{\calN}$.
The critical points of $h$ are in 1-to-1 correspondence with the constant solutions of \eqref{eq:TWEnl}.
Assuming these are all nondegenerate, we take the critical points of $h$ to be the generators of the homology, graded using the classical Morse index of $h$.
In Section \ref{sec:TWN_functional} we relate these indices with Fredholm properties of \eqref{eq:TWEnl}, along with presenting various other functional properties of \eqref{eq:TWEnl}.
Moreover, in Section \ref{sec:TWN_heteroclinics} it is shown that solutions of \eqref{eq:TWEnl} satisfy uniform bounds on kinetic energy and possess a gradient-like dichotomy.
This allows us to define the boundary operator $\partial$ of the homology as the \emph{binary count of isolated heteroclinic solutions} of \eqref{eq:TWEnl}.
The gradient-like dichotomy is one of the essential ingredients in deriving the fundamental relation $\partial_n \partial_{n+1} = 0$.
We thus obtain a chain complex
\[
\begin{tikzcd}
  \cdots \arrow[r, "\partial_{n+1}"] & C_n \arrow[r, "\partial_{n}"] & C_{n-1} \arrow[r, "\partial_{n-1}"] & \cdots
\end{tikzcd}
\]
We may perform this construction whilst restricting our attention to solutions of \eqref{eq:TWEnl} which lie within certain subsets $E \subset L^\infty(\R,\R^d)$, 
which have nice isolating properties for \eqref{eq:TWEnl}; we refer to these sets as \emph{isolating trajectory neighbourhoods}, which are defined in Section \ref{subsec:TWN_isolating_nbhds}.

\paragraph{Transversality.}
There is one large caveat.
The counting of heteroclinic solutions of \eqref{eq:TWEnl}, needed in order to define the boundary operator $\partial$ and verify $\partial_n \partial_{n+1} = 0$,
is only defined under certain abstract perturbations of \eqref{eq:TWEnl}.
In the setting of Morse theory, where one has a global flow on a phase space, this condition corresponds to the transverse intersection of stable and unstable manifolds.
By analogy, we refer to abstract perturbations of \eqref{eq:TWEnl} which make the boundary operator well-defined over an isolating trajectory neighbourhood $E$ as \emph{$E$-transverse perturbations}.
Fundamental is the equivariance of these perturbations with respect to translations $u \mapsto \tau \cdot u$,
as the translation invariance of \eqref{eq:TWEnl} should be retained.

One of the main technical contributions of this paper is the derivation of a generic class of perturbations $\Psi$ of \eqref{eq:TWEnl}, equivariant with respect to translation, which are $E$-transverse.
These perturbations appear as an additive term in \eqref{eq:TWEnl}, leading to the perturbed equation
\begin{equation}
  \label{eq:TWN}
  \tag{TWN}
    u' + \Phi(u) + \Psi(u) = 0, \qquad
    \Phi(u) = \nabla_g S(u)^T \calN[S(u)] + \nabla_g F(u).
\end{equation}
It appears that the mechanism behind these perturbations does not rely on the specifics of \eqref{eq:TWEnl}.
We are therefore of the conviction that such perturbations may find a wide range of applications in other Floer-type theories.

These perturbations are defined in Section \ref{sec:TWN_setup}.
In Section \ref{subsec:sojourn} we derive estimates which allow us to control the size of $\Psi(u)$ in terms of the kinetic energy of $u$.
These are essential to ensure uniform bounds on the kinetic energy and the gradient-like dichotomy of \eqref{eq:TWEnl} over $E$, 
derived in Section \ref{sec:TWN_heteroclinics}, is retained by \eqref{eq:TWN} for sufficiently small $\Psi$,
a condition we refer to as $E$-tameness of $\Psi$.
The $E$-tameness of small perturbations $\Psi$ is proven in Section \ref{sec:TWN_generic}, along with the the genericity of $E$-tranverse perturbations.
This then accumulates into the following result, which paraphrases Theorem \ref{thm:TWN_Psi_generic_transverse} and Theorem \ref{thm:TWN_generic_Floer_homology}.
\begin{mythm*}
  There exists a generic class of $E$-transverse perturbations $\Psi$ of \eqref{eq:TWEnl}.
  For each such perturbation $\Psi$, the boundary operator $\partial_n$ satisfies the fundamental relation
  \[
  \partial_n \partial_{n+1} = 0.
  \]
  Consequently, the Conley--Floer homology groups for \eqref{eq:TWN}
\[
\HF_n(E,\Phi,\Psi;\Z_2) := \frac{ \ker{ \partial_n } }{ \img{ \partial_{n+1} } }
\]
are well-defined.
\end{mythm*}

The homology satisfies various continuation principles, as is typical for Morse--Floer theories.
In particular, up to isomorphism, the homology is independent of the perturbation $\Psi$.
Furthermore, the homology is invariant under homotopies of $\Phi$ which are stable with respect to $E$.
These results are collected in Section \ref{sec:TWN_homology}.

\paragraph{Multivalued dynamics and Morse isomorphism.}
In Section \ref{sec:TWN_Morse_isomorphism} we interpret $\Phi$ as a multivalued vector field on $\R^d$, 
and define isolating blocks for these multivalued vector fields.
Roughly speaking, these are compact manifolds with corners $B \subset \R^d$ so that the multivalued vector field $\Phi$ is transverse to $\bdy B$.
The boundary then decomposes as $\bdy B = \bdy B_+ \cup \bdy B_- \cup \bdy B_0$, where $\Phi$ is pointing inwards (respectively, outwards) on $\bdy B_-$ (respectively, on $\bdy B_+$),
and curves away from $B$ at the corners $\bdy B_0$.
Such an isolating block induces an isolating trajectory neighbourhood $E(B)$, which consists of all maps $u : \R \to \inter{B}$.

In the absence of nonlocal coupling, $\Phi$ becomes the single-valued vector field $\nabla_g F$, for which isolating blocks $B$ are well known, and in applications, easily constructed using Lyapunov functions.
When, in the transverse direction on $\bdy B$, the nonlocal coupling term $\nabla S(\cdot)^T \calN[ S(\cdot) ]$ is (when interpreted as a multivalued vector field)
uniformly smaller than $\nabla_g F$, the set $B$ is an isolating block for the nonlocal equation, as well.
We can then formulate a Morse isomorphism (paraphrasing Theorem \ref{thm:TWN_Morse_isomorphism}).
\begin{mythm*}
  Suppose that, in the transverse direction on $\bdy B$, the nonlocal coupling term $\nabla S(\cdot)^T \calN[ S(\cdot) ]$ is 
  (when interpreted as a multivalued vector field) uniformly smaller than $\nabla_g F$.
  Then
  \[
  \HF_*(E(B),\Phi;\Z_2) \iso \homology_*( B , \bdy B_- ; \Z_2 ),
  \]
  where $\homology_*$ denotes the relative singular homology with $\Z_2$ coefficients.
\end{mythm*}
This result indicates the Conley--Floer homology may be interpreted an extension of the homological Conley index to multivalued vector fields.

In applications the set of bounded solutions $\Sigma$ is often compact.
In that case a global Conley--Floer homology may be defined as the Conley--Floer homology associated with arbitrary large isolating blocks $B$.
A similar Morse isomorphism holds for the global Conley--Floer homology, but this time the nonlocal coupling term only has to be ``small at $\infty$'' when compared to $\nabla_g F$;
see Theorem \ref{thm:TWN_global_Morse_isomorphism}.

\paragraph{Forcing theorems.}
By construction of the Conley--Floer homology, when (part of) the chain groups are known, we can exploit the relation $\partial_n \partial_{n+1} = 0$
to prove existence and multiplicity of heteroclinic solutions, firstly, 
in the perturbed equation \eqref{eq:TWN}, and secondly, using compactness arguments, also in the unperturbed equation \eqref{eq:TWEnl}.
We refer to results of this type as forcing theorems.
The following Morse inequality is an example of such a result (paraphrasing Theorem \ref{thm:TWN_Morse_inequality}).
\begin{mythm*}
Let $E$ be an isolating trajectory neighbourhood for \eqref{eq:TWEnl}.
Then
  \begin{multline*}
      2 \cdot \# \big\{ \text{ heteroclinic solutions of \eqref{eq:TWEnl} in $E$, modulo translation } \big\} \\
      \geq \# \big\{ \text{ nondegenerate constant solutions of \eqref{eq:TWEnl} in $E$ } \big\} - \rank \HF_*(E,\Phi;\Z_2).
  \end{multline*}
Each of these bounded solutions $u$ is a heteroclinic, accumulating onto constant solutions of \eqref{eq:TWEnl} as $x\to \pm\infty$.
Furthermore, at least one of the limits of $u(x)$ as $x \to - \infty$ or $x \to \infty$ is nondegenerate.
\end{mythm*}

\begin{figure}[h]
  \centering
  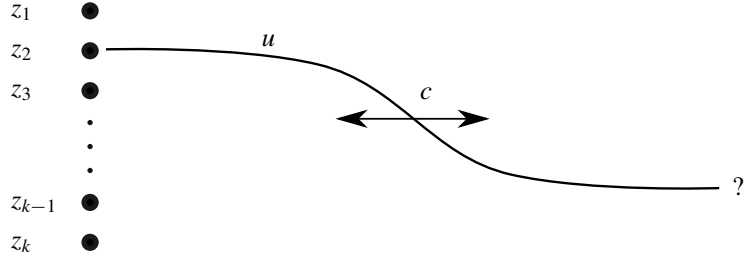
  \caption{For each wave speed $c \neq 0$, the existence of $k$ nondegenerate constant solutions $z_1, \dots, z_k$, where $k > \rank \HF_*(E,\Phi;\Z_2)$, 
    forces existence of a travelling front $u$, modelling the propagation (with linear wave speed $c$) of a spatially homogeneous state $z_+$ into an unstable spatially homogeneous state $z_-$.
    At least one of these states ($z_+$ or $z_-$) is one of the specified states $z_1, \dots, z_k$,
  whilst the other state ($z_-$ or $z_+$, respectively) may be degenerate.}
  \label{fig:TWN_forcing}
\end{figure}

Finally, in Section \ref{subsec:TWN_wave_application} we return to Examples \ref{ex:TWN_NFE} and \ref{ex:TWN_nonlocal_diffusion}.
We compute the Conley--Floer homology explicitly for specific nonlinearities $S$ and $F$.
For the neural field equations consisting of $d$ layers and symmetric interaction kernel, as in Example \ref{ex:TWN_NFE}, 
the Conley--Floer homology computes the reduced homology of a $d$-sphere, hence it is always of rank $1$.
For the nonlocal reaction-diffusion equations in Example \ref{ex:TWN_nonlocal_diffusion} with superlinear growth of $\|\nabla F\|$, 
the rank of the homology can be any natural number, and is determined by the asymptotics of the potential $F$.
This allows us to derive existence and multiplicity of travelling front solutions.
A typical application is depicted in Figure \ref{fig:TWN_forcing}.

\paragraph{Outline of the paper.}
In Section \ref{sec:TWN_setup} we specify restrictions upon $\Phi$ under which our theory works.
The perturbations $\Psi$ are defined in this section, as well.
Section \ref{sec:TWN_functional} collects various functional analytic results about \eqref{eq:TWN}.
In Section \ref{sec:TWN_compactness} we derive properties of solutions to \eqref{eq:TWN} which do not rely on the gradient-like structure.
These include compactness, asymptotics, and kinematic estimates on $\Psi(u)$.
Isolating trajectory neighbourhoods are also defined in this section.
A Lyapunov function for \eqref{eq:TWN} is derived in Section \ref{sec:TWN_heteroclinics}, under a technical assumption on $\Psi$, which we refer to as $E$-tameness.
Using this Lyapunov function we obtain uniform bounds on the kinetic energy and a gradient-like dichotomy for \eqref{eq:TWN}.
We then introduce the moduli spaces of heteroclinic solutions, under the hypothesis that the perturbation $\Psi$ is $E$-tame and $E$-transverse.
The genericity of $E$-tame and $E$-transverse perturbations $\Psi$ is verified in Section \ref{sec:TWN_generic}.
Section \ref{sec:TWN_homology} contains the definition of the Conley--Floer homology as well as various continuation results.
Isolating blocks for multivalued vector fields and a Morse isomorphism are discussed in Section \ref{sec:TWN_Morse_isomorphism}.
Finally, Section \ref{sec:TWN_application} contains the derivation of a Morse inequality and various examples of applications.


\section{Problem setup}
\label{sec:TWN_setup}

The goal of this section is to provide a concise description of the restrictions we put upon the various parameters occurring in \eqref{eq:TWN}.
More precisely, we specify conditions on the nonlocal operator $\calN$ and the local nonlinear maps $S$, $F$, and the Riemannian metric $g$, all occurring in the definition of $\Phi$.
We also construct the spaces of perturbations $\Psi$.

\subsection{Description of the class of equations}

Throughout this paper, we denote by $| \cdot |$ the Euclidean norm on $\F^d$, where $\F = \R$ or $\F = \C$.
We write $\Mat_{n_1 \times n_2}(\F)$ for the space of $n_1 \times n_2$ matrices with coefficients in the field $\F$,
and denote by $\| \cdot \|$ the operator norm on $\Mat_{n_1 \times n_2}(\F)$.
Norms on function spaces, consisting of functions taking values on $\F^d$ or $\Mat_{n_1 \times n_2}(\F)$, will always be defined in terms of these norms.

We now specify restrictions upon the function
\[
\Phi(u) := \nabla_g S(u)^T \calN[ S(u) ] + \nabla_g F(u)
\]
appearing in \eqref{eq:TWN}.
Typically, in the applications we have in mind, the linearisation $DS(u)$ is positive definite.
Consequently, we could choose the Riemannian metric $g_u(v_1,v_2) = \big\langle v_1 , DS(u) v_2 \big\rangle$ for $v_1, v_2 \in T_u \R^d = \R^d$,
so that equation \eqref{eq:TWN} reduces to
\[
u' + \calN[S(u)] + \nabla_g F(u) = 0.
\]
We allow for other types maps $S$ and Riemannian metrics $g$ as well, since it gives more flexibility once we wish to compute the homology groups.

The map $S$ is assumed to map from $\R^d$ into $\R^D$, for some ``intermediate dimension'' $D$.
The nonlocal operator $\calN$ acts on maps $\R \to \R^D$.
The local term $F$ is a potential on $\R^d$, and $g$ is a Riemannian metric on $\R^d$.
The $g$-gradient of $S$ is defined element-wise, thus, the term $\nabla_g S(u)$ is a $D \times d$ matrix whose $i^{\text{th}}$ row is given by $\nabla_g S_i(u)^T$, where $S = (S_1,\dots,S_D)$.

 \begin{myhyp}[$\calN$]
\namedlabel{hyp:calN}{$\calN$}
   The nonlocal operator $\calN$ decomposes into continuous and discrete convolution operators
   \[
   \calN[v](x) = \int_\R N(y) v(x-y) \d y + \sum_{j\in\Z} N_j v(x - \lambda_j).
   \]
   Here $N : \R \to \Mat_{D\times D}(\R)$, $N_j \in \Mat_{D\times D}(\R)$ and $\lambda_j \in \R$ are such that
   \begin{enumerate}
   \item $N(x) = N(-x)$, and $N(x) = N(x)^T$,
   \item $N_j = N_{-j}$ and $\lambda_j = - \lambda_{-j}$, and $N_j = N_j^T$,
   \item there exists $\eta_0 > 0$ such that
     \[
     \int_\R e^{\eta_0 |x|} \| N(x) \| \d x < \infty, \qquad \sum_{j\in\Z} e^{\eta_0 |\lambda_j|} \| N_j \|  < \infty.
     \]
   \end{enumerate}
 \end{myhyp}

We now define the matrix $\widetilde{\calN}$ by
\[
\widetilde{\calN} := \int_\R N(x) \d x + \sum_{j\in\Z} N_j.
\]
Then $\calN[z] = \widetilde{\calN} z$ for any constant $z \in \R^d$.

\begin{myhyp}[$S$, $g$, $F$]
\namedlabel{hyp:SF}{$S$, $g$, $F$}
  We impose the following restrictions upon $S$, $g$ and $F$.
  \begin{enumerate}
\item The function $S : \R^d \to \R^D$ is $C^4$.
\item The Riemannian metric $g$ on $\R^d$ is $C^3$ and trivialisable. 
More precisely, there exists a $C^3$ smooth map $G : \R^d \to \Mat_{d\times d}(\R)$,
such that for each $u \in \R^d$ the matrix $G(u)$ is positive definite, 
and $g$ is defined by
\[
g_u(v_1,v_2) := \big\langle v_1 , G(u) v_2 \big\rangle, \qquad v_1, v_2 \in T_u \R^d = \R^d.
\]
\item The function $F : \R^d \to \R$ is $C^4$.
\end{enumerate}
\end{myhyp}

With the choice of the Riemannian metric $g$ comes a natural notion of kinetic energy.
\begin{mydef}[Kinetic energy]
  We define the kinetic energy of $u \in W^{1,2}_{\text{loc}}(\R,\R^d)$ by
\[
\calE_{\text{kin}}(u) := \frac 1 2 \int_\R g_{u(x)}\big(u'(x),u'(x)\big) \d x,
\]
which may be infinite.
\end{mydef}

Finally, we remark here that the class of nonlinearities $\Phi$ is chosen in such a way that, formally, $\Phi$ is the gradient of a functional $H$
in the $L^2$ inner product induced by $g$, where $H$ is given by
\[
H(u) = \frac 1 2 \int_\R S(u(x)) \cdot \calN[S(u)](x) \d x + \int_\R F(u(x)) \d x,
\]
but we caution that \eqref{eq:TWN} should not be thought of as a formal gradient flow of $H$,
since the derivative appearing in \eqref{eq:TWN} acts on the spatial variable $x$.

\subsection{Constant solutions and the Morse function $h$}
\label{subsec:Morse_function_h}
Define $h : \R^d \to \R$ by
\[
h(z) := \frac 1 2 S(z) \cdot \widetilde{\calN} S(z) + F(z), \qquad \widetilde{\calN} = \int_\R N(x) \d x + \sum_{j\in\Z} N_j.
\]
Denote by $\crit(h) := \set{ z \in \R^d }{ D h(z) = 0 }$ the set of critical points of $h$.
The perturbations $\Psi$, which we define in the next subsection, are assumed to be uniformly small, and vanish in a neighbourhood of $\crit(h)$.
Hence it follows that critical points $z \in \crit(h)$ correspond to the constant solutions of \eqref{eq:TWN}.

\begin{mydef}[Hyperbolicity]
\label{def:TWN_hyperbolic}
  A critical point $z \in \crit(h)$ is said to be hyperbolic if the Hessian matrix $\Hess_g h(z)$ is nondegenerate.
  The map $h$ is said to be a Morse function if all critical points are hyperbolic.
\end{mydef}
By the Thom transversality theorem any $h$ can be made into a Morse function by arbitrarily small perturbations of $F$.
Throughout sections 2--7 we will assume the following hypothesis.
\begin{myhyp}
  We assume the map $h$ is a Morse function.
\end{myhyp}
Towards the end of Section \ref{sec:TWN_homology}, after constructing the Conley--Floer homology, we will drop this hypothesis on $h$.

\begin{mydef}[Morse index]
  Given a hyperbolic critical point $z \in \crit(h)$ of $h$, define its Morse index $m_h(z)$ as the
  dimension of the maximal subspace $E_- \subseteq \R^d$ on which the Hessian $\Hess_g h(z)$ is negative definite.
\end{mydef}

\subsection{Construction of perturbations $\Psi$}
\label{subsec:TWN_perturbations}
We begin with a construction of functions $\psi_{\theta;\xi} : C^0_{\text{loc}}(\R,\R^d) \to C^0_{\text{loc}}(\R,\R^d)$ which will form a basis for the space of perturbations $\bPsi_h$.

\paragraph{The spatial localiser $\sigma_{\ell,z}$.}
We first define a family of cutoff functions on $\R$.
Given $\ell \in \N$ and $\rho > 0$, let $\chi_{\ell,\rho} : \R \to \R$ be given by
\[
\chi_{\ell,\rho}(y) :=
\left\{
\begin{array}{l l}
  e^{- \ell^{-1} \left[ 1 - \left( 2 y / \rho - 3 \right)^2 \right]^{-1} } & \rho  < y < 2 \rho, \\
  0 & \text{otherwise},
\end{array}
\right.
\]
We also define $\chi_{\ell,\infty} = 0$.
Note that $\chi_{\ell,\rho}$ is smooth, and $\lim_{\ell\to\infty} \chi_{\ell,\rho}(y) = \1_{(\rho,2\rho)}(y)$.

Now define $\rho : \crit(h) \to (0,\infty]$ by
\begin{equation}
  \label{eq:TWN_Psi_rho_z}
  \rho(z) := \frac 1 3 \inf_{ \substack{ z' \in\crit(h) \\ z' \neq z } } | z' - z |, \qquad z \in \crit(h),
\end{equation}
with the convention that the infimum over an empty set equals $\infty$.
Since $h$ is assumed to be a Morse function, the critical points $\crit(h)$ form a discrete subset of $\R^d$, hence $\rho(z) > 0$.
For given $\ell \in \N$ and $z \in \crit(h)$, define $\sigma_{\ell,z} : \R^d \to \R$ by
\[
\sigma_{\ell,z}(u) := \chi_{\ell , \rho(z)}( | u - z | ).
\]
This construction is made such that that the support of $\sigma_{\ell,z}(u)$ is bounded away from $\crit(h)$ by a distance $\rho(z)$.
It will be useful to identify $\sigma_{\ell,z}$ with the induced smooth map
\[
\sigma_{\ell,z} : C^0_{\text{loc}}(\R,\R^d) \to C^0_{\text{loc}}(\R,\R), \qquad u(x) \mapsto \sigma_{\ell,z}(u(x)).
\]

\paragraph{The basic perturbations $\psi_{\theta;\xi}$.}
Let
\[
\Theta_h := \set{ ( \ell , z_1 , z_2 ) \in \N \times \crit(h) \times \crit(h) }{ z_1 \neq z_2 }.
\]
For a given $\xi =(\xi_1,\xi_2) \in \R^2$ and $\theta = (\ell,z_1,z_2) \in \Theta_h$, we now define a basic perturbation $\psi_{\theta;\xi} : C^0_{\text{loc}}(\R,\R^d) \to C^0_{\text{loc}}(\R,\R)$ by
\[
\rst{ \psi_{\theta;\xi}(u) }{ x } := \rst{ \sigma_{\ell,z_1}(u) }{ x - \xi_1 } \rst{ \sigma_{\ell,z_2}(u) }{ x - \xi_2 }.
\]
It will also be convenient to define $\psi_\theta : C^0_{\text{loc}}(\R,\R^d) \to C^0_{\text{loc}}(\R^2,\R)$ by
\[
\rst{ \psi_\theta(u) }{ \xi } := \rst{ \sigma_{\ell,z_1}(u) }{ \xi_1 } \rst{ \sigma_{\ell,z_2}(u) }{ \xi_2 },
\]
so that
\[
\rst{ \psi_{\theta;\xi}(u) }{ x } = \rst{ \psi_\theta(u) }{ (x,x) - \xi }.
\]

\begin{figure}
  \centering
  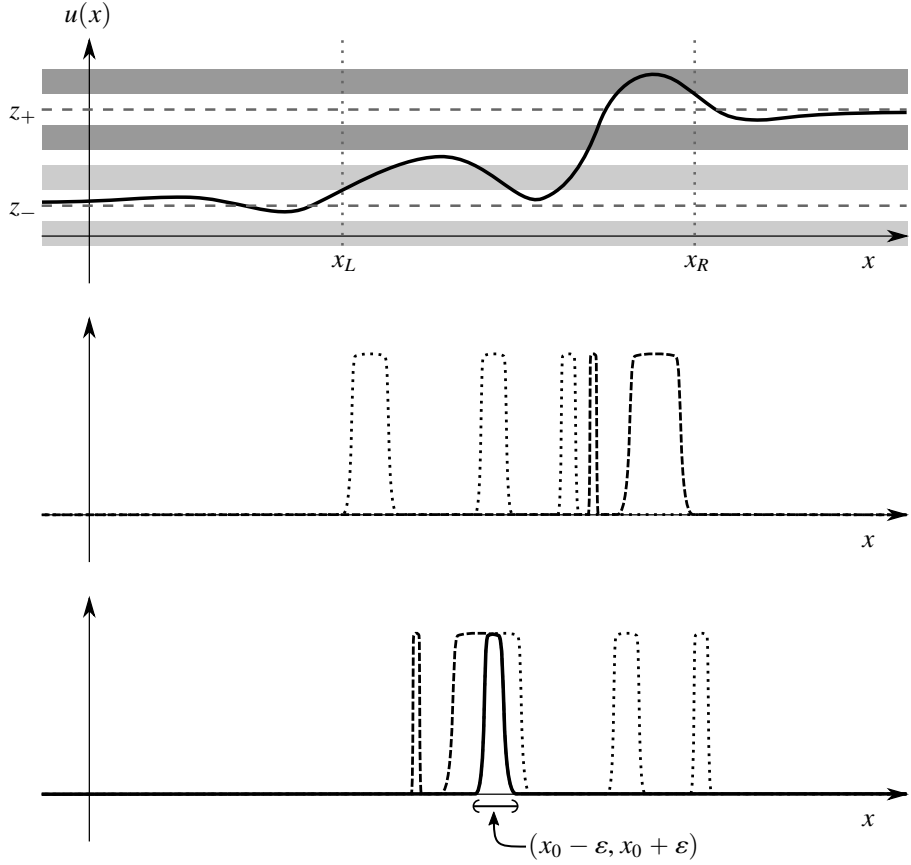
  \caption{Construction of the basic perturbation $\psi_{\theta;\xi}$.
Top: Graph of a curve $u : \R \to \R^d$ with $u(x) \to z_\pm$ as $x\to\pm\infty$. Lightly shaded area indicates support of $\sigma_{\ell,z_-}$, 
and darkly shaded area indicates support of $\sigma_{\ell,z_+}$.
Middle: The spatial localisers $\sigma_{\ell,z_-}$ (dotted line) and $\sigma_{\ell,z_+}$ (dashed line).
Bottom: The construction of the basic perturbation $\psi_{\theta,(\xi_1(x_0,\epsilon),\xi_2(x_0,\epsilon))}$ (solid line), with $\theta = (\ell,z_-,z_+)$.
Dotted line depicts $\sigma_{\ell,z_-}(x-\xi_1(x_0,\epsilon))$.
Dashed line depicts $\sigma_{\ell,z_+}(x-\xi_2(x_0,\epsilon))$.}
  \label{fig:TWN_perturbation_intuition}
\end{figure}

\paragraph{Some intuition.}
At this stage we can motivate the usefulness of this class of perturbations.
Figure \ref{fig:TWN_perturbation_intuition} gives a graphical depiction of the situation.
Suppose $u \in C^1_{\text{loc}}(\R,\R^d)$ is a solution of \eqref{eq:TWN} such that $u(x) \to z_\pm$ as $x\to\pm\infty$, where $z_\pm \in \crit(h)$.
Then the curve $x \mapsto u(x)$ must tranverse the annulus with inner and outer radii $\rho(z_-)$ and $2 \rho(z_-)$ centred at $z_-$, 
so that as $\ell \to \infty$, the spatial localiser $\sigma_{\ell,z_-}(u)$ converges to an indicator function with compact and nonempty support.
Likewise, the spatial localiser $\sigma_{\ell,z_+}(u)$ converges to an indicator function with compact and nonempty support.
Let $x_L := \essinf{\big( \supp\big( \sigma_{\ell,z_-}(u) \big) \big)}$ and $x_R := \esssup{\big( \supp\big( \sigma_{\ell,z_+}(u) \big) \big)}$.
Then, given $x_0 \in \R$ and sufficiently small $\epsilon > 0$, let
\[
\xi_1(x_0,\epsilon) := x_0 - x_L - \epsilon, \qquad \xi_2(x_0,\epsilon) := x_0 - x_R + \epsilon.
\]
It follows that
\[
 \lim_{\ell\to\infty} \psi_{ ( \ell , z_- , z_+ ) ; (\xi_1(x_0,\epsilon),\xi_2(x_0,\epsilon))}(u) = \1_{(x_0-\epsilon,x_0+\epsilon)},
\]
so that
\[
\frac{1}{2\epsilon} \psi_{( \ell , z_- , z_+ ); (\xi_1(x_0,\epsilon),\xi_2(x_0,\epsilon))}(u) \rightharpoonup \delta_{x_0}
\]
as $\ell \to \infty$ and $\epsilon \to 0$.
Therefore, if $u$ is a solution of \eqref{eq:TWN} in the unperturbed case $\Psi = 0$,
which is singular in the sense that the set of all bounded solution of \eqref{eq:TWN} lacks regularity (with respect to some appropriately chosen topology) in the vicinity of $u$,
we can smoothen the singularity by adding a perturbation $\Psi$ which approximates a delta peak near $u$.
A precise argument will be postponed to Section \ref{sec:TWN_generic}.

\paragraph{Spaces of perturbations $\Psi$.}
We now define a module of perturbations $\Psi$ by taking $\R^d$-linear combinations of the basic perturbations $\psi_{\theta;\xi}$.
If $\Theta_h = \emptyset$, equation \eqref{eq:TWN} has at most $1$ constant solution.
In that case we will not require any perturbation, and set $\Psi = 0$.
Henceforth assume $\Theta_h \neq \emptyset$.
Given a map $\alpha : \Theta_h \times \R^2 \to \R^d$, we define
\begin{equation}
  \label{eq:Psi_perturbation}
  \Psi(u) := \sum_{\theta \in \Theta_h} \iint_{\R^2} \alpha(\theta,\xi) \psi_{\theta;\xi}(u) \d \xi
\end{equation}
and
\[
\| \Psi \|_{\bPsi_h} := \sum_{\theta = (\ell,z_1,z_2) \in \Theta_h} \iint_{\R^d} \ell^3 e^{\eta_0 | \xi |} \bigg( | \alpha(\theta,\xi) | + \| D_\xi \alpha(\theta,\xi) \| \bigg) \d \xi,
\]
where $\eta_0$ is as in Hypothesis \ref{hyp:calN}.
We note that
\[
\sup_{u \in C^0_{\text{loc}}(\R,\R^d)} \| \psi_{\theta;\xi}(u) \|_{L^\infty(\R,\R)} \leq 1.
\]
Hence, whenever $\| \Psi \|_{\bPsi_h} < \infty$, the sum and integrals in \eqref{eq:Psi_perturbation} converge in $C^1_{\text{loc}}(\R,\R^d)$,
in the sense of Bochner integration, with rates of convergence uniform in $u \in C^0_{\text{loc}}(\R,\R^d)$.

We now define a vector space $\bPsi_h$ of perturbations, consisting of all maps $\Psi$ of the form \eqref{eq:Psi_perturbation} with $\| \Psi \|_{\bPsi_h} < \infty$.
Then $( \bPsi_h , \|\cdot\|_{\bPsi_h} )$ is a separable Banach space.

It will also be convenient to introduce the notation $\alpha_\theta(\xi) := \alpha(\theta,\xi)$, 
and then define $\Psi_\theta : C^0_{\text{loc}}(\R,\R^d) \to C^1_{\text{loc}}(\R,\R^d)$ by
\[
\rst{ \Psi_\theta(u) }{ x } := \rst{ \alpha_\theta * \psi_\theta(u) }{(x,x)} = \iint_{\R^2} \alpha(\theta,(x,x)-\xi) \rst{ \psi_{\theta}(u) }{\xi} \d \xi.
\]
Thus, elements $\Psi \in \bPsi_h$ are of the form
\begin{equation}
  \label{eq:Psi_ell_perturbation}
  \Psi(u) = \sum_{\theta \in \Theta_h} \Psi_\theta(u).
\end{equation}

\paragraph{$R$-admissible perturbations.}
We now return to the relation between $\crit(h)$ and constant solutions of \eqref{eq:TWN}.
By construction $\Psi$ vanishes in a neighbourhood of $\crit(h)$, so that any critical point of $h$ is a constant solution of \eqref{eq:TWN}.
The reverse inclusion requires further restrictions upon $\Psi$, as given in the following hypothesis.
\begin{mydef}[$R$-admissible perturbations]
  Given $R \geq 0$, a perturbation $\Psi \in \bPsi_h$ is said to be $R$-admissible if
  the constant solutions $z$ of \eqref{eq:TWN} with $|z| \leq R$ are in one-to-one correspondence with the critical points of $h$ inside the ball of radius $R$.
\end{mydef}

The following lemma demonstrates that $R$-admissible perturbations form open subsets of the spaces of perturbations.
In particular, as $\Psi = 0$ is obviously $R$-admissible, this shows that $R$-admissible perturbations exist, for any given $R \geq 0$.
\begin{mylemma}
\label{lemma:R-admissible_open}
  Let $R \geq 0$ and suppose $\Psi_0 \in \bPsi_h$ is an $R$-admissible perturbation.
Then there exists $\epsilon > 0$ such that any $\Psi \in \bPsi_h$ with $\| \Psi - \Psi_0 \|_{\bPsi_h} < \epsilon$ is $R$-admissible.
\end{mylemma}
\begin{myproof}
  If $z \in \crit(h)$, then $\Phi(z) + \Psi(z) = \Phi(z) = \nabla_g h(z) = 0$.
  On the other hand, suppose $\Phi(z) + \Psi(z) = 0$ and $|z| \leq R$.
  Let $\rho : \crit(h) \to (0,\infty]$ be as defined in \eqref{eq:TWN_Psi_rho_z}, and set
  \[
  \rho_R := \inf_{ \substack{ z_0 \in \crit(h) \\ | z_0 | \leq R } } \rho(z_0).
  \]
  The definition of $\Psi$ ensures that $\Psi(z)=0$ if $z$ is at most a distance $\rho_R$ away from $\crit(h)$.
  We will now specify a choice of $\epsilon > 0$ for which this is the case.
  Denote by $B_r(v_0) = \set{ v \in \R^d }{ |v-v_0| \leq r }$ the closed ball of radius $r$ around $v_0$ in $\R^d$.
  Now set
  \[
  \epsilon := \inf\set{ |\Phi(v) + \Psi_0(v)| }{ v \in B_R(0) \setmin \bigcup_{z\in\crit(h)} B_{\rho_R}(z) },
  \]
  and note that since $\Psi_0$ is $R$-admissible it follows that $\epsilon > 0$.
  If $\| \Psi - \Psi_0 \|_{\bPsi_h} < \epsilon$, we have 
  \begin{align*}
    | \Phi(z) + \Psi_0(z) | &\leq | \Phi(z) + \Psi(z) | + | \Psi(z) - \Psi_0(z) | = | \Psi(z) - \Psi_0(z) | \\
    &\leq \| \Psi - \Psi_0 \|_{\bPsi_h} < \epsilon,
  \end{align*}
  so that per definition of $\epsilon$ the point $z$ is at most a distance $\rho_R$ away from $\crit(h)$.
  The definition of the perturbation $\Psi$ then ensures that $\Psi(z) = 0$, from which it follows that $\nabla_g h(z) = \Phi(z) = \Phi(z) + \Psi(z) = 0$.
\end{myproof}



\section{Functional properties}

\label{sec:TWN_functional}

Note that solutions of \eqref{eq:TWN} are the zeroes of the map 
\[
\partial_x + \Phi(\cdot) + \Psi(\cdot),
\]
where
\[
\Phi(u) = \nabla_g S(u)^T \calN[S(u)] + \nabla_g F(u).
\]
In this section we will study various functional analytic properties of this map.


\subsection{Continuity in compact-open topology}

We outline some properties of the map
\[
\partial_x + \Phi(\cdot) + \Psi(\cdot) : C^1_{\text{loc}}(\R,\R^d) \to C^0_{\text{loc}}(\R,\R^d).
\]
We start out with a general continuity result for convolution operators.
\begin{mylemma}
\label{lemma:convolution_weak_continuity}
  Let $m , n \in \N$, suppose $K \in L^1(\R^m,\Mat_{n\times n}(\R))$, $(K_j)_j \in \ell^1(\Z,\Mat_{n\times n}(\R))$, and $(\xi_j)_{j\in \Z} \subset \R^m$, and for $v \in L^\infty(\R^m,\R^n)$ define
  \[
  \calK[v](x) := \int_{\R^m} K(y) v(x-y) \d y + \sum_{j\in\Z} K_j v(x-\xi_j).
  \]
  Let $k\in \N\cup\{0\}$, 
  and suppose $(v_j)_j \subset C^k_{\text{loc}}(\R^m,\R^n)$ is a sequence converging in $C^k_{\text{loc}}(\R^m,\R^n)$ to $v_\infty$.
  Assume further that $\sup_j \| v_j \|_{W^{k,\infty}(\R^m,\R^n)} < \infty$.
  Then it follows that $\calK[v_j] \to \calK[v_\infty]$ in $C^k_{\text{loc}}(\R^m,\R^n)$ as $j\to\infty$.
\end{mylemma}
\begin{myproof}
  We prove the result for $k=0$ and $n=1$; the general result follows from this, first for $k=0$ and arbitrary $n$ by applying this result elementwise to the matrix-vector multiplications,
  and then for arbitrary $k$ by differentiation under the integral sign and repeatedly applying the convergence result for $C^0_{\text{loc}}(\R^m,\R^n)$.
  The result for $k=0$ and $n=1$ follows from the following slightly more general result.
  Suppose $\mu$ is a signed Radon measure of bounded variation on $\R^m$, and for $w \in L^\infty(\R^m,\R)$ define
  \[
  (w * \mu)(x) := \int_{\R^m} w(x-y) \d \mu(y).
  \]
  We will show that if $(w_j)_j \subset C^0_{\text{loc}}(\R^m,\R)$ is a sequence converging in $C^0_{\text{loc}}(\R^m,\R)$ to $w_\infty$, and $\sup_j \| w_j \|_{L^\infty(\R^m,\R)} < \infty$,
  then $w_j * \mu \to w_\infty * \mu$ in $C^0_{\text{loc}}(\R^m,\R)$.

  Since $\mu$ is a signed Radon measure of bounded variation, given $\delta > 0$ we may find a compact set $F \subset \R^m$ such that
\[
| \mu |(\R^m \setmin F) < \delta,
\]
hence, as $\sup_j \| w_j \|_{L^\infty(\R^m,\R)} < \infty$, we find that given $\epsilon > 0$ there exists compact $F \subset \R^m$ such that for all $j \in \N \cup \{\infty\}$
\[
 \sup_{x\in\R^m} \left| \int_{\R^m \setmin F} w_j(x-y) \d \mu(y) \right| < \epsilon / 4.
\]
Now let $E \subset \R^m$ be another compact subset, and set
\[
E - F := \set{ x - y }{ x \in E,\; y \in F } \subset \R^m.
\]
Then $E-F$ is compact, hence $\| w_\infty - w_j \|_{L^\infty(E-F,\R)} \to 0$ as $j \to \infty$.
Now
\[
\sup_{x\in E} \left| \int_{F} \big( w_\infty(x-y) - w_j(x-y) \big) \d \mu(y) \right| \leq \| w_\infty - w_j \|_{L^\infty(E-F,\R)} |\mu|(F),
\]
hence there exists $j_0 \in \N$ so that
\[
\sup_{x\in E} \left| \int_{F} w_\infty(x-y) \d \mu(y) - \int_F w_j(x-y) \d \mu(y) \right| < \epsilon / 2, \qquad j \geq j_0.
\]
Combining these estimates, we conclude that
\[
\| w_\infty * \mu - w_j * \mu \|_{L^\infty(E,\R)} < \epsilon \qquad \text{for all} \quad j \geq j_0.
\]
Since $E$ and $\epsilon$ where chosen arbitrary, we conclude that $( w_j * \mu )_j$ converges in $C^0_{\text{loc}}(\R^m,\R)$ towards $w_\infty * \mu$.
\end{myproof}

\paragraph{The term $\Phi$.}
We first focus on the mapping properties of the map
\[
\Phi : C^k_{\text{loc}}(\R,\R^d) \to C^k_{\text{loc}}(\R,\R^d),
\]
for given $k \in \{ 0 , \dots , 3\}$.

\begin{mylemma}
  \label{lemma:TWN_Phi_continuity}
Let $k \in \{ 0 , \dots , 3\}$ and let $(u_n)_n \subset C^k_{\text{loc}}(\R,\R^d)$ be a sequence converging $C^k_{\text{loc}}(\R,\R^d)$ to $u_\infty$.
Suppose $\sup_n \| u_n \|_{W^{k,\infty}(\R,\R^d)} < \infty$.
Then $\Phi(u_n) \to \Phi(u_\infty)$ in $C^k_{\text{loc}}(\R,\R^d)$.
\end{mylemma}
\begin{myproof}
  The continuity of the nonlocal term $\calN$ follows from Lemma \ref{lemma:convolution_weak_continuity}.
  Since $\Phi$ consists of the composition of pointwise evaluations by $C^3$ smooth functions and linear nonlocal operations $\calN$,
  the result follows.
\end{myproof}

\paragraph{The term $\Psi$.}
We now collect some mapping properties of
\[
\Psi : C^k_{\text{loc}}(\R,\R^d) \to C^{k+1}_{\text{loc}}(\R,\R^d),
\]
for given $k \in \N \cup \{0\}$.

\begin{mylemma}
\label{lemma:Psi_weak_continuity}
Let $\Psi \in \bPsi_h$ and $k \in \{0,\dots,3\}$, and let $(u_n)_n \subset C^k_{\text{loc}}(\R,\R^d)$ be a sequence converging to $u_\infty$ in $C^k_{\text{loc}}(\R,\R^d)$.
Suppose furthermore that $\sup_n \| u_n \|_{W^{k,\infty}(\R,\R^d)} < \infty$.
Then $\Psi(u_n) \to \Psi(u_\infty)$ in $C^{k+1}_{\text{loc}}(\R,\R^d)$ as $n \to \infty$.
\end{mylemma}
\begin{myproof}
    By definition of $\Psi$, we have the decomposition
  \[
  \Psi(u) = \sum_{\theta \in \Theta_h} \Psi_\theta(u), \qquad \rst{\Psi_\theta(u)}{x} = \rst{ \alpha_\theta * \psi_\theta(u) }{(x,x)}.
  \]
  Differentiation under the sum and convolution term are possible by the  rapid decay rate of $\alpha_\theta$ and uniform convergence of the sum.
  Indeed, inspecting the definition of $\psi_\theta$, we see that given $R > 0$ there exists $C_R > 0$ so that for $k \in \{0,\dots,3\}$
  \[
  \sup_{  \| u \|_{W^{k,\infty}(\R,\R^d)} \leq R  } \| \partial_x^{k} \psi_\theta(u) \|_{L^\infty(\R^2,\R)} \leq C_R \ell^{k}, \qquad \text{where} \quad \theta = (\ell,z_1,z_2) \in \Theta_h.
  \]
  It is precisely the uniform estimate, combined with the imposed decay rate of $(\alpha_\theta)_\theta$, which justifies the differentiation under the sum and convolution,
  yielding
  \[
  \partial_x^{k+1} \rst{ \Psi(u) }{x} = \sum_{\theta\in \Theta_h} \rst{ (\nabla \alpha_\theta) * \partial_x^{k} \psi_\theta(u) }{(x,x)},
  \]
  where the series converges uniformly in $\| u \|_{W^{k,\infty}(\R,\R^d)} \leq R$.
  The conclusion of the lemma then follows by applying Lemma \ref{lemma:convolution_weak_continuity}.
\end{myproof}

\paragraph{Continuity of \eqref{eq:TWN}.}
Finally, we combine all these results to arrive at the following continuity theorem.
\begin{mythm}
\label{thm:TWN_weak_continuity}
Let $k \in \{ 0 , \dots , 3\}$.
  Suppose $(u_n)_n$ is convergent in $C^{k+1}_{\text{loc}}(\R,\R^d)$,
  and $\sup_n \| u_n \|_{W^{k,\infty}(\R,\R^d)} < \infty$,
 and $( \Psi_n )_n$ is convergent in $\bPsi_h$.
 Then the sequence
 \[
 u_n' + \Phi(u_n) + \Psi_n(u_n)
 \]
  is convergent in $C^k_{\text{loc}}(\R,\R^d)$.
\end{mythm}
\begin{myproof}
  Convergence of the sequence $( u_n' + \Phi(u_n) )_n$ follows immediately from Lemma \ref{lemma:TWN_Phi_continuity}.
  Let $\Psi_\infty$ be the limit point of $( \Psi_n )_n$.
  Then, for $x \in \R$,
  \begin{align*}
    | \Psi_\infty(u_\infty)(x) - \Psi_n(u_n)(x) | &\leq | \Psi_\infty(u_\infty)(x) - \Psi_\infty(u_n)(x) | + | \Psi_\infty(u_n)(x) - \Psi_n(u_n)(x) | \\
    &\leq | \Psi_\infty(u_\infty)(x) - \Psi_\infty(u_n)(x) | + \| \Psi_\infty - \Psi_n \|_{\bPsi_h}.
  \end{align*}
  By Lemma \ref{lemma:Psi_weak_continuity}, the first term tends to $0$ as $n \to \infty$, uniformly for $x$ in compact subsets of $\R$.
  The final term converges by assumption.
  This proves that $( \Psi_n(u_n) )_n$ converges in $C^0_{\text{loc}}(\R,\R^d)$ towards $\Psi_\infty(u_\infty)$.
  Convergence in $C^k_{\text{loc}}(\R,\R^d)$ is obtained similarly, exploiting the linear dependence on $\Psi$ of the term $\Psi(u)$.
\end{myproof}

\subsection{Uniform continuity}
We now study the restriction of the map $\partial_x + \Phi(\cdot) + \Psi(\cdot)$ to spaces of uniformly bounded functions.
Denote by $L^\infty(\R,\R^d)$ the Banach space of essentially bounded measurable functions $u : \R \to \R^d$.
For $k \in \N$, the Sobolev space $W^{k,\infty}(\R,\R^d)$ consists of measurable functions $u : \R \to \R^d$ for which the norm
\[
\| u \|_{W^{k,\infty}(\R,\R^d)} = \| u \|_{L^\infty(\R,\R^d)} + \| \partial_x u \|_{L^\infty(\R,\R^d)} + \cdots + \| \partial_x^k u \|_{L^\infty(\R,\R^d)}
\]
is finite.

\paragraph{The term $\Phi$.}
We first focus on the properties of the map $\Phi$.

\begin{mylemma}
\label{lemma:TWN_Phi_uniform_bound}
Given $R \geq 0$ there exists $C_R > 0$ such that for $k \in \{0,\dots,3\}$ we have the uniform bound
\[
\| \Phi(u) \|_{W^{k,\infty}(\R,\R^d)} \leq C_R \qquad \text{whenever} \quad \| u \|_{W^{k,\infty}(\R,\R^d)} \leq R.
\]
\end{mylemma}
\begin{myproof}
  By Young's convolution inequality, the nonlocal operator $\calN$ maps bounded subsets of $W^{k,\infty}(\R,\R^D)$ into bounded subsets of $W^{k,\infty}(\R,\R^D)$, 
  for any $k \in \N \cup \{0\}$.
  Since the map $\Phi$ consists of the composition of pointwise evaluations by $C^3$ smooth functions and linear nonlocal operations $\calN$, the result follows readily.
\end{myproof}

In particular, we obtain a map $\Phi : L^\infty(\R,\R^d) \to L^\infty(\R,\R^d)$.
We now inspect the regularity of this function.
\begin{mylemma}
\label{lemma:TWN_Phi_uniform_smooth}
  The map $\Phi : L^\infty(\R,\R^d) \to L^\infty(\R,\R^d)$ is $C^3$.
\end{mylemma}
\begin{myproof}
  The continuity of the nonlocal term $\calN$ follows from Young's convolution inequality.
  Since $\Phi$ consists of the composition of pointwise evaluations by $C^3$ smooth functions and linear nonlocal operations $\calN$,
  the result follows.
\end{myproof}

\paragraph{The term $\Psi$.}
We now inspect the map $\Psi$.

\begin{mylemma}
\label{lemma:TWN_Psi_uniform_bound}
  Given $R \geq 0$ there exists $C_R > 0$ such that for $k \in \{0,\dots,3\}$ and $\Psi \in \bPsi_h$ we have the uniform bound
  \[
  \| \Psi(u) \|_{W^{k+1,\infty}(\R,\R^d)} \leq C_R \| \Psi \|_{\bPsi_h}, \qquad \text{whenever} \quad \| u \|_{W^{k,\infty}(\R,\R^d)} \leq R.
  \]
\end{mylemma}
\begin{myproof}
  In the proof of Lemma \ref{lemma:Psi_weak_continuity} we derived that, for $k \in \{0,\dots,3\}$, we have the identity
  \[
  \partial_x^{k+1} \rst{ \Psi(u) }{x} = \sum_{\theta\in \Theta_h} \rst{ (\nabla \alpha_\theta) * \partial_x^{k} \psi_\theta(u) }{(x,x)}.
  \]
  By Young's convolution inequality it follows that
  \begin{align*}
    \sup_{  \| u \|_{W^{k,\infty}(\R,\R^d)} \leq R  } \| \partial_x^{k+1} \Psi(u)\|_{L^\infty(\R,\R^d)} &\leq C_R \sum_{\theta = (\ell,z_1,z_2) \in \Theta_h} \ell^k \| \nabla \alpha_\theta \|_{L^1(\R^2,\Mat_{d\times d}(\R))}  \\
    &\leq C_R \| \Psi \|_{\bPsi_h}.
  \end{align*}
\end{myproof}

We thus obtain the function $\Psi : L^\infty(\R,\R^d) \to L^\infty(\R,\R^d)$.
\begin{mylemma}
\label{lemma:TWN_Psi_uniform_smooth}
  The map $\Psi : L^\infty(\R,\R^d) \to L^\infty(\R,\R^d)$ is $C^3$.
\end{mylemma}
\begin{myproof}
  As in the proof of Lemma \ref{lemma:Psi_weak_continuity}, we make use of the representation
  \[
  \Psi(u) = \sum_{\theta \in \Theta_h} \Psi_\theta(u), \qquad \rst{\Psi_\theta(u)}{x} = \rst{ \alpha_\theta * \psi_\theta(u) }{(x,x)}.
  \]
We note that $\psi_\theta : L^\infty(\R,\R^d) \to L^\infty(\R^2,\R)$ is smooth, and straightforward computations show that, 
for any given $k \in \{0,\dots,3\}$ and a radius $R > 0$, there exists $C_R > 0$ such that
\[
\sup_{ \| u \|_{L^\infty(\R,\R^d)} \leq R }\| D^k \psi_\theta(u) \|_{\calL( L^\infty(\R,\R^d)^k , L^\infty(\R^2,\R))} \leq C_R l^k, \qquad \text{where} \quad \theta = (\ell,z_1,z_2).
\]
The uniformity of this estimate justifies the differentiation under the sum and convolution term appearing in $\Psi$. It follows that
\[
\rst{ D^k \Psi(u)[ v_1 , \dots , v_k ] }{ x } = \sum_{\theta\in \Theta_h} \rst{ \alpha_\theta * D^k \psi_\theta(u)[ v_1 , \dots , v_k ] }{(x,x)}
\]
for $v_1,\dots,v_k \in L^\infty(\R,\R^d)$, with the estimate
\begin{align*}
  \sup_{ \| u \|_{L^\infty(\R,\R^d)} \leq R }\| D^k \Psi(u) \|_{\calL( L^\infty(\R,\R^d)^k , L^\infty(\R,\R^d))} &\leq C_R \sum_{\theta = (\ell,z_1,z_2) \in \Theta_h} \ell^k \| \alpha_\theta \|_{L^1(\R^2,\R^d)} \\
  &\leq C_R \| \Psi \|_{\bPsi_h}.
\end{align*}
\end{myproof}

\paragraph{Continuity of \eqref{eq:TWN}.}
Finally, we combine all these results to arrive at the following continuity theorem.
\begin{mythm}
\label{thm:TWN_uniform_continuity}
The map
\[
W^{1,\infty}(\R,\R^d) \times \bPsi_h \to L^\infty(\R,\R^d), \qquad (u,\Psi) \mapsto u' + \Phi(u) + \Psi(u)
\]
is $C^3$.
\end{mythm}
\begin{myproof}
  The argument is identical to the proof of Theorem \ref{thm:TWN_weak_continuity}, 
  with the only modification being the uniform estimate holding over $\R$ instead of over compact intervals.
\end{myproof}

\subsection{Restriction to path space}

Given $z_- , z_+ \in \crit(h)$, i.e., hyperbolic constant solutions of \eqref{eq:TWN}, 
we define the path space $\calP(z_-,z_+)$ to be the $W^{1,2}$ Sobolev completion of smooth maps $u : \R \to \R^d$ with $u(x) \to z_\pm$ and $u'(x) \to 0$ as $x \to \pm\infty$.
Stated differently, the path space can be defined as the affine Hilbert manifold
\[
\calP(z_-,z_+) = z + W^{1,2}(\R,\R^d),
\]
where $z : \R \to \R^d$ is a chosen smooth path with $z(x) = z_-$ for $x \leq -1$ and $z(x) = z_+$ for $x \geq 1$.

We now study some of the properties of the map $\partial_x + \Phi(\cdot) + \Psi(\cdot)$ restricted to the path space $\calP(z_-,z_+)$.

\paragraph{The map $\Phi$.}
We consider the restriction of the map $\Phi : C^0_{\text{loc}}(\R,\R^d) \to C^0_{\text{loc}}(\R,\R^d)$ to the path space $\calP(z_-,z_+)$.

\begin{mylemma}
\label{lemma:TWN_Phi_path_space}
  The function $\Phi$ restricts to a map $\Phi : \calP(z_-,z_+) \to L^2(\R,\R^d)$.
\end{mylemma}
\begin{myproof}
  Recall that
  \[
  \nabla_g h(u) = \nabla_g S(u)^T \widetilde{\calN} S(u) + \nabla_g F(u), \qquad \widetilde{\calN} = \int_\R N(x) \d x + \sum_{j\in\Z} N_j.
  \]
  Let $z : \R \to \R^d$ be a smooth path with $z(x) = z_-$ for $x \leq -1$ and $z(x) = z_+$ for $x \geq 1$.
  Then
  \begin{align*}
    \nabla_g S(u)^T \bigg( \calN[ S(u) ] - \widetilde{\calN} S(u) \bigg) &= \nabla_g S(u)^T \big( \calN - \widetilde{\calN} \big)\big[ S(u) - S(z) \big] + \nabla_g S(u)^T P_z, \\
    P_z &= \big( \calN - \widetilde{\calN} \big) S(z).
  \end{align*}
  Now
  \[
  \Phi(u) = \nabla_g S(u)^T \big( \calN - \widetilde{\calN} \big)\big[ S(u) - S(z) \big]  + \nabla_g h(u) + \nabla_g S(u)^T P_z.
  \]
  Since $z_-$, $z_+$ are hyperbolic critical points of $h$, the map $u \mapsto \nabla_g h(u)$ is $C^3$ smooth from $\calP(z_-,z_+)$ into $L^2(\R,\R^d)$,
  see also \cite{schwarz1993morse}.
  So we just have to prove the mapping property for $\nabla_g S(u)^T \big( \calN - \widetilde{\calN} \big)\big[ S(u) - S(z) \big]$ and the remainder term $\nabla_g S(u)^T P_z$.

  For both these terms, we remark that the prefactor $\nabla_g S(u)^T$ is a $C^3$ smooth map from $\calP(z_-,z_+)$ into $L^\infty(\R,\Mat_{d\times D}(\R))$.
  Indeed, by Morrey's inequality the inclusion $\calP(z_-,z_+) \hookrightarrow L^\infty(\R,\R^d)$ is a smooth map.
  By the hypotheses on $g$ and $S$, the map $\nabla_g S(\cdot)^T : \R^d \to \Mat_{d \times D}(\R)$ is $C^3$ smooth.
  Hence the Nemytskii operator
  \[
  \nabla_g S(\cdot)^T : \calP(z_-,z_+) \to L^\infty(\R,\Mat_{d \times D}(\R)), \qquad u(x) \mapsto \nabla_g S(u(x))^T
  \]
  is $C^3$ smooth.
  It thus suffices to check that the functions $\big( \calN - \widetilde{\calN} \big)\big[ S(u) - S(z) \big]$ and $P_z$ map $\calP(z_-,z_+)$ into $L^2(\R,\R^D)$.
  
  Let $R := \| u - z \|_{L^\infty(\R,\R^d)} + \| z \|_{L^\infty(\R,\R^d)}$, and $L_{S,R} := \sup_{|v| \leq R} \| DS(v) \|$, so that by the mean value theorem
    we obtain the Lipschitz estimate
    \[
    | S(u(x)) - S(z(x)) | \leq L_{S,R} | u(x) - z(x) |.
    \]
  By Young's convolution inequality, we have
  \begin{equation}
    \label{eq:TWN_nonlocal_part_L2_continuous}
    \begin{split}
     \big\| \big( \calN - \widetilde{\calN} \big)\big[ S(u) - S(z) \big]  \big\|_{L^2(\R,\R^D)} &\leq C \| S(u) - S(z) \|_{L^2(\R,\R^D)} \\
    &\leq C L_{S,R} \| u - z \|_{L^2(\R,\R^d)} < \infty,
    \end{split}
  \end{equation}
  where 
  \[
  C = \| N(\cdot) \|_{L^1(\R,\Mat_d(\R))} + \| (N_j)_j \|_{\ell^1(\Z,\Mat_d(\R))} + \| \widetilde{\calN} \|.
  \]
  Hence $u \mapsto \big( \calN - \widetilde{\calN} \big)\big[ S(u) - S(z) \big]$ maps $\calP(z_-,z_+)$ into $L^2(\R,\R^D)$.
  
  For the remainder term $P_z$, remark that for any $x \leq -1$ one has
  \begin{align*}
    P_z(x) &= \calN[S(z)](x) - \widetilde{\calN}S(z_-) = \calN[S(z) - S(z_-)](x) \\
    &= \int_{-\infty}^{x+1} N(y) \big( S(z(x-y)) - S(z_+) \big) \d y + \sum_{\substack{ j\in\Z \\ \lambda_j < x+1}} N_j \big( S(z(x-\lambda_j)) - S(z_+) \big),
  \end{align*}
  so that, in light of the exponential localisation of $N(\cdot)$ and $(N_j)_j$, we obtain the estimate $\| P_z \|_{L^2((-\infty,-1),\R^D)} < \infty$.
  Similarly we obtain $\| P_z \|_{L^2((1,\infty),\R^D)} < \infty$.
  Clearly $\| P_z \|_{L^2((-1,1),\R^D)} \leq 2 \| P_z \|_{L^\infty((-1,1),\R^D)} < \infty$.
  We conclude that $P_z \in L^2(\R,\R^D)$, as desired.
\end{myproof}

\begin{mylemma}
\label{lemma:Phi_L2_smooth}
  The map $\Phi : \calP(z_-,z_+) \to L^2(\R,\R^d)$ is $C^3$.
\end{mylemma}
\begin{myproof}
  Again we make use of the decomposition
  \[
  \Phi(u) = \nabla_g S(u)^T \big( \calN - \widetilde{\calN} \big)\big[ S(u) - S(z) \big]  + \nabla_g h(u) + \nabla_g S(u)^T P_z.
  \]
  We already remarked that $\nabla_g h(\cdot) : \calP(z_-,z_+) \to L^2(\R,\R^d)$ is $C^3$ smooth.
  Likewise, we argued that $\nabla_g S(\cdot)^T : \calP(z_-,z_+) \to L^\infty(\R,\Mat_{d \times D}(\R))$ is $C^3$ smooth.
  The remainder $P_z$ is independent of $u$ hence trivially depends smoothly on $u$.

  Finally, consider the map $\big( \calN - \widetilde{\calN} \big)\big[ S(\cdot) - S(z) \big] : \calP(z_-,z_+) \to L^2(\R,\R^D)$.
  Let $R := 2 \| u_0 \|_{L^\infty(\R,\R^d)}$, and $L_{S,R} := \sup_{|v| \leq R} \| DS(v) \|$.
  Using an estimate identical to \eqref{eq:TWN_nonlocal_part_L2_continuous}, 
  we find that whenever $u \in \calP(z_-,z_+)$ is such that $\| u - u_0 \|_{L^\infty(R,\R^d)} \leq \| u_0 \|_{L^\infty(\R,\R^d)}$, it follows that
    \begin{align*}
    \bigg\| \big( \calN - \widetilde{\calN} \big)&\big[ S(u) - S(z) \big]  - \big( \calN - \widetilde{\calN} \big)\big[ S(u_0) - S(z) \big] \bigg\|_{L^2(\R,\R^d)} \\
   &= \| \big( \calN - \widetilde{\calN} \big)\big[ S(u) - S(u_0) \big] \|_{L^2(\R,\R^d)} \\
   &\leq C L_{S,R} \| u - u_0 \|_{L^2(\R,\R^d)}.
    \end{align*}
  This established continuity.
  Finally, $C^3$ smoothness of the map $\big( \calN - \widetilde{\calN} \big)\big[ S(\cdot) - S(z) \big] : \calP(z_-,z_+) \to L^2(\R,\R^D)$
  is established by differentiation under the integral sign, which is justified by the rapid decay of the convolution kernels.
\end{myproof}

\paragraph{The map $\Psi$.}
By definition, any $\Psi \in \bPsi_h$ is a map $\Psi : C^0_{\text{loc}}(\R,\R^d) \to C^1_{\text{loc}}(\R,\R^d)$.
As such, it can be restricted to the path space $\calP(z_-,z_+)$.
We now consider some of the properties of this restriction.

\begin{mylemma}
\label{lemma:Psi_L2_map}
  Any $\Psi \in \bPsi_h$ restricts to a map $\Psi : \calP(z_-,z_+) \to L^2(\R,\R^d)$.
\end{mylemma}
\begin{myproof}
  Recall from \eqref{eq:Psi_ell_perturbation} that $\Psi$ can be written as
  \[
  \Psi(u) = \sum_{\theta \in \Theta_h} \Psi_\theta(u), \qquad \rst{\Psi_\theta(u)}{x} = \rst{ \alpha_\theta * \psi_\theta(u) }{(x,x)},
  \]
  where
   \[
  \rst{ \psi_\theta(u) }{ (\xi_1, \xi_2) } = \rst{ \sigma_{\ell,z_1}(u) }{ \xi_1 } \rst{ \sigma_{\ell,z_2}(u) }{ \xi_2 }, \qquad \theta = (\ell,z_1,z_2) \in \Theta_h.
  \]
  By definition we have $0 \leq \sigma_{\ell,z_2}(u) \leq 1$, which allows us to estimate
  \[
  | \rst{ \Psi_\theta(u) }{x} | \leq \rst{ \beta_\theta * \sigma_{\ell,z_1}(u)}{x}, \qquad \beta_\theta(r) := \int_\R  | \alpha_\theta( r,s ) | \d s.
  \]
  Applying Young's convolution inequality, we have
  \[
  \| \Psi_\theta(u) \|_{L^2(\R,\R^d)} \leq \| \beta_\theta * \sigma_{\ell,z_1}(u) \|_{L^2(\R,\R^d)} \leq \| \alpha_\theta \|_{L^1(\R^2,\R^d)} \| \sigma_{\ell,z_1}(u) \|_{L^2(\R,\R)},
  \]
  where we used $\| \beta_\theta \|_{L^1(\R,\R^d)} = \| \alpha_\theta \|_{L^1(\R^2,\R^d)}$.
  We will now bound $\| \sigma_{\ell,z_1}(u) \|_{L^2(\R,\R)}$, uniformly in $\ell$ and $z_1$.
  Set $R := \| u \|_{L^\infty(\R,\R^d)}$, let $\rho : \crit(h) \to (0,\infty]$ be as given in \eqref{eq:TWN_Psi_rho_z}, and define
  \[
  \rho_R := \inf_{ \substack{ z \in \crit(h) \\ | z | \leq R } } \rho(z).
  \]
  The definition of $\sigma_{\ell,z_1}$ ensures that, whenever $\inf_{z\in\crit(h)} |u(x) - z| < \rho_R$, it follows that $\rst{ \sigma_{\ell,z_1}(u) }{ x } = 0$.
  Therefore,
  \[
    \| \sigma_{\ell,z_1}(u) \|_{L^2(\R,\R)} \leq \sqrt{ \vol\left( \set{ x \in \R }{ \inf_{z\in\crit(h)} |u(x) - z| \geq \rho_R } \right) }.
  \]
  The right hand side is finite since $u \in \calP(z_-,z_+)$.
  Hence the series $\Psi(u) = \sum_{\theta \in \Theta_h} \Psi_\theta(u)$ converges in $L^2(\R,\R^d)$, and 
   \[
    \| \Psi(u) \|_{L^2(\R,\R^d)} \leq \| \Psi \|_{\bPsi_h} \sqrt{ \vol\left( \set{ x \in \R }{ \inf_{z\in\crit(h)} |u(x) - z| \geq \rho_R } \right) }.
  \]
  In fact, this estimate holds uniformly in $u$, whenever $\| u \|_{L^\infty(\R,\R^d)} \leq R$.
\end{myproof}

\begin{mylemma}
\label{lemma:Psi_L2_smooth}
  The map $\Psi : \calP(z_-,z_+) \to L^2(\R,\R^d)$ is $C^3$.
\end{mylemma}
\begin{myproof}
Let us first address continuity.
Let $(u_n)_n \subset \calP(z_-,z_+)$ be a sequence, converging to $u_\infty \in \calP(z_-,z_+)$.
Then
\[
\rst{ \Psi(u_\infty) }{x} - \rst{ \Psi(u_n) }{x} = \sum_{\theta \in \Theta_h} \rst{ \alpha_\theta * ( \psi_\theta(u_\infty) - \psi_\theta(u_n) ) }{(x,x)}.
\]
Note that, with $\theta = (\ell,z_1,z_2)$,
\begin{align*}
  \rst{( \psi_\theta(u_\infty) - \psi_\theta(u_n) )}{(\xi_1,\xi_2)} &= \rst{( \sigma_{\ell,z_1}(u_\infty) - \sigma_{\ell,z_1}(u_n) )}{\xi_1} \rst{\sigma_{\ell,z_2}(u_\infty)}{\xi_2} \\
&\quad + \rst{\sigma_{\ell,z_1}(u_n)}{\xi_1} \rst{( \sigma_{\ell,z_2}(u_\infty) - \sigma_{\ell,z_2}(u_n) )}{\xi_2}.
\end{align*}
As each term in the right hand side is the product of bounded functions in $\xi_1$ with bounded functions in $\xi_2$, 
we may use an identical estimate as the one used in the proof of Lemma \ref{lemma:Psi_L2_map} to derive that
\begin{equation}
  \label{eq:Psi_continuity_estimate}
  \begin{split}
    \| \Psi(u_\infty) - \Psi(u_n) \|_{L^2(\R,\R^d)} & \leq \sum_{\theta = (\ell,z_1,z_2) \in \Theta_h} \| \alpha_\theta \|_{L^1(\R^2,\R^d)}  \| \sigma_{\ell,z_1}(u_\infty) - \sigma_{\ell,z_1}(u_n) \|_{L^2(\R,\R)} \\
    &\quad + \sum_{\theta = (\ell,z_1,z_2) \in \Theta_h} \| \alpha_\theta \|_{L^1(\R^2,\R^d)} \| \sigma_{\ell,z_2}(u_\infty) - \sigma_{\ell,z_2}(u_n) \|_{L^2(\R,\R)}.
  \end{split}
\end{equation}
Now remark that the Nemytskii operator
\[
\sigma_{\ell,z_i}(\cdot) : \calP(z_-,z_+) \to L^2(\R,\R), \qquad u(x) \mapsto \sigma_{\ell,z_i}(u(x))
\]
is smooth.
This follows from the pointwise smoothness of $\sigma_{\ell,z_i} : \R^d \to \R$, and the fact that $\sigma_{\ell,z_i}$ maps $\calP(z_-,z_+)$ into compactly supported functions.
We now conclude from \eqref{eq:Psi_continuity_estimate} that $\Psi(u_n) \to \Psi(u_\infty)$ in $L^2(\R,\R^d)$ as $n\to\infty$. 

As for the smoothness of $\Psi$, recall from Lemma \ref{lemma:TWN_Psi_uniform_smooth} that $\Psi : L^\infty(\R,\R^d) \to L^\infty(\R,\R^d)$ is $C^3$ smooth, with for $k \in \{0,\dots,3\}$
\begin{equation}
  \label{eq:Psi_differentiation_identity}
  \rst{ D^k \Psi(u)[ v_1 , \dots , v_k ] }{ x } = \sum_{\theta\in \Theta_h} \rst{ \alpha_\theta * D^k \psi_\theta(u)[ v_1 , \dots , v_k ] }{(x,x)}
\end{equation}
for $v_1,\dots,v_k \in L^\infty(\R,\R^d)$.
It thus suffices to check that the series in the right hand side of \eqref{eq:Psi_differentiation_identity} converges in $L^2(\R,\R^d)$,
whenever $u \in \calP(z_-,z_+)$ and $v_1, \dots , v_k \in W^{1,2}(\R,\R^d)$, and depends continuously on $u$.
Here we note that, just like $\rst{ \psi_\theta(u) }{(\xi_1,\xi_2)}$, the expression $\rst{ D^k \psi_\theta(u)[v_1,\dots,v_k] }{(\xi_1,\xi_2)}$ can be written as (linear combination of)
products of compactly supported functions in $\xi_1$ and compactly supported functions in $\xi_2$.
We can thus apply estimates identical to those used in proving the continuity of $\Psi$ to conclude $\Psi : \calP(z_-,z_+) \to L^2(\R,\R^d)$ is $C^3$ smooth.
\end{myproof}

\paragraph{Smoothness of \eqref{eq:TWN} on path space.}
Combining the previous lemmata, we obtain the following.

\begin{mythm}
  \label{thm:TWN_path_space}
  The map
  \[
  \calP(z_-,z_+) \times \bPsi_h \to L^2(\R,\R^d), \qquad (u,\Psi) \mapsto u' + \Phi(u) + \Psi(u)
  \]
  is $C^3$.
\end{mythm}
\begin{myproof}
  Smoothness of the map $u \mapsto u' + \Phi(u)$ is already established in Lemma \ref{lemma:Phi_L2_smooth}.
  By Lemma \ref{lemma:Psi_L2_smooth} the map
  \[
  \calP(z_-,z_+) \times \bPsi_h \to L^2(\R,\R^d), \qquad (u,\Psi) \mapsto \Psi(u)
  \]
  is $C^3$ in $u$.
  Since this map depends linearly on $\Psi$, the $C^3$ smoothness follows from a diagonal argument,
  comparable to the proof of Lemma \ref{thm:TWN_weak_continuity}.
\end{myproof}

\subsection{Fredholm indices}

We will now investigate Fredholm properties of the map
\[
\partial_x + \Phi(\cdot) + \Psi(\cdot) : \calP(z_-,z_+) \to L^2(\R,\R^d).
\]
We first consider the case where $z_- = z_+$.
We recall from Definition \ref{def:TWN_hyperbolic} that a critical point $z \in \crit(h)$ is said to be hyperbolic
if the Hessian matrix $\Hess_g h(z)$ is nondegenerate.
\begin{mylemma}
  \label{lemma:hyperbolicity_implies_invertible}
  A critical point $z \in \crit(h)$ is hyperbolic if and only if 
  \[
  \partial_x + D\Phi(z) : W^{1,2}(\R,\R^d) \to L^2(\R,\R^d)
  \]
  has a bounded inverse.
\end{mylemma}
\begin{myproof}
  Using the Fourier transform, solving $v' + D\Phi(z) v = w$ is equivalent to solving
  \begin{equation}
    \label{eq:hyperbolicity_analytic_Fourier}
    L  V = W, \qquad \langle \cdot \rangle V(\cdot) \in L^2(\R,\C^d), \; W \in L^2(\R,\C^d),
  \end{equation}
where $\langle \xi \rangle =  \sqrt{ 1 + |\xi|^2 }$ and $L : \R \to \calL(\C^d , \C^d)$ is the multiplication operator defined by
\[
L(\xi) v = i \xi v + \nabla_g S(z)^T \widehat \calN(\xi)[ DS(z) v ] + P_z v,
\]
with $P_z \in \calL(\C^d,\C^d)$ defined by
\[
P_z v = \big( D[ \nabla_g S^T ](z)v \big) \widetilde{\calN} S(z) + D [ \nabla_g F ](z)v,
\]
and
\[
\widehat \calN(\xi) = \int_\R N(x) e^{- i x \xi} \d x + \sum_{j\in\Z} N_j e^{- i \lambda_j \xi}.
\]
The exponential localisation in Hypothesis \ref{hyp:calN} ensures that $\widehat \calN \in W^{1,2}(\R,\Mat_{D\times D}(\C))$.
Thus by Morrey's inequality $\widehat \calN$ is continuous and $\| \widehat \calN(\xi) \|_{L^\infty(\R,\Mat_{D \times D}(\C))} < \infty$.
Consequently, $L$ is continuous.
Invertibility of \eqref{eq:hyperbolicity_analytic_Fourier} thus implies that $L(0)$ is invertible.
Now note that
\[
L(0) = \nabla_g S(z)^T \widetilde{\calN} D S(z) + \big( D[ \nabla_g S^T ](z)[\cdot] \big) \widetilde{\calN} S(z) + D[\nabla_g F](z) = D[ \nabla_g h ](z),
\]
which is invertible if and only if $\Hess_g h(z)$ is invertible.
Thus invertibility of
\[
  \partial_x + D\Phi(z) : W^{1,2}(\R,\R^d) \to L^2(\R,\R^d)
\]
implies hyperbolicity of $z$ as a critical point of $h$.

Now suppose $z$ is a hyperbolic critical point of $h$.
As we saw, this is equivalent to $L(0)$ being invertible.
Then, by a Von Neumann series, there exists $\epsilon$ so that $L(\xi)$ is invertible, for $|\xi| < \epsilon$, with $\sup_{|\xi| < \epsilon} \| L(\xi)^{-1} \| < \infty$.
Since by Hypothesis \ref{hyp:calN} we have $N(r) = N(-r)$, $N_j = N_{-j}$, and $\lambda_j = -\lambda_{-j}$, 
we find that $\imag{ \widehat \calN(\xi) } = 0$ for any $\xi \in \R$.
Hence $\imag{ L(\xi) } = \xi \id$ is a nonzero multiple of the identity for $\xi \neq 0$.
Thus $L(\xi)$ is invertible for all values of $\xi \in \R$.
Since $L(\xi) - i \xi \id$ is uniformly bounded in $\xi \in \R$, it follows that for large values of $|\xi|$ the operator $A(\xi) := \xi^{-1} L(\xi)$ is a small perturbation of the identity.
Hence $\| A(\xi)^{-1} \| \leq 2$ for sufficiently large values of $|\xi|$.
It follows that $\| L(\xi)^{-1} \| \leq 2 / |\xi|$ for sufficiently large $\xi$.
We conclude that
\[
\| L(\xi)^{-1} \| \leq \frac{ C }{ 1 + |\xi| }, \qquad \xi \in \R
\]
for some $C > 0$.
This implies \eqref{eq:hyperbolicity_analytic_Fourier} is solvable, hence proves that hyperbolicity of $z$ as a critical point of $h$ is sufficient to ensure
invertibility of the map $\partial_x + D\Phi(z)$.
\end{myproof}

Since we assume $h$ is a Morse function, all critical points of $h$ are hyperbolic, 
so that by the preceding lemma the linearisation of \eqref{eq:TWN} around a constant solution is always invertible.

The following lemma indicates that $\partial_x + \Phi(\cdot) + \Psi(\cdot)$ and $\partial_x + \Phi(\cdot)$ belong to the same Fredholm region.
\begin{mylemma}
\label{lemma:DPsi_compact}
  For any $u \in \calP(z_-,z_+)$, the linearisation $D \Psi(u) : W^{1,2}(\R,\R^d) \to L^2(\R,\R^d)$ is a compact operator.
\end{mylemma}
\begin{myproof}
  By (the proof of) Lemma \ref{lemma:Psi_L2_smooth}, we have $D \Psi(u) = \sum_{\theta \in \Theta_h} D \Psi_\theta(u)$, where the series converges in the operator norm.
  Hence it suffices to see that $D \Psi_\theta(u)$ is compact.
  Now $D \Psi_\theta(u)$ is given by
  \[
  \rst{ D\Psi_\theta(u) v }{x} = \rst{ \alpha_\theta * \big[ D \psi_\theta(u) v \big] }{ (x,x) }, \qquad v \in W^{1,2}(\R,\R^d),\; x\in\R,
  \]
  where the convolution with $\alpha_\theta$ defines a bounded linear map from $L^2(\R^2,\R)$ into $L^2(\R^2,\R^d)$.
  Hence it suffices to see that $D \psi_\theta(u) : W^{1,2}(\R,\R^d) \to L^2(\R^2,\R)$ is compact.
  With $\theta = (\ell , z_1 , z_2)$, we have
  \begin{equation}
    \label{eq:psi_compact}
    \rst{ D \psi_\theta(u) v }{(\xi_1,\xi_2)} = \rst{ \sigma_{\ell,z_2}(u) }{\xi_2} \rst{ \big[ D \sigma_{\ell,z_1}(u) v \big] }{\xi_1} + \rst{ \sigma_{\ell,z_1}(u) }{\xi_1} \rst{\big[ D \sigma_{\ell,z_2}(u) v \big] }{\xi_2}.
  \end{equation}
  
  By Morrey's inequality, the function $u : \R \to \R^d$ is $\frac 1 2$-H\"older continuous.
  Let $B \subset \R^d$ be a bounded neighbourhood of $u(\R)$.
  Then the function $D \sigma_{\ell,z_i} : \R^d \to \calL(\R^d,\R)$, with $i \in \{1,2\}$, restricts to a smooth bounded Lipschitz map $\rst{D \sigma_{\ell,z_i}}{B} : B \to \calL(\R^d,\R)$.
  Hence the composition $D \sigma_{\ell,z_i}(u) : \R \to \calL(\R^d,\R)$ is $\frac 1 2$-H\"older continuous and uniformly bounded.
  Moreover, the fact that $u(x) \to z_\pm$ as $x \to \pm \infty$, combined with the definition of $\sigma_{\ell,z_i}$, ensures that $D \sigma_{\ell,z_i}(u)$ is compactly supported.

  Now let $(v_n)_n \subset W^{1,2}(\R,\R^d)$ be a bounded sequence.
  By Morrey's inequality, there exists $R, L > 0$ so that
  \[
  \sup_n \| v_n \|_{L^\infty(\R,\R^d)} \leq R, \qquad \sup_n | v_n(x) - v_n(y) | \leq L | x - y |^{1/2}.
  \]
  With $i \in \{1,2\}$, we now have $\sup_n \| D \sigma_{\ell,z_i}(u) v_n \|_{L^\infty(\R,\R^d)} < \infty$ and
  \begin{align*}
    | D \sigma_{\ell,z_i}(u(x)) v_n(x) - D \sigma_{\ell,z_i}(u(y)) v_n(y) | &\leq \| D \sigma_{\ell,z_i}(u(x)) \| | v_n(x) - v_n(y) | \\
    &\quad + \| D \sigma_{\ell,z_i}(u(x)) - D \sigma_{\ell,z_i}(u(y)) \| |v_n(y)| \\
    &\leq C | x - y |^{1/2},
  \end{align*}
  for some $C > 0$ independent of $n$.
  Since $\supp( D \sigma_{\ell,z_i}(u) v_n ) \subset [-r,r]$ for some $r > 0$ independent of $n$, 
  we now find by the Arzel\`a-Ascoli theorem that $(D \sigma_{\ell,z_i}(u) v_n)_n$ converges uniformly over a subsequence.
  Consequently $(D \sigma_{\ell,z_i}(u) v_n)_n$ converges over a subsequence in $L^2(\R,\R)$.
  This shows that the map $v \mapsto D \sigma_{\ell,z_i}(u) v$ is compact from $W^{1,2}(\R,\R^d)$ into $L^2(\R,\R)$.
  Consequently, recalling \eqref{eq:psi_compact}, we see that
  \[
  D\psi_\theta(u) : W^{1,2}(\R,\R^d) \to L^2(\R^2,\R)
  \]
  is compact.
  It follows that $D\Psi(u) : W^{1,2}(\R,\R^d) \to L^2(\R,\R^d)$ is compact, as claimed.
\end{myproof}

We can now establish the Fredholm alternative for $\partial_x + \Phi(\cdot) + \Psi(\cdot)$.
We recall that a nonlinear map $f : X \to Y$, with $X$, $Y$ connected smooth manifolds,
is said to be a Fredholm map if it is $C^1$ and $Df(x) : T_x X \to T_{f(x)} Y$
is a Fredholm operator, for each $x \in X$.
By continuity the index is then independent of the basepoint $x \in X$.
\begin{mythm}
\label{thm:TWN_Fredholm}
  The function
  \[
  \partial_x + \Phi(\cdot) + \Psi(\cdot) : \calP(z_-,z_+) \to L^2(\R,\R^d)
  \]
  is a Fredholm map, with index given by
  \[
  \ind\big( \partial_x + \Phi(\cdot) + \Psi(\cdot) \big) = m_h(z_-) - m_h(z_+),
  \]
  where $m_h$ denotes the Morse index defined in Section \ref{subsec:Morse_function_h}.
\end{mythm}
\begin{myproof}
Pick arbitrary $u \in \calP(z_-,z_+)$.
By Lemma \ref{lemma:DPsi_compact} the Fredholm properties of the map $\partial_x + D \Phi(u) + D\Psi(u)$
are the same as those of $\partial_x + D\Phi(u)$.
The operator 
\[
\partial_x + D \Phi(u) : W^{1,2}(\R,\R^d) \to L^2(\R,\R^d) 
\]
fits within the framework of \cite{faye2013fredholm}.
The asymptotic operators $\partial_x + D\Phi(z_\pm)$ are invertible by Lemma \ref{lemma:hyperbolicity_implies_invertible},
hence it now follows from \cite{faye2013fredholm} that $\partial_x + D\Phi(u)$ is Fredholm.
Furthermore, the index depends only on the asymptotic operators $D\Phi(z_\pm)$, i.e.,
there exists a relative index  $\nu(\cdot,\cdot)$ such that
\[
\ind( \partial_x + D \Phi(u) ) = \nu( D\Phi(z_-) , D\Phi(z_+) ), \qquad \text{for all} \quad u \in \calP(z_-,z_+).
\]
This shows that linear map $\partial_x + D\Phi(u) + D\Psi(u)$ is Fredholm, and as the index only depends on the asymptotic operators $D\Phi(z_\pm)$,
the function $\partial_x + \Phi(\cdot) + \Psi(\cdot)$ is a Fredholm map. 

Left to prove is the expression of the Fredholm index in terms of classical Morse indices.
We will show this by means of a continuation argument.
For $\beta \in [0,1]$ define the linear map $\calT^\pm_\beta : W^{1,2}(\R,\R^d) \to L^2(\R,\R^d)$ by
\[
\calT^\pm_\beta v := \nabla_g S(z)^T \big( \beta \calN + (1 - \beta ) \widetilde{\calN} \big)\big[ D S(z_\pm) v \big] + P_{z_\pm} v,
\]
where
\[
P_{z_\pm} v = \big( D[ \nabla_g S^T ](z_\pm)v \big) \widetilde{\calN} S(z_\pm) + D [ \nabla_g F ](z_\pm)v.
\]
Then $\calT^\pm_0 = D \Phi(z_\pm)$ and $\calT^\pm_1 = D[ \nabla_g h(z_\pm) ]$.

Now choose a smooth function $\beta : \R \to [0,1]$ with $\beta(x) = 1$ for $x \leq -1$ and $\beta(x) = 0$ for $x \geq 1$,
and define the operators $\calT^\pm : W^{1,2}(\R,\R^d) \to L^2(\R,\R^d)$ by
\[
\calT^\pm[ v ](x) = \calT^\pm_{\beta(x)}[v](x).
\]
Now the results from \cite{faye2013fredholm} are applicable to operators $\partial_x + \calT^\pm$, showing that $\partial_x + \calT^\pm$ are Fredholm with
\[
\ind( \partial_x + \calT^\pm ) = \nu( \calT^\pm_{-\infty} , \calT^\pm_{\infty} ) = \nu( D \Phi(z_\pm) , D[ \nabla_g h ](z_\pm) ).
\]
Classical Morse theory tells us that
\[
\ind( \partial_x + D[\nabla_g h](u) ) = m_h(z_-) - m_h(z_+),
\]
see for instance \cite{schwarz1993morse, robbin1995spectral}.
By transitivity of the relative index we have
\begin{align*}
  \ind( \partial_x + D \Phi(u) ) &=  \nu ( D[\nabla_g h](z_-) , D[\nabla_g h](z_+) ) + \nu( D \Phi(z_-) , D[\nabla_g h](z_-) )\\
  &\quad + \nu( D[\nabla_g h](z_+) , D \Phi(z_+) ) \\
  &= \ind( \partial_x + D[\nabla_g h](u) ) + \ind( \partial_x + \calT^- ) - \ind(\partial_x + \calT^+).
\end{align*}
Combining there observations yields
\[
\ind( \partial_x + D \Phi(u) ) = m_h(z_-) - m_h(z_+) + \ind(\partial_x + \calT^-) - \ind(\partial_x + \calT^+).
\]

Left to prove is that $\ind(\partial_x + \calT^\pm) = 0$.
The spectral flow formula developed in \cite{faye2013fredholm} shows that the Fredholm index of $\partial_x + \calT^\pm$ 
is computable by counting the spectral crossings of the homotopy of operators $( \calT^\pm_\beta )_\beta$.
That is, as we vary $\beta$ we need to count roots $\xi$ of the equation
\begin{equation}
  \label{eq:TWN_SF}
  \det\bigg( i \xi \id + \widehat{ \calT^\pm_\beta }(\xi) \bigg) = 0, \qquad \xi \in \R,
\end{equation}
where
\[
\widehat{ \calT^\pm_\beta }(\xi) = \nabla_g S(z_\pm)^T \big( \beta \widehat \calN(\xi) + (1 - \beta) \widetilde{\calN} \big) D S(z_\pm) + P_{z_\pm},
\]
with
\[
\widehat \calN(\xi) = \widehat N(\xi) + \sum_{j\in\Z} N_j e^{-\lambda_j \xi}.
\]
Since by Hypothesis \ref{hyp:calN} we have $N(r) = N(-r)$, $N_j = N_{-j}$, and $\lambda_j = -\lambda_{-j}$, 
we find that $\imag{ \widehat \calN(\xi) } = 0$ for any $\xi \in \R$.
Hence the only possible solution of \eqref{eq:TWN_SF} would be $\xi = 0$.
Now we note
\[
\widehat{ \calT^\pm_\beta }(0) = \nabla_g S(z_\pm)^T \widetilde{\calN} D S(z_\pm) + P_{z_\pm} = D[ \nabla_g h ](z_\pm),
\]
so that by hyperbolicity of $z_-$ and $z_+$ we conclude that \eqref{eq:TWN_SF} does not have any solutions.
This is true for arbitrary $\beta \in [0,1]$, which implies there are no spectral crossings.
It follows from the spectral flow formula that $\ind( \partial_x + \calT^\pm ) = 0$.
\end{myproof}


\section{Some convergence results}

\label{sec:TWN_compactness}

In this section we will study various convergence and compactness properties of solutions of \eqref{eq:TWN} which are bounded in $L^\infty$.

\subsection{Approximate solutions and Palais--Smale sequences}
We start out by defining the notion of approximate solutions of \eqref{eq:TWN}
\begin{mydef}[Approximate solutions]
  Given $\delta > 0$, a function $u \in C^1_{\text{loc}}(\R,\R^d)$ is said to be a $\delta$-approximate solution of \eqref{eq:TWN} if
  \[
  \| u' + \Phi(u) + \Psi(u) \|_{L^\infty(\R,\R^d)} \leq \delta.
  \]
\end{mydef}

We may also consider sequences of approximate solutions.
If $\delta \to 0$ along such a sequence, we call it a Palais--Smale sequence.
\begin{mydef}[Palais--Smale sequences]
  A sequence $(u_n)_n \subset C^1_{\text{loc}}(\R,\R^d)$ is said to be a Palais--Smale sequence provided that
  \[
  \sup_n \| u_n \|_{L^\infty(\R,\R^d)} < \infty
  \]
  and
  \[
  \| u_n' + \Phi(u_n) + \Psi(u_n) \|_{L^\infty(\R,\R^d)}  \to 0, \qquad \text{as} \quad n \to \infty.
  \]
\end{mydef}

We now state a regularity result for approximate solutions of \eqref{eq:TWN}.

\begin{mythm}
\label{thm:TWN_weak_compactness}
  Let $\delta > 0$ and suppose $(u_n)_n \subset C^1_{\text{loc}}(\R,\R^d)$ is a sequence of $\delta$-approximate solutions of \eqref{eq:TWN}.
  Assume further that
  \[
    \sup_n \| u_n \|_{L^\infty(\R,\R^d)} < \infty.
  \]
  Then
  \[
  \sup_n \| u_n \|_{W^{2,\infty}(\R,\R^d)} < \infty.
  \]
  Furthermore, there exists a subsequence $(u_{n_k})_k$ such that $u_{n_k} \to u_\infty$ in $C^1_{\text{loc}}(\R,\R^d)$ as $k \to \infty$,
  where $u_\infty$ is another $\delta$-approximate solution of \eqref{eq:TWN}.
  If we assume moreover that $(u_n)_n$ is Palais--Smale, then the limit $u_\infty$ is a solution of \eqref{eq:TWN}.
\end{mythm}
\begin{myproof}
    In light of Lemma \ref{lemma:TWN_Phi_uniform_bound} and Lemma \ref{lemma:TWN_Psi_uniform_bound}, 
  we find that 
  \[
  \sup_n \| \Phi(u_n) + \Psi(u_n) \|_{L^\infty(\R,\R^d)} < \infty.
  \]
   Since $- u_n' = \Phi(u_n) + \Psi(u_n)$ we thus obtain a uniform bound
\[
\sup_n \| u_n \|_{W^{1,\infty}(\R,\R^d)} < \infty.
\]
  Applying Lemmata \ref{lemma:TWN_Phi_uniform_bound} and \ref{lemma:TWN_Psi_uniform_bound} once again,
  we now obtain a uniform bound 
\[
\sup_n \| u_n \|_{W^{2,\infty}(\R,\R^d)} < \infty.
\]
  Combining these bootstrap estimates with the Arzel\`a--Ascoli theorem,
  we extract a subsequence $(u_{n_k})_k$ converging in $C^1_{\text{loc}}(\R,\R^d)$ to $u_\infty$.
  Theorem \ref{thm:TWN_weak_continuity} implies that, for any bounded interval $I \subset \R$, we have
  \begin{align*}
    \| u_\infty' + \Phi(u_\infty) + \Psi(u_\infty) \|_{L^\infty(I,\R^d)} &= \lim_{k\to\infty} \| u_{n_k}' + \Phi(u_{n_k}) + \Psi(u_{n_k}) \|_{L^\infty(I,\R^d)} \\
    &\leq \limsup_{k\to\infty} \| u_{n_k}' + \Phi(u_{n_k}) + \Psi(u_{n_k}) \|_{L^\infty(\R,\R^d)} \\
    &\leq \delta.
  \end{align*}
  Hence $u_\infty$ is a $\delta$-approximate solution of \eqref{eq:TWN}.
  The statement about Palais--Smale sequences readily follows from this estimate.
\end{myproof}

\subsection{Asymptotics}
\label{subsec:TWN_asymptotics}
As we lack the notion of a phase space, we first introduce the notion of $\alpha$- and $\omega$-limit sets in our setup.
\begin{mydef}[$\alpha$- and $\omega$-limit sets]
\label{def:TWN_limit_sets}
  Given $E \subset C_{\text{loc}}^1(\R,\R^d)$,
we define $\alpha(E)$ to be the $\alpha$-limit set of $E$ with respect to the shift dynamics on $C_{\text{loc}}^1(\R,\R^d)$.
Thus, $\alpha(E)$ consists of the accumulation points of $\set{ \tau \cdot u }{ \tau < 0,\; u \in E }$.
Similarly, the $\omega$-limit set $\omega(E)$ is defined as the set of accumulation points of  $\set{ \tau \cdot u }{ \tau > 0,\; u \in E }$.
\end{mydef}

Given a measurable subset $B \subseteq \R$, we define the restricted kinetic energy $\calE_{\text{kin}}(\cdot | B)$ by
\[
\calE_{\text{kin}}(u | B) := \frac 1 2 \int_B g_{u(x)}(u'(x),u'(x)) \d x.
\]
We are now prepared to study the asymptotic behaviour of approximate solutions with finite kinetic energy.
\begin{mylemma}
\label{lemma:TWN_limit_sets_approximate}
Suppose $u \in C^1_{\text{loc}}(\R,\R^d)$ is a $\delta$-approximate solution of \eqref{eq:TWN},
which satisfies $\calE_{\text{kin}}(u|(0,\infty)) < \infty$. 
Then $\omega(u)$ is nonempty and consists of $\delta$-approximate constant solutions of \eqref{eq:TWN}.
Likewise, if $\calE_{\text{kin}}(u|(-\infty,0)) < \infty$, then $\alpha(u)$ is nonempty and consists of $\delta$-approximate constant solutions of \eqref{eq:TWN}.
\end{mylemma}
\begin{myproof}
  We discuss here the case where $\calE_{\text{kin}}(u|(0,\infty)) < \infty$; the argument for the other case is identical.
  By definition, an element $z \in \omega(u)$ is the limit in $C^1_{\text{loc}}(\R,\R^d)$ of a sequence $u_n(x) := u(x + \tau_n)$, 
  for some sequence $(\tau_n)_n\subset \R$ with  $\tau_n \to \infty$ as $n\to\infty$.
  In light of Theorem \ref{thm:TWN_weak_compactness} the sequence $(u_n)_n$ has a subsequence converging to a $\delta$-approximate solution $u_\infty$ of \eqref{eq:TWN}.
  Since $C^1_{\text{loc}}(\R,\R^d)$ is Hausdorff it follows that $z = u_\infty$.
  Thus any element of $\omega(u)$ is a $\delta$-approximate solutions of \eqref{eq:TWN}.
  The same argument (using, say, the sequence $\tau_n = n$) shows that $\omega(u) \neq \emptyset$.

  To see that any $z \in \omega(u)$ is constant, we further exploit the energy bound.
  First we set $R := \| u \|_{L^\infty(\R,\R^d)}$.
  After recalling that
  \[
  g_u(v_1,v_2) := \big\langle v_1 , G(u) v_2 \big\rangle, \qquad v_1, v_2 \in T_u \R^d = \R^d,
  \]
  we find that
  \[
    |u'(x)|^2 = g_{u(x)}\big( u'(x) ,  G(u(x))^{-1} u'(x) \big) \leq \sup_{|v| \leq R} \| G(v)^{-1} \| \; g_{u(x)}(u'(x),u'(x)).
  \]
  Consequently
  \[
  \| u' \|_{L^2((0,\infty),\R^d)}^2 \leq 2 \sup_{|v| \leq R} \| G(v)^{-1} \| \; \calE_{\text{kin}}(u|(0,\infty)) < \infty.
  \]
  Now, for any $L > 0$
  \[
  \int_{-L}^L | z'|^2 \d x = \lim_{n\to\infty} \int_{-L}^L | u'(x+\tau_n) |^2 \d x \leq \lim_{n\to\infty} \int_{ \tau_n - L }^\infty |u'|^2 \d x = 0.
  \]
  As $z$ is $C^1$, it follows that $z' = 0$ everywhere.
  Hence $z$ is a constant $\delta$-approximate solution of \eqref{eq:TWN}.
\end{myproof}

When $u$ is a solution, instead of an approximate solution, we obtain a slightly more detailed description of the limit sets of $u$
\begin{mylemma}
\label{lemma:TWN_limit_sets}
Let $R \geq 0$ and suppose $\Psi$ is $R$-admissible.
Suppose $u \in C^1_{\text{loc}}(\R,\R^d)$ is solution of \eqref{eq:TWN} with $\| u \|_{L^\infty(\R,\R^d)} \leq R$ and $\calE_{\text{kin}}(u|(0,\infty)) < \infty$.
Then there exists a constant solution $z_+$ of \eqref{eq:TWN} such that
\[
\big( u(x) , u'(x) \big) \to \big( z_+,0 \big), \qquad \text{as} \quad x\to \infty.
\]
Likewise, if $\calE_{\text{kin}}(u|(-\infty,0)) < \infty$, there exists a constant solution $z_-$ of \eqref{eq:TWN} such that
\[
\big( u(x) , u'(x) \big) \to \big( z_-,0 \big), \qquad \text{as} \quad x\to -\infty.
\]
\end{mylemma}
\begin{myproof}
  We consider only the case there $\calE_{\text{kin}}(u|(0,\infty)) < \infty$; the argument for the other case is identical.
  The proof of Lemma \ref{lemma:TWN_limit_sets_approximate} shows that $\omega(u)$ consists of constant solutions of \eqref{eq:TWN}.
  We now observe that $\omega(u)$ is necessarily a path connected subset of $C^1_{\text{loc}}(\R,\R^d)$.
  As $\omega(u)$ consist of constant solutions of \eqref{eq:TWN}, the evaluation map
  \[
  u \mapsto u(0), \qquad C^1_{\text{loc}}(\R,\R^d) \to \R^d
  \]
  sets up a continuous bijection from $\omega(u)$ onto its image.
  Path connectedness is preserved under this identification.
  Henceforth we consider $\omega(u)$ as subsets of $\R^d$.
  Since $\Psi$ is $R$-admissible we find that $\omega(u) \subseteq \crit(h)$, and since $h$ is assumed to be a Morse function, $\crit(h)$ is a discrete subset of $\R^d$.
  Hence $\omega(u)$ is a path connected subset of a discrete space, which implies that $\omega(u) = \{z_+\}$.
\end{myproof}

\subsubsection{Decay rates}
\label{subsec:TWN_exp_decay}

We have established that bounded solutions of \eqref{eq:TWN} with bounded kinetic energy converge towards constant solutions of \eqref{eq:TWN}.
We now further investigate the rate at which this convergence takes place.
We start out with a convergence result for the linearised equation.

For the moment, assume $u \in L^\infty(\R,\R^d)$, without imposing asymptotics as $x\to\pm\infty$.
We define $L_u : W^{1,\infty}(\R,\R^d) \to L^\infty(\R,\R^d)$ to be the linearisation of \eqref{eq:TWN} around $u$, 
\begin{equation}
  \label{eq:TWN_formal_linearisation}
  \begin{split}
    L_u[v](x) &:= v'(x) + D[\nabla_g F](u(x))v(x) + \big( D[ \nabla_g S^T ](u(x))v(x) \big) \calN[S(u)](x) \\
    &\quad + \nabla_g S(u)^T \calN[ DS(u) v ](x) +\sum_{\theta\in \Theta_h} \rst{ \alpha_\theta *\big(  D\psi_\theta(u)v \big)}{(x,x)}.
  \end{split}
\end{equation}

Before stating the next lemma, we introduce exponentially weighted $L^p$ spaces.
Given $\eta \in \R$ and a normed vector space $V$, we denote by $L^p_\eta(\R^m,V) = L^p(\d\mu^m,V)$ 
the space of functions $f \in L^1_{\text{loc}}(\R^m,V)$ which are 
$p$-integrable with respect to the measure $\d \mu^m(x) = e^{\eta |x|} \d \lambda^m(x)$, where $\lambda^m$ is the $m$-dimensional Lebesgue measure on $\R^m$.
Exponentially weighted $\ell^p$ spaces $\ell^p_\eta(\Z,V)$ and Sobolev spaces $W^{k,p}_\eta(\R^m,V) = W^{k,p}(\d\mu^m,V)$ are defined in a similar manner.

\begin{mylemma}
  For $\eta \in [-\eta_0,\eta_0]$, where $\eta_0$ is as in Hypothesis \ref{hyp:calN},
  the linear map $L_u$ defined in \eqref{eq:TWN_formal_linearisation} 
  extends to a bounded linear operator $L_{u,\eta} : W^{1,2}_\eta(\R,\R^d) \to L^2_\eta(\R,\R^d)$.
\end{mylemma}
\begin{myproof}
  The term $\partial_x$ maps $W^{1,2}_\eta(\R,\R^d)$ continuously into $L^2_\eta(\R,\R^d)$.
  The exponential weight factors through the multiplication operators $D[ \nabla_g F ](u(\cdot))$ and $D[ \nabla_g S^T ](u(\cdot))$, 
  so that these map $W^{1,2}_\eta(\R,\R^d)$ into $L^2_\eta(\R,\R^d)$ and $L^2_\eta(\R,\Mat_{d \times D}(\R))$, respectively.
 
  Left to verify is the mapping property of the last two terms appearing in \eqref{eq:TWN_formal_linearisation}.
  This is slightly more involved, as these terms depend nonlocally on $v$.
  Using the triangle inequality, for any $\eta' \in \R$ we have
  \[
  \eta' |x| \leq \eta' |x-y| + |\eta'| |y|.
  \]
  Hence, for any $w \in W^{1,2}_{\eta'}(\R,\R^D)$ we find
  \begin{equation}
    \label{eq:TWN_nonloc_exponential_estimate}
      \begin{split}
    e^{\eta' |x|} | \calN[w](x) | &\leq \int_\R \bigg( e^{|\eta'| |y|} \| N(y) \| \bigg) \bigg( e^{\eta' |x - y|} |w(x-y)| \bigg) \d y \\
    &\quad + \sum_{j \in \Z} \bigg( e^{|\eta'| |\lambda_j|} \| N_j \| \bigg) \bigg( e^{\eta' |x - \lambda_j| } |w(x - \lambda_j)| \bigg).
  \end{split}
  \end{equation}
  Let $- \eta_0 \leq \eta \leq \eta_0$, set $\eta' = \eta / 2$, and apply Young's convolution inequality to \eqref{eq:TWN_nonloc_exponential_estimate} to obtain
  \[
  \| \calN[w] \|_{L^2_\eta(\R,\R^D)} \leq C \| w \|_{L^2_{\eta}(\R,\R^D)}
  \]
where
\[
C =  \| N \|_{L^1_{|\eta|/2}(\R,\Mat_{D \times D}(\R))} + \| (N_j)_n \|_{\ell^1_{|\eta|/2}(\Z,\Mat_{D \times D}(\R))}.
\]
  The decay rate imposed on $N(\cdot)$ and $(N_j)_j$ in Hypothesis \ref{hyp:calN} ensures the constant $C$ is finite.
  With $w = DS(u) v$, we find
  \[
  \| \nabla_g S(u)^T \calN[DS(u) v] \|_{L^2_\eta(\R,\R^d)} \leq C' \| v \|_{L^2_{\eta}(\R,\R^d)}
  \]
  for some $C' \geq 0$ depending continuously on $u \in L^\infty(\R,\R^d)$.
  
  To deal with the final term in \eqref{eq:TWN_formal_linearisation}, we note that, with $\theta = (\ell,z_1,z_2)$,
  \[
  \rst{ [ D \psi_\theta(u) v ] }{ (\xi_1, \xi_2) } = \rst{ [ D \sigma_{\ell,z_1}(u) v ] }{ \xi_1 } \rst{ \sigma_{\ell,z_2}(u) }{ \xi_2 } + \rst{ \sigma_{\ell,z_1}(u) }{ \xi_1 } \rst{ [ D \sigma_{\ell,z_2}(u) v ] }{ \xi_2 }.
  \]
  By definition we have $0 \leq \sigma_{\ell,z_i}(u) \leq 1$, which allows us to estimate
  \[
  | \rst{ \alpha_\theta *(  D\psi_\theta(u)v )}{(x,x)} | \leq \rst{ \beta_{1,\theta} * | D \sigma_{\ell,z_1}(u) v |}{x} + \rst{ \beta_{2,\theta} * | D \sigma_{\ell,z_2}(u) v |}{x}, 
\]
where
\[
\beta_{1,\theta}(r) := \int_\R  | \alpha_\theta( r,s ) | \d s, \qquad \beta_{2,\theta}(r) := \int_R | \alpha_\theta( s , r ) | \d s.
  \]
  Incorporating exponential weights into these convolution operators, as we did in \eqref{eq:TWN_nonloc_exponential_estimate}, and applying Young's convolution inequality, we find that
  \begin{align*}
    \| \rst{ \alpha_\theta *(  D\psi_\theta(u)v )}{(x,x)} \|_{L^2_{\eta}(\R,\R^d)} &\leq \| \alpha_\theta \|_{L_{|\eta|/2}^1(\R^2,\R^d)} \| D \sigma_{\ell,z_1}(u) v \|_{L^2_{\eta}(\R,\R)} \\
    &\quad +  \| \alpha_\theta \|_{L_{|\eta|/2}^1(\R^2,\R^d)}  \| D \sigma_{\ell,z_2}(u) v \|_{L^2_{\eta}(\R,\R)},
  \end{align*}
  where we used that $\| \beta_{i,\theta} \|_{L_{|\eta|/2}^1(\R,\R^d)} \leq \| \alpha_\theta \|_{L_{|\eta|/2}^1(\R^2,\R^d)}$ for $i \in \{1,2\}$.
  Recalling the definition of $\sigma_{\rho,z_i}$, we find that there exists $C_R \geq 0$, depending on $R := \| u \|_{L^\infty(\R,\R^d)}$, so that
  \[
  \| \rst{ \alpha_\theta *(  D\psi_\theta(u)v )}{(x,x)} \|_{L^2_{\eta}(\R,\R^d)}  \leq 2 C_R \ell \| \alpha_\theta \|_{L_{|\eta|/2}^1(\R^2,\R^d)} \| v \|_{L^2_\eta(\R,\R^d)}
  \]
  Hence
  \[
  \left\| \sum_{\theta\in \Theta_h} \rst{ \alpha_\theta *\big(  D\psi_\theta(u)v \big)}{(x,x)} \right\|_{L^2_\eta(\R,\R^d)} \leq 2 C \| \Psi \|_{\bPsi_h} \| v \|_{L^2_\eta(\R,\R^d)}.
  \]
  This establishes the bounded mapping of the nonlocal parts appearing in \eqref{eq:TWN_formal_linearisation}.
\end{myproof}

Let $m_\eta : W^{k,2}_\eta(\R,\R^d) \to W^{k,2}(\R,\R^d)$ denote the multiplication operator $u \mapsto e^{\eta \langle \cdot \rangle} u$, where $\langle x \rangle := \sqrt{ 1 + |x|^2 }$.
This defines a linear isomorphism, and as such we can relate $L_{u,\eta}$ to the conjugate operator
\[
L_u(\eta) := m_\eta \circ L_{u,\eta} \circ m_{- \eta} : W^{1,2}(\R,\R^d) \to L^2(\R,\R^d).
\]
The following is then a straighforward consequence of the exponential localisation of the convolution terms appearing in \eqref{eq:TWN_formal_linearisation}.
\begin{mylemma}
  The map 
  \[
  L^\infty(\R,\R^d) \times (-\eta_0,\eta_0) \to \calL(W^{1,2}(\R,\R^d) , L^2(\R,\R^d)), \qquad (u,\eta) \mapsto L_u(\eta)
  \]
  is continuous.
  Furthermore, for fixed $u$, the operator $L_u(\eta)$ depends analytically on $\eta$.
\end{mylemma}
\begin{myproof}
  Continuous dependence of $L_{u,\eta}$ on $u$ follows with the aid of Young's inequality, using estimates identical as those in the proofs of 
  Lemma \ref{lemma:TWN_Phi_path_space}, Lemma \ref{lemma:Psi_L2_smooth}, and Lemma \ref{lemma:Psi_L2_map}.
  Since $L_u(\eta)$ is obtained by conjugating $L_{u,\eta}$ with a map which is independent of $u$, the continuous dependence of $L_u(\eta)$ on $u$ readily follows.
  As for analytic dependence on $\eta$, observe that
  \[
  \rst{ L_{u,-\eta_0} }{ W^{1,2}_{\eta}(\R,\R^d) } = L_{u,\eta}, \qquad \text{for any} \quad \eta \in (-\eta_0,\eta_0),
  \]
  so that
  \[
  L_u(\eta) = m_{\eta} \circ L_{u,-\eta_0} \circ m_{-\eta}.
  \]
  Now choose arbitrary $v \in W^{1,2}(\R,\R^d)$, $w \in L^2(\R,\R^d)$, and $\eta' \in (-\eta_0,\eta_0)$.
  Then note that
  \begin{align*}
    \langle &L_u(\eta) v , w \rangle_{L^2(\R,\R^d)} = \int_\R m_{\eta} \big( L_{u,-\eta_0}[ m_{-\eta} v ](x) \big) \cdot w(x) \d x \\
    &= \int_\R \sum_{n_1=0}^\infty \frac{\langle x \rangle^{n_1} }{n_1!} ( \eta - \eta' )^{n_1} 
      \left( L_{u,-\eta_0}\left[ \sum_{n_2=0}^\infty \frac{ - \langle \cdot \rangle^{n_2} }{n_2!} ( \eta - \eta' )^{n_2} v \right](x) \right) \cdot w(x)  \d x \\
    &=\sum_{(n_1,n_2)\in \N^2} \frac{ 1 }{ n_1! n_2! } (\eta - \eta')^{n_1 + n_2} \int_\R \langle x \rangle^{n_1} \big( L_{u,-\eta_0}\left[ - \langle \cdot \rangle^{n_2} v \right](x) \big) \cdot w(x) \d x,
  \end{align*}
  hence $L_u(\eta)$ depends analytically on $\eta$ in the weak operator topology.
  Analyticity in the operator norm then follows from standard results in analytic perturbation theory \cite{kato2012perturbation}.
\end{myproof}

Now suppose that $u \in L^\infty(\R,\R^d)$ satisfies
\[
\big( u(x) , u'(x) \big) \to \big( z_\pm,0 \big), \qquad \text{as} \quad x\to\pm\infty,
\]
where $z_-$, $z_+$ are critical points of $h$.
Then for any $\epsilon$ there exists $u_* \in \calP(z_-,z_+)$ such that $\| u - u_* \|_{L^\infty(\R,\R^d)} < \epsilon$.
In light of Theorem \ref{thm:TWN_Fredholm} the operator $L_{u_*}(0)$ is Fredholm,
hence by continuity of the Fredholm index, $L_u(0)$ is also Fredholm.
As $L_{u,\eta}$ is conjugate to $L_u(\eta)$, and the latter depends continuous on $\eta$, 
we find that for small values of $\eta$ the operator $L_{u,\eta}$ belongs to the same Fredholm region as $L_{u}(0)$.
We let $\calR_u \subset \R$ denote this Fredholm region, i.e., 
it is the connected neighbourhood around $0$ consisting of values $\eta$ for which the Fredholm index of $L_{u,\eta}$ is constant.

In fact we can show that the kernel $\ker L_{u,\eta}$ is independent of $\eta \in \calR_u$.
To make this statement precise, note that for any $\eta_1 > \eta_2$ we have the natural inclusion map
\[
\iota_{\eta_1,\eta_2} : W^{1,2}_{\eta_1}(\R,\R^d) \hookrightarrow W^{1,2}_{\eta_2}(\R,\R^d).
\]
\begin{mylemma}
\label{lemma:TWN_linear_exponential}
Let $\eta_1 , \eta_2 \in \calR_u$, and suppose $\eta_1 > \eta_2$. 
Then the inclusion map $\iota_{\eta_1,\eta_2}$ defines a bijection from $\ker L_{u,\eta_1}$ onto $\ker L_{u,\eta_2}$.
\end{mylemma}
\begin{myproof}
  First remark that $\iota_{\eta_1,\eta_2}$ is injective from $\ker L_{u,\eta_1}$ into $\ker L_{u,\eta_2}$.
Since $\calR_u$ is connected, it thus suffices to show that $\dim \ker L_{u,\eta}$ is locally constant.
Essentially this now follows from analytic dependence of $\ker L_{u,\eta}$ on $\eta$.

To obtain such an analytic description of $\ker L_{u,\eta}$, we make use of a Lyapunov--Schmidt reduction on the conjugate operator $L_u(\eta)$.
Pick any $\eta' \in \calR_u$, and let $P$ and $Q$ be orthogonal projections only $\ker L_u(\eta')$ and $\coker L_u(\eta')$, respectively.
Letting $v_0 := P v$ and $v_h := (1-P) v$, we find that solving $L_u(\eta) v = 0$ for $\eta$ near $\eta'$ is equivalent to solving
\[
Q L_u(\eta)( v_0 + v_h ) = 0, \qquad (1-Q) L_u(\eta)( v_0 + v_h ) = 0.
\]
The second equation allows us to solve $v_h = A_u(\eta) v_0$, where $A_u(\eta)$ is analytic and $A_u(\eta') = 0$.
Substituting into the first equation, we find that
\[
\ker L_u(\eta) = \set{ v_0 + A_u(\eta) v_0 }{ v_0 \in \ker B_u(\eta) },
\]
where $B_u(\eta) = Q L_u(\eta)( 1 + A(\eta) )$ is a linear map from $\ker L_u(\eta')$ into $\coker L_u(\eta')$.
This shows that
\[
\dim \ker L_u(\eta) \leq \dim \ker L_u(\eta'), \qquad \eta \text{ near } \eta'.
\]
However, by the injectivity of the inclusion $\iota_{\eta',\eta}$ for $\eta \leq \eta'$ we have
\[
\dim \ker L_u(\eta) \geq \dim \ker L_u(\eta'), \qquad \eta \leq \eta'.
\]
Hence $\dim \ker L_u(\eta) = \dim \ker L_u(\eta')$ for $\eta \leq \eta'$.
This can only occur if $B_u(\eta) = 0$ for $\eta \leq \eta'$.
By analyticity we find that $B_u(\eta) = 0$ in a neighbourhood of $\eta'$, hence $\dim \ker L_u(\eta)$ is constant in a neighbourhood of $\eta'$.
Equivalently, $\dim \ker L_{u,\eta}$ is constant in a neighbourhood of $\eta'$.
As $\eta' \in \calR_u$ was chosen arbitrary the claim follows.
\end{myproof}

Since $\ker L_{u,\eta_1}$ and $\ker L_{u,\eta_2}$ are finite dimensional, the inclusion map between the kernels $\iota_{\eta_1,\eta_2} : \ker L_{u,\eta_1} \hookrightarrow \ker L_{u,\eta_2}$
is bijective, with bounded linear inverse.
The norm of the inverse of this map depends, in principle, on the choices of $u$, $\eta_1$, and $\eta_2$.
We now set out to further analyse this dependence.

Denote by $V_u \subset W^{1,\infty}(\R,\R^d)$ the kernel of $L_u$ considered as a map on $W^{1,\infty}(\R,\R^d)$, i.e.,
\[
V_u := \set{ v \in W^{1,\infty}(\R,\R^d) }{ L_u v = 0}.
\]
Given $\eta \in \calR_u$ with $\eta > 0$, define the inclusion $W_{\eta}(u) : V_u \hookrightarrow L^\infty_{\eta}(\R,\R^d)$
by first choosing $\eta' \in \calR_u$ with $\eta' < 0$, and then considering the commuting diagram
\[
\begin{tikzcd}
  V_u
    \arrow[rrr, hook, bend left=20, start anchor={[yshift=1ex]}, end anchor={[yshift=-0.5ex]}, "W_{\eta}(u)"]
    \arrow[r, hook]
  &
  \ker L_{u,\eta'} 
    \arrow[r, "\iota_{\eta,\eta'}^{-1}"]
  &
  \ker L_{u,\eta} 
    \arrow[r]
  &
  L^\infty_{\eta}(\R,R^d).
\end{tikzcd}
\]
The final arrow on the bottom row is induced by the Sobolev embedding $W^{1,2}(\R,\R^d) \hookrightarrow L^\infty(\R,\R^d)$
combined with the isomorphism of weighted Sobolev spaces induced by the multiplication map $m_{\eta}$.
\begin{mylemma}
  The map $W_{\eta}(u)$ is well-defined, i.e., independent of $\eta'$.
\end{mylemma}
\begin{myproof}
  Let $\eta' , \eta'' \in \calR_u$ with $\eta'' < \eta' < 0$.
  A diagram chase then shows the following diagram commutes
  \[
  \begin{tikzcd}
    & \ker L_{u,\eta'} \arrow[dr, "\iota_{\eta,\eta'}^{-1}"] \arrow[dd, "\iota_{\eta',\eta''}"] & \\
   V_u \arrow[ur, hook] \arrow[dr, hook] & & \ker L_{u,\eta} \\
    & \ker L_{u,\eta''} \arrow[ur, "\iota_{\eta,\eta''}^{-1}"'] &
  \end{tikzcd}
  \]
  from which the independence of $W_\eta(u)$ from $\eta'$ follows.
\end{myproof}

We now define $w_{\eta} : W^{1,\infty}(\R,\R^d) \to \R$ by equating $w_{\eta}(u)$ to the operator norm of the linear map $W_{\eta}(u)$,
and investigate the dependence of $w_{\eta}(u)$ on $u$.
This is easier than to study the dependence of $W_{\eta}(u)$ on $u$ directly, since the domain of definition $V_u$ of $W_{\eta}(u)$ depends on $u$,
possibly in a discontinuous manner, as the dimension of $V_u$ may vary with $u$.
Suppose $(u_n)_n \subset L^\infty(\R,\R^d)$ converges uniformly towards $u_\infty$, with for each $n$
\[
\big( u_n(x) , u_n(x) \big) \to \big( z_\pm,0 \big), \qquad \text{as} \quad x\to\pm\infty,
\]
where $z_-$, $z_+$ are critical points of $h$.
Let $\eta\in \calR_{u_\infty}$ with $\eta > 0$.
By continuity of the Fredholm index, we then have $\eta \in \calR_{u_n}$ for sufficiently large $n$.
Thus for each such $n$, we have associated a map $w_{\eta}(u_n)$.
This map turns out to be upper semicontinuous.
\begin{mylemma}
\label{lemma:Weta_upper_semicontinuous}
  Let $(u_n)_n$ and $u_\infty$ be as described above.
  Then, for each $\eta \in \calR_{u_\infty}$ with $\eta > 0$, we have
  \[
  \limsup_{n\to\infty} w_{\eta}(u_n) \leq w_{\eta}(u_\infty).
  \]
\end{mylemma}
\begin{myproof}
  The norms of the maps
  \[
  V_u \hookrightarrow \ker L_{u,\eta'}, \qquad \ker_{u,\eta} \to L^\infty_\eta(\R,\R^d)
  \]
  are both independent of $u$,
  so that upper semicontinuity of the norm of $W_\eta(u)$ is equivalent to 
  upper semicontinuity of the norm of $\iota_{\eta,\eta'}^{-1} : \ker L_{u,\eta'} \to \ker L_{u,\eta}$.
  This in turn is equivalent to upper semicontinuity of the norm of the linear map $f_u : \ker L_u(\eta') \to \ker L_u(\eta)$
  defined by $f_u := m_\eta \circ \iota_{\eta,\eta'}^{-1} \circ m_{-\eta'}$.
  Using a Lyapunov--Schmidt reduction around $u = u_\infty$, $\eta = 0$, with notation as in the proof of Lemma \ref{lemma:TWN_linear_exponential}, 
  we obtain maps 
  \[
    A_u(\eta) , B_u(\eta) : \ker L_{u_\infty}(0) \to \coker L_{u_\infty}(0)
    \]
  and orthogonal projections $P$ and $Q$ onto $\ker L_{u_\infty}(0)$ and $\coker L_{u_\infty}(0)$, respectively,
  such that the following diagram commutes
  \[
  \begin{tikzcd}
    & X_{u}(\eta') \arrow[rr, dashrightarrow, "\phi_{u}"] & & X_{u}(\eta) & \\
    \ker L_{u}(\eta') \arrow[ur, hook] \arrow[rrrr, bend left=60, "f_{u}", "\iso"'] && 
    \ker L_{u_\infty}(0) \arrow[ul, "\id + A_{u}(\eta')", shift left] \arrow[ul, "P"', leftarrow, shift right] \arrow[ur, "\id + A_{u}(\eta)"', shift right] \arrow[ur, "P", leftarrow, shift left] && \ker L_{u}(\eta) \arrow[ul, hook'] \\
    & \ker B_{u}(\eta') \arrow[ul, "\id + A_{u}(\eta')", shift left] \arrow[ul, "P"', leftarrow, shift right] \arrow[rr, equal] \arrow[ur, hook] 
    && \ker B_{u}(\eta) \arrow[ur, "\id + A_{u}(\eta)"', shift right] \arrow[ur, "P", leftarrow, shift left] \arrow[ul, hook'] &
  \end{tikzcd}
  \]
  where
  \[
  X_u(\eta) = \img{ \big( \id + A_{u}(\eta) \big) }. 
  \]
  Here all spaces in the diagram are to be considered subspaces of $W^{1,2}(\R,\R^d)$.
  The map $\phi_u$ is defined by this commuting diagram.
  The spaces $X_u(\eta)$ and the norm of $\phi_{u}$ depend continuously on $u$.
  Note that $f_{u_\infty} = \phi_{u_\infty}$.
  Since $\ker L_{u}(\eta')$ can be a proper subspace of $X_u(\eta')$, we have $\| f_u \| \leq \| \phi_u \|$, with possibly a strict inequality.
  It follows that
  \[
  \limsup_{n\to\infty} \| f_{u_n} \| \leq \limsup_{n\to\infty} \| \phi_{u_n} \| = \| \phi_{u_\infty} \| = \| f_{u_\infty} \|.
  \]
  The upper semicontinuity of $w_\eta$ follows from this.
\end{myproof}

We are now prepared to prove exponential decay of bounded solutions to \eqref{eq:TWN}.
\begin{mythm}
\label{thm:TWN_nonlinear_decay_rates}
Let $z_-$, $z_+$ be critical points of $h$, 
and suppose $u$ is a solution of \eqref{eq:TWN} with $\big( u(x) , u'(x) \big) \to \big( z_\pm , 0 \big)$ as $x\to\pm\infty$.
  Then for any $\eta \in \calR_u$ with $\eta > 0$, there exists $C > 0$ such that
  \begin{align*}
      | u(x) - z_- | + | u'(x) | \leq C e^{-\eta |x|}, & \qquad x < 0, \\
    |u(x) - z_+ | + |u'(x) | \leq C e^{-\eta|x|}, & \qquad x \geq 0.
  \end{align*}
  Furthermore, if $(u_n)_n$ is a sequence of solutions of \eqref{eq:TWN}, each of which satisfying $\big( u_n(x) , u_n'(x) \big) \to \big( z_\pm , 0 \big)$ as $x\to\pm\infty$,
  and the sequence $(u_n)_n$ is convergent in $L^\infty(\R,\R^d)$,
  then the constant $C$ can be chosen independent of $n$.
\end{mythm}
\begin{myproof}
  By Theorem \ref{thm:TWN_weak_compactness}, the derivative $u'$ is bounded in $W^{1,\infty}(\R,\R^d)$.
  Hence $u'$ lies in the domain of definition of $W_\eta(u)$, and we find that
  \[
  \| u' \|_{L^\infty_\eta(\R,\R^d)} = \| W_\eta(u)[u'] \|_{L^\infty_\eta(\R,\R^d)} \leq w_\eta(u) \| u' \|_{W^{1,\infty}(\R,\R^d)}.
  \]
  Hence $u' \in W^{1,2}_{\eta_2}(\R,\R^d)$, so that $u' \in \ker L_{u,\eta_2}$.
  By Lemma \ref{lemma:TWN_linear_exponential} we conclude that $u' \in \ker L_{u,\eta_1}$,
  This establishes that
  \[
  | u'(x) | \leq e^{-\eta |x|} w_\eta(u) \| u' \|_{W^{1,\infty}(\R,\R^d)}.
  \]
  In light of this decay estimate the function $u'$ is integrable over $\R$, hence
  \[
  u(x) - z_- = \int_{-\infty}^x u'(s) \d s.
  \]
  Consequently, for $x < 0$ we find
  \begin{align*}
      |u(x) - z_-| &\leq \int_{-\infty}^x |u'(s)| \d s \leq \bigg( \sup_{x\in\R} e^{\eta |x|} | u'(x) | \bigg)  \int_{-\infty}^x e^{-\eta|s|} \d s \\
  &\leq \eta^{-1} e^{-\eta |x|} w_\eta(u) \| u' \|_{W^{1,\infty}(\R,\R^d)},
  \end{align*}
  and similarly $|u(x) - z_+| \leq \eta^{-1} e^{-\eta |x|} w_\eta(u) \|u'\|_{W^{1,\infty}(\R,\R^d)}$ for $x \geq 0$.
  This establishes that
  \begin{align*}
      | u(x) - z_- | + | u'(x) | \leq C_\eta(u) e^{-\eta |x|}, & \qquad x < 0, \\
    |u(x) - z_+ | + |u'(x) | \leq C_\eta(u) e^{-\eta|x|}, & \qquad x \geq 0,
  \end{align*}
  where $C_\eta(u) = ( \eta^{-1} + 1 ) w_\eta(u) \| u' \|_{W^{1,\infty}(\R,\R^d)}$.
  This proves the theorem for a fixed $u$.
  The fact that $C$ can be chosen independent of $n$ for a convergent sequence $(u_n)_n$ 
  follows from the upper semicontinuity of $w_\eta(u)$ as established in Lemma \ref{lemma:Weta_upper_semicontinuous}.
\end{myproof}

\subsection{A kinematic estimate on $\Psi(u)$}
\label{subsec:sojourn}

The aim of this section is to estimate the $L^2$ norm of $\Psi(u)$ in terms of the kinetic energy of $u$, where $u$ is a $\delta$-approximate solution of \eqref{eq:TWN}.
Such a kinematic estimate of $\Psi(u)$ turns out to be fundamental 
in order to ensure persistence of uniform kinetic energy bounds and gradient-like behaviour of \eqref{eq:TWN} under small perturbations $\Psi$, discussed in Section \ref{sec:TWN_heteroclinics}
and Section \ref{sec:TWN_generic}.
Looking back at the proof of Lemma \ref{lemma:Psi_L2_map}, we see that we are naturally led to study, 
for given $\rho > 0$ and a $\delta$-approximate solution $u$ of \eqref{eq:TWN}, sets of the type
\[
  A_\rho(u) := \set{ x \in \R }{ \inf_{z \in \crit(h)} | u(x) - z | \geq \rho }.
\]
Interpreting $x$ as a time-like variable, 
the set $A_\rho(u)$ consists of the ``times'' at which $u$ sojourns a distance larger than $\rho$ away from the constant solutions of \eqref{eq:TWN}.
Hence we refer to $A_\rho(u)$ as a sojourn set.
In this section we will derive estimates on the sizes of the sojourn sets in terms of kinetic energy.
This then allows us to obtain $L^2(\R,\R^d)$ bounds on $\Psi(u)$ in terms of kinetic energy.

\subsubsection{Breaking of solutions}

Let us now introduce a notion of noncompactness of solutions.

\begin{mydef}[Breaking sequences]
\label{def:TWN_breaking}
    Let $(u_n)_n \subset C^1_{\text{loc}}(\R,\R^d)$. 
  We say the sequence $(u_n)_n$ breaks along $(\tau^0_n)_n$ towards $z_0$ if the following holds:
  \begin{enumerate}
  \item the translated sequence $(\tau^0_n \cdot u_n)_n$ converges to a constant function $z_0$, and
  \item there exists sequences $(\tau^-_n)_n , (\tau^+_n)_n \subset \R$, with $\tau^-_n < \tau^0_n < \tau^+_n$,
  so that the sequences $(\tau^\pm_n \cdot u_n)_n$ converge over a subsequence towards $v_\pm$, with $\calE_{\text{kin}}(v_\pm) > 0$, and
  \item $\big( v_-(x) , v_-'(x) \big) \to \big( z_0 , 0 \big)$ as $x \to \infty$, whilst $\big( v_+(x) , v_+'(x) \big) \to \big( z_0 , 0 \big)$ as $x \to -\infty$.
  \end{enumerate}
\end{mydef}

The breaking of solutions is of importance in the construction of Conley-Floer homology, 
as it is precisely this type of noncompactness which is encoded in the boundary operator of the chain complex.
Usually the breaking of orbits is studied in the context of gradient-like systems.
This we will do towards the end of Section \ref{sec:TWN_heteroclinics}.
Here, we study the breaking of orbits in a nongeneric setup without exploiting gradient-like dichotomies.

Given a measurable subset $B \subseteq \R$, define the convex hull of $B$ by
\[
\conv( B ) := \bigcap_{\substack{(a,b) \subseteq \R \\ B \subseteq (a,b)}} (a,b). 
\]
We now present a result relating the breaking of orbits with the sets $A_\rho(u_n)$.
\begin{mylemma}
\label{lemma:Bn_diameter_blowup_breaking}
Let $R \geq 0$ and suppose $\Psi$ is an $R$-admissible perturbation.
  Let $(u_n)_n \subset C^1_{\text{loc}}(\R,\R^d)$ be a Palais--Smale sequence.
  For each $n \in \N$, let $B_n \subseteq A_\rho(u_n)$ be a measurable subset.
  Assume that
  \[
    \sup_n \| u_n \|_{L^\infty(\R,\R^d)} \leq R, \qquad \sup_n \calE_{\text{kin}}(u_n | \conv(B_n)) < \infty.
  \]
  Suppose $\lim_{n\to\infty} \diam B_n = \infty$.
  Then there exists $(\tau^0_n)_n \subset \R$, with $\tau^0_n \in \conv( B_n )$,
  such that $(u_n)_n$ breaks along a subsequence of $(\tau^0_n)_n$.
\end{mylemma}
\begin{myproof}
  Let $B_n \subseteq A_\rho(u_n)$ with $\lim_{n\to\infty} \diam B_n = \infty$.
  Denote by
  \[
  \tau^L_n := \inf B_n, \qquad \tau^R_n := \sup B_n,
  \]
  so that $\cl{\conv(B_n)} = [\tau^L_n , \tau^R_n]$ and $\tau^R_n - \tau^L_n \to \infty$ as $n\to\infty$.

  First assume there exists a sequence $(\tau^0_n)_n \subset \R$ such that $\tau^L_n < \tau^0_n < \tau^R_n$
  and the translated sequence $(\tau^0_n \cdot u_n)_n$ converges over a subsequence to a constant solution $z_0$ of \eqref{eq:TWN}.
  We will then construct sequences $(\tau^\pm_n)_n$ and functions $v_\pm$ satisfy Definition \ref{def:TWN_breaking}, thus showing that $(u_n)_n$ breaks along a subsequence of $(\tau^0_n)_n$.
  Choose
  \[
  0 < \delta < \min\left\{ \rho , \inf_{\substack{ z_-, z_+ \in \crit(h) \\  |z_-|, |z_+| \leq R}} | z_- - z_+ | \right\}.
  \]
  Such a $\delta$ exists by hyperbolicity of the critical points of $h$.
  Then let
  \[
  \tau^-_n := \sup\set{ \tau \in A_\delta(u_{n}) }{ \tau < \tau^0_n }, \qquad \tau^+_n := \inf\set{ \tau \in A_\delta(u_n) }{ \tau^0_n < \tau }.
  \]
  Since
  \begin{align*}
    A_\delta(u_{n}) \cap ( \tau^L_n , \tau^0_{n} ) &\supset A_\rho(u_{n}) \cap ( \tau^L_n , \tau^0_{n} ) \neq \emptyset, \\
    A_\delta(u_{n}) \cap ( \tau^0_n , \tau^R_{n} ) &\supset A_\rho(u_{n}) \cap ( \tau^0_n , \tau^R_{n} ) \neq \emptyset,
  \end{align*}
  the supremum and infimum in the definition of $\tau^-_n$ and $\tau^+_n$, respectively, are finite.
  Recall that, per assumption, $(\tau^0_n \cdot u_n)_n$ converges to $z_0$ over a subsequence.
  By Theorem \ref{thm:TWN_weak_compactness} the translated sequences $(\tau^\pm_n \cdot u_n)_n$
  converges over a further subsequence $(n_k)_k$ towards solutions $v_\pm$ of \eqref{eq:TWN}.

  Notice that
  \[
  \inf_{z\in\crit(h)} | v_\pm(0) - z | \geq \limsup_{k\to\infty} \inf_{z\in\crit(h)} | u_{n_k}(\tau^\pm_{n_k}) - z | \geq \delta > 0,
  \]
  hence $\calE_{\text{kin}}(v_\pm) > 0$.
  Also note that $\lim_{k\to\infty} | \tau^\pm_{n_k} - \tau^0_{n_k} | = \infty$.
  Indeed, if this where not the case, then at least one of the sequences $(\tau^-_{n_k} \cdot u_{n_k})_k$ or $(\tau^+_{n_k} \cdot u_{n_k})_k$ would accumulate (over a subsequence) onto $z_0$ as $k\to\infty$,
  so that $v_- = z_0$ or $v_+ = z_0$, contradicting $\calE_{\text{kin}}(v_\pm) > 0$.

  Since $\tau^0_{n_k} - \tau^-_{n_k} \to \infty$ as $k \to \infty$ and $(\tau^-_{n_k} , \tau^0_{n_k}) \subset \conv(B_{n_k})$, using Fatou's lemma we obtain the estimate
  \begin{align*}
    \calE_{\text{kin}}( v_- | (0,\infty) ) &= \frac 1 2 \int_\R \lim_{k\to\infty} \1_{(0,\tau^0_{n_k} - \tau^-_{n_k})} g_{\tau^-_{n_k} \cdot u_{n_k}}( \tau^-_{n_k} \cdot u_{n_k} , \tau^-_{n_k} \cdot u_{n_k} ) \d x \\
    &\leq \liminf_{k\to\infty} \calE_{\text{kin}}( \tau^-_{n_k} \cdot u_{n_k} | (0,\tau^0_{n_k} - \tau^-_{n_k}) ) \\
    &= \liminf_{k\to\infty} \calE_{\text{kin}}( u_{n_k} | (\tau^-_{n_k}, \tau^0_{n_k}) ) \\
    &\leq \liminf_{k\to\infty} \calE_{\text{kin}}( u_{n_k} | \conv(B_{n_k})).
  \end{align*}
  It thus follows from the hypothesis of the lemma that $\calE_{\text{kin}}( v_- | (0,\infty) ) < \infty$, and similarly, $\calE_{\text{kin}}( v_+ | (-\infty,0) ) < \infty$.
  By Lemma \ref{lemma:TWN_limit_sets} $\big( v_-(x) , v_-'(x) \big) \to \big(z_1,0\big)$ as $x\to \infty$ and $\big( v_+(x) , v_+'(x) \big) \to \big(z_2,0\big)$ as $x\to -\infty$,
  where $z_1$ and $z_2$ are constant solutions of \eqref{eq:TWN}.

  Left to verify is that $z_0 = z_1 = z_2$.
  For any $x > 0$, we have $\tau^-_{n_k} + x < \tau^0_{n_k}$ for sufficiently large $k$.
  By definition of $(\tau^-_{n_k})_k$ we then have $\tau^-_{n_k} + x \not\in A_\delta(u_{n_k})$.
  Consequently
  \[
  | v_-(x) - z_0 | = \lim_{k\to\infty} | u_{n_k}(\tau^-_{n_k} + x) - z_0 | \leq \delta.
  \]
  Letting $x\to\infty$, we find that $|z_1 - z_0| \leq \delta$.
  Similarly one shows that $| z_2 - z_0 | \leq \delta$.
  Since $\delta$ is chosen such that the only constant solution of \eqref{eq:TWN} in a $\delta$-neighbourhood of $z_0$ is $z_0$ itself,
  we conclude that $z_0 = z_1 = z_2$.
  This shows that $(\tau^\pm_n)_n$, $(\tau^0_n)_n$, and $v_\pm$ satisfy Definition \ref{def:TWN_breaking}.

  We will now construct the sequence $(\tau^0_n)_n$.
  Define $v_n := \tau^L_n \cdot u_n$.
  By Theorem \ref{thm:TWN_weak_compactness} the sequence $(v_n)_n$ converges over a subsequence (which we again denote by $(v_n)_n$)
  to a solution $v_\infty$ of \eqref{eq:TWN}.
  By the max-min inequality
  \[
  \inf_{z\in\crit(h)} |v_\infty(0) - z| \geq \limsup_{n\to\infty} \inf_{z\in\crit(h)} |u_n(\tau^L_n) - z| \geq \rho,
  \]
  so that $v_\infty$ is not a constant solution of \eqref{eq:TWN}.
  Since $\tau^R_n - \tau^L_n \to \infty$ as $k \to \infty$ and $(\tau^L_n , \tau^R_n) \subset \conv(B_n)$, using Fatou's lemma we obtain the estimate
  \begin{align*}
    \calE_{\text{kin}}( v_\infty | (0,\infty) ) &= \frac 1 2 \int_\R \lim_{n\to\infty} \1_{(0,\tau^R_n - \tau^L_n)} g_{\tau^L_n  \cdot u_{n_k}}( \tau^L_n \cdot u_n , \tau^L_n \cdot u_n ) \d x \\
    &\leq \liminf_{n\to\infty} \calE_{\text{kin}}( \tau^L_n \cdot u_n | (0,\tau^R_n - \tau^L_n) ) \\ 
    &\leq \liminf_{n\to\infty} \calE_{\text{kin}}( u_n | \conv(B_n)).
  \end{align*}
  By hypothesis of the lemma it now follows that $\calE_{\text{kin}}( v_\infty | (0,\infty) ) < \infty$, hence by Lemma \ref{lemma:TWN_limit_sets}
  we find that $\big( v_\infty(x) , v_\infty'(x) \big) \to \big( z_0 , 0 \big)$ as $x \to \infty$, for some $z_0 \in \crit(h)$.

    For given $x\geq 0$, by the $C^1_{\text{loc}}(\R,\R^d)$ convergence of $v_n$ towards $v_\infty$ 
  we may choose a sequence $(\tau_n(x))_n \subset \R$ such that 
  \[
  \tau^L_n < \tau_n(x), \qquad \sup_n | \tau_n(x) - \tau^L_n | < \infty, \qquad \lim_{n\to\infty} u(\tau_n(x)) = v(x).
  \]
  Then since $\tau^R_n - \tau^L_n \to \infty$ as $n\to\infty$, for fixed $x \geq 0$ one has $\tau_n(x) \in \conv(B_n)$ for all sufficiently large $n$.
  Hence, letting $x = k$ where $k \in \N$,
  we may choose a sequence $(n_k)_k \subset \N$ so that $\tau_n(k) \in \conv( B_{n_k} )$ for all $k \in \N$.
  We now define the sequence $(\tau^0_n)_n$ on the subsequence $(n_k)_k$ by setting $\tau^0_{n_k} := \tau_{n_k}(k)$.
  On the complement of the subsequence $(n_k)_k$ we may choose $\tau^0_n$ arbitrary.
  This establishes the sequence $(\tau^0_n)_n$ with
  \[
  \tau^L_n < \tau^0_n < \tau^R_n, \qquad \lim_{n\to\infty} | \tau^0_n - \tau^L_n | = \lim_{n\to\infty} | \tau^0_n - \tau^R_n | = \infty,
  \]
  such that $(\tau^0_n \cdot u_n)_n$ converges over the subsequence $(n_k)_k$ towards the constant solution $z_0$ of \eqref{eq:TWN}.
\end{myproof}

\begin{mylemma}
\label{lemma:B_diameter_bound}
Given $R_1 , R_2 \geq 0$, $\rho > 0$, and an $R_1$-admissible perturbation $\Psi$, there exist $D > 0$, $\delta > 0$ such that the following holds.
  Let $u$ be a $\delta$-approximate solution of \eqref{eq:TWN}.
  Suppose $B \subseteq A_\rho(u)$ is a connected subset, and
  \[
    \| u \|_{L^\infty(\R,\R^d)} \leq R_1, \qquad \calE_{\text{kin}}(u|B) \leq R_2.
  \]
  Then $\diam(B) \leq D$.
\end{mylemma}
\begin{myproof}
  Given $\delta > 0$, let $D_\delta$ denote the supremum of $\diam(B)$, ranging over all sets $B$ as in the hypothesis of the lemma.
  Now, arguing by contradiction, suppose $\limsup_{\delta \to 0} D_\delta = \infty$.
  We then obtain a Palais--Smale sequence $(u_n)_n$, and connected subsets $B_n \subseteq A_\rho(u_n)$,
  such that
  \[
  \sup_n \| u_n \|_{L^\infty(\R,\R^d)} \leq R_1, \qquad \sup_n \calE_{\text{kin}}(u_n | B_n) \leq R_2,
  \]
  but $\lim_{n\to\infty} \diam(B_n) = \infty$.
  Then, by Lemma \ref{lemma:Bn_diameter_blowup_breaking}, we obtain a sequence $(\tau^0_n)_n$ with $\tau^0_n \in B_n$,
  such that $( \tau^0_n \cdot u_n )_n$ converges over a subsequence $(n_k)_k$ towards a constant solution $z_0$ of \eqref{eq:TWN}.
  However,
  \[
  \inf_{z\in\crit(h)} | z_0- z | \geq \limsup_{k\to\infty} \inf_{z\in\crit(h)} | u_{n_k}(\tau^0_{n_k}) - z | \geq \rho > 0,
  \]
  since $\tau^0_n \in A_\rho(u_n)$.
  This contradicts with $z_0$ being a constant solution of \eqref{eq:TWN}.
  We conclude that $\limsup_{\delta \to 0} D_\delta < \infty$, from which the conclusion of the lemma follows.
\end{myproof}

\subsubsection{Minimal energy quanta}
In Lemma \ref{lemma:B_diameter_bound} we saw that bounds on kinetic energy of a $\delta$-approximate $u$ over an interval $B\subseteq A_\rho(u)$ lead to bounds on the diameter of $B$.
In this section we will show, loosely speaking, that given $\rho > 0$ there is a minimal amount of energy $\hbar$ required for any $\delta$-approximate solution of \eqref{eq:TWN}
to sojourn a distance $\rho$ away from the constant solutions of \eqref{eq:TWN}.
This then provides us with ``minimal energy quanta'' into which we can decompose the sojourn sets.
\begin{mylemma}
\label{lemma:kinetic_energy_lower_bound}
  Given $R \geq 0$, $r_1 , r_2 > 0$ with $r_1 \neq r_2$, and an $R$-admissible perturbation $\Psi$,
  there exist $\delta > 0$ and $\epsilon_{r_1,r_2} > 0$ such that the following holds.
  Suppose $u$ is a $\delta$-approximate solution of \eqref{eq:TWN}, such that $\| u \|_{L^\infty(\R,\R^d)} \leq R$.
  Suppose $-\infty \leq a < b \leq \infty$ are such that 
  \[
  \lim_{x \to a} \inf_{z\in\crit(h)} | u(x) - z | = r_1, \qquad \lim_{x \to b} \inf_{z\in\crit(h)} | u(x) - z | = r_2.
  \]
  Then
  \[
  \calE_{\text{kin}}(u | (a,b) ) \geq \epsilon_{r_1,r_2}.
  \]
\end{mylemma}
\begin{myproof}
  First note, if $(a,b) = \R$, it follows from Lemma \ref{lemma:TWN_limit_sets} that $\calE_{\text{kin}}(u|(a,b)) = \infty$.
  Now assume $(a,b) \neq \R$.
  We will prove the energy estimate assuming that $a > - \infty$, which we may do without loss of generality, 
  since in case $a = -\infty$ but $b < \infty$ the result can be applied to $\check{u}(x) := u(-x)$.
  Then, since $\calE_{\text{kin}}(u | (a,b)) = \calE_{\text{kin}}( a \cdot u | (0,b-a))$, we may as well assume $a = 0$.
  Let
  \[
  \Sigma := \set{ (u,s) \in C^1_{\text{loc}}(\R,\R^d) \times [0,\infty] }{
    \begin{array}{c}
     \| u' + \Phi(u) + \Psi(u) \|_{L^\infty(\R,\R^d)} \leq \delta \\
      \| u \|_{L^\infty(\R,\R^d)} \leq R \\
      \inf_{z\in\crit(h)} | u(0) - z | = r_1 \\
      \lim_{x \to s} \inf_{z\in\crit(h)} | u(x) - z | = r_2
    \end{array}
    },
    \]
    where $[0,\infty] = [0,\infty) \cup \{\infty\}$ is the one-point compactification of $[0,\infty)$.
    In light of Theorem \ref{thm:TWN_weak_compactness},
    the set $\Sigma$ is precompact in $C^1_{\text{loc}}(\R,\R^d) \times [0,\infty]$;
    we denote by $\cl{\Sigma}$ the closure of $\Sigma$.
    We will show that $\inf_{(u,s)\in \cl{\Sigma}} \calE_{\text{kin}}(u|(0,s)) > 0$.
    
    The key properties for the remainder of the argument are that, in light of Theorem \ref{thm:TWN_weak_compactness}, 
    $\cl{\Sigma}$ is compact and $\sup_{(u,s) \in \cl{\Sigma}} \| u \|_{W^{1,\infty}(\R,\R^d)} < \infty$.
    Furthermore, given $(u,s) \in \cl{\Sigma}$, we claim that $u' \neq 0$ on some open subset of $(0,s)$.
    To see why, we need to consider two cases, $s < \infty$ and $s=\infty$.
    
    When $s < \infty$, we use the fact that $u$ is the limit in $C^1_{\text{loc}}(\R,\R^d)$ of a sequence $( u_n )_n$, 
    whereas $s$ is the limit in $(0,\infty)$ of a bounded sequence $(s_n)_n$.
    This ensures that $u_n(s_n) \to u(s)$, hence
    \[
    \inf_{z\in\crit(h)} | u(s) - z | = \lim_{n\to\infty} \inf_{z\in\crit(h)} | u_n(s_n) - z | = r_2.
    \]
    Here we used that $u$ is uniformly bounded and $h$ is Morse, so that the infima are in fact minima over finite sets, hence they depend continuously on $u(s)$.
    Since it also holds that $\inf_{z \in \crit(h)} | u(0) - z | = r_1 \neq r_2$, we find that $u(0) \neq u(s)$, so that $u' \neq 0$ on some open subset of $(0,s)$.
    
    If, on the other hand, $s = \infty$, it may happen that $\lim_{x\to s} \inf_{z\in\crit(h)} | u(x) - z | \neq r_2$, or even that this limit does not exist.
    Thus, in this case, we need a different argument to conclude that $u' \neq 0$ on an open subset of $(0,s)$.
    It is at this stage that we specify the constant $\delta > 0$.
    Using the fact that $\Psi$ is $R$-admissible, we may choose $\delta > 0$ sufficiently small such that for any $v\in \R^d$
    \[
    |v| \leq R \quad \text{and} \quad | \nabla_g h(v) + \Psi(v) | \leq \delta \qquad \Longrightarrow \qquad \inf_{z\in\crit(h)} | v - z | < r_1.
    \]
    With this choice of $\delta > 0$ established, we will now show that $u' \neq 0$ on an open subset of $(0,s) = (0,\infty)$.
    We argue by contradiction, assuming $u' = 0$ on $(0,\infty)$.
    By exponential decay of the convolution kernels appearing in the nonlocal operators $\Phi$ and $\Psi$, we have
    \[
    \Phi(u(0)) = \lim_{x\to\infty} \Phi(u)(x), \qquad \Psi(u(0)) = \lim_{x\to\infty} \Psi(u)(x).
    \]
    Since $\nabla_g h(u(0)) = \Phi(u(0))$, and using that $u$ is a $\delta$-approximate solution of \eqref{eq:TWN}, we thus have
    \[
     | \nabla_g h(u(0)) + \Psi(u(0)) | 
      = \lim_{x  \to \infty} | u'(x) + \Phi(u)(x) + \Psi(u)(x)| \leq \delta.
    \]
    But then $\inf_{z \in \crit(h)} | u(0) - z | < r_1$, contradicting the fact that $(u , \infty) \in \cl{\Sigma}$.
    Hence we conclude that $u' \neq 0$ on an open subset of $(0,\infty)$.

    Since $\calE_{\text{kin}}$ is not continuous with respect to the topology of $C^1_{\text{loc}}(\R,\R^d)$, we introduce the regularised energy
\[
  \calE_n : \cl{\Sigma} \to [0,\infty], \qquad \calE_n(u,s) := \frac 1 2 \int_0^s \frac{ g_{u(x)}(u'(x),u'(x)) }{ 1 + |x|^2/n } \d x.
\]
    By continuity and compactness there exists $(u_n,s_n)_n \subset \cl{ \Sigma }$ such that
    \[
    \calE_n(u_n,s_n) = \inf_{(u,s) \in \cl{ \Sigma }} \calE_n(u,s) \qquad \text{for all} \quad n \in \N.
    \]
    Since $u_n' \neq 0$ on some open subset of $(0,s_n)$, we find that $\calE_n(u_n,s_n) > 0$.
    Furthermore, the minimising property of $(u_n,s_n)$ and the definition of $\calE_n$ ensure that
    \[
    \calE_n(u_n,s_n) \leq \calE_n(u_{n+1},s_{n+1}) \leq \calE_{n+1}(u_{n+1},s_{n+1}) \qquad \text{for all} \quad n \in \N.
    \]
    We thus find that
    \[
    \inf_{(u,s)\in \Sigma} \calE_{\text{kin}}(u|(0,s)) \geq \lim_{n\to\infty} \calE_n(u_n,s_n) = \sup_n \calE_n(u_n,s_n) > 0.
    \]
\end{myproof}

\begin{mylemma}
\label{lemma:kinetic_energy_lower_bound_B}
  Given $R \geq 0$, $\rho_1 > \rho_2 > 0$, and an $R$-admissible perturbation $\Psi$, there exists $\delta > 0$ and $\hbar > 0$ such that the following holds.
  Suppose $u$ is a $\delta$-approximate solution of \eqref{eq:TWN}, with $\| u \|_{L^\infty(\R,\R^d)} \leq R$ and $\calE_{\text{kin}}(u) < \infty$.
  Let $B \subseteq A_{\rho/2}(u)$ be a connected component and suppose $B' := B \cap A_\rho(u) \neq \emptyset$.
  Then
  \[
  \calE_{\text{kin}}(u|B \setmin B') \geq \hbar.
  \]
\end{mylemma}
\begin{myproof}
    By Lemma \ref{lemma:TWN_limit_sets_approximate} the limit sets $\alpha(u)$ and $\omega(u)$ consist of $\delta$-approximate constant solutions of \eqref{eq:TWN}.
    Using the fact that $\Psi$ is $R$-admissible, we may choose $\delta > 0$ sufficiently small such that for any $v\in \R^d$
    \[
    |v| \leq R \quad \text{and} \quad | \nabla_g h(v) + \Psi(v) | \leq \delta \qquad \Longrightarrow \qquad \inf_{z\in\crit(h)} | v - z | < \rho_2.
    \]
    With this choice of $\delta$, we thus find
    \[
    \alpha(u) \cap A_{\rho_2}(u) = \emptyset, \qquad \omega(u) \cap A_{\rho_2}(u) = \emptyset.
    \]
  The definition of $B$, $B'$ now ensures there exist $-\infty < a < b \leq a' < b' < \infty$ such that
  \[
  (a,b) \cup (a',b') \subseteq B' \setmin B,
  \]
  and
  \begin{align*}
     \inf_{z\in\crit(h)} |u(b) - z| &= \inf_{z\in\crit(h)} |u(a') - z| = \rho_1, \\
    \inf_{z\in\crit(h)} |u(a) - z| &= \inf_{z\in\crit(h)}  |u(b') - z| = \rho_2.
  \end{align*}
  After decreasing the value of $\delta > 0$, we may apply Lemma \ref{lemma:kinetic_energy_lower_bound}, which yields
  \[
  \calE_{\text{kin}}(u|B\setmin B') \geq \calE_{\text{kin}}(u|(a,b)) + \calE_{\text{kin}}(u|(a',b')) \geq \epsilon_{\rho_2,\rho_1} + \epsilon_{\rho_1,\rho_2},
  \]
  where $\epsilon_{r_1,r_2}$ are as in Lemma \ref{lemma:kinetic_energy_lower_bound}.
  This proves the lemma, with $\hbar = \epsilon_{\rho_2,\rho_1} + \epsilon_{\rho_1,\rho_2}$.
\end{myproof}

\subsubsection{Bounding the volume of sojourn sets}

We are now prepared to derive a fundamental relation between sojourn sets and kinetic energy.

\begin{mythm}
\label{thm:Arho_volume_bound}
  Given $R \geq 0$, $\rho > 0$, and an $R$-admissible perturbation $\Psi$, there exist $C_{R,\rho} > 0$ and $\delta > 0$ such that the following holds.
  Suppose $u$ is a $\delta$-approximate solution of \eqref{eq:TWN}, and $\| u \|_{L^\infty(\R,\R^d)} \leq R$.
  Then
  \[
  \vol( A_\rho(u) ) \leq C_{R,\rho} \calE_{\text{kin}}(u).
  \]
\end{mythm}
\begin{figure}
  \centering
  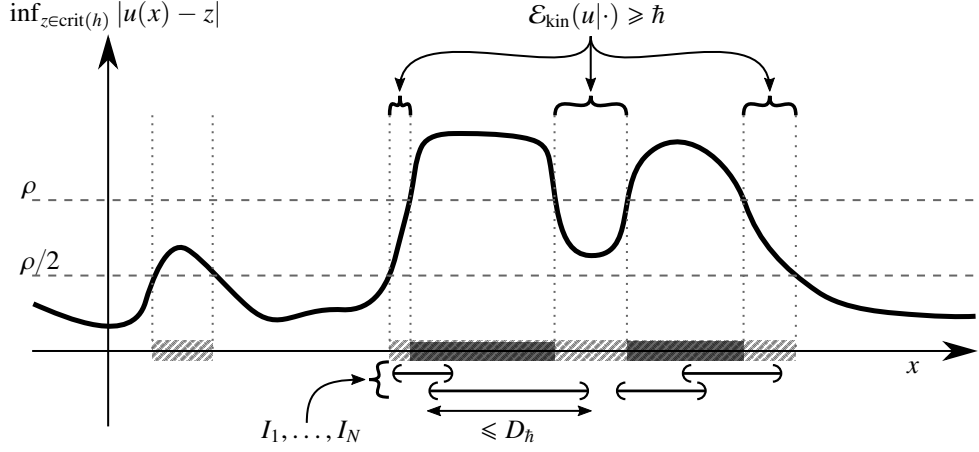
  \caption{Covering the sojourn set by minimal energy quanta $I_1,\dots,I_n$ with $\calE_{\text{kin}}(u | I_k) = \hbar$.}
  \label{fig:TWN_sojourn}
\end{figure}
\begin{myproof}
Without loss of generality, assume $\calE_{\text{kin}}(u) < \infty$.
Let $\delta > 0$ and $\hbar > 0$ be as specified by Lemma \ref{lemma:kinetic_energy_lower_bound_B} with $\rho_1 = \rho$ and $\rho_2 = \rho/2$.
Let $\calB$ consist of the connected components of $A_{\rho/2}(u)$.
Pick $B \in \calB$, suppose $B \cap A_\rho(u) \neq \emptyset$, and let $I_1,\dots,I_n \subset A_{\rho/2}(u)$ be intervals such that
\[
B \cap A_\rho(u) \subset \bigcup_{k=1}^n I_k \subset B, \qquad \calE_{\text{kin}}(u | I_k) = \hbar.
\]
See also Figure \ref{fig:TWN_sojourn}.
With regards to the maximum number $n$ of intervals required for such a covering, we see from Lemma \ref{lemma:kinetic_energy_lower_bound_B} that
\[
n = \left\lceil \frac{1}{\hbar} \calE_{\text{kin}}(u|B \cap A_\rho(u)) \right\rceil \leq \left\lceil \frac{1}{\hbar} \bigg( \calE_{\text{kin}}(u|B) - \hbar \bigg) \right\rceil 
\leq \frac{1}{\hbar} \calE_{\text{kin}}(u|B).
\]
By Lemma \ref{lemma:B_diameter_bound}, applied with $\rho/2$ instead of $\rho$ and $R_2 = \hbar$, there exists a constant $D_\hbar > 0$, independent of $u$, such that
\[
\diam( I_k ) \leq D_\hbar.
\]
Combining these estimates, we obtain
\[
\vol(B \cap A_\rho(u) ) \leq \sum_{k=1}^n \diam(I_k) \leq \frac{D_\hbar}{\hbar} \calE_{\text{kin}}(u|B).
\]
Note that this estimate trivially extends to $B \in \calB$ for which $B \cap A_\rho(u) = \emptyset$.
Finally, since $A_\rho(u) \subset A_{\rho/2}(u)$, we obtain
\[
\vol(A_\rho(u)) = \sum_{ B \in \calB } \vol(B \cap A_\rho(u)) \leq \frac{D_\hbar}{\hbar} \calE_{\text{kin}}(u).
\]
We remark that the essential ingredient which allows for uniformity in this estimate is the existence of the lower bound $\hbar$ for the kinetic energy.
\end{myproof}

\subsubsection{The estimate on  $\Psi(u)$}

Using the bounds on the sojourn set, we now obtain the following kinematic estimate on $\Psi(u)$.
\begin{mythm}
\label{thm:kinematic_estimate_Psi}
  Given $R \geq 0$ and an $R$-admissible perturbation $\Psi_0$, there exist $C_{R,\Psi_0} > 0$ and $\delta > 0$ such that the following holds.
  Suppose $u$ is a $\delta$-approximate solution of \eqref{eq:TWN} with $\Psi = \Psi_0$, and $\| u \|_{L^\infty(\R,\R^d)} \leq R$.
  Then, for any $\Psi \in \bPsi_h$ we have
  \[
  \| \Psi(u) \|_{L^2(\R,\R^d)} \leq C_{R,\Psi_0} \| \Psi \|_{\bPsi_h} \sqrt{ \calE_{\text{kin}}(u) }.
  \]
\end{mythm}
\begin{myproof}
We recall from Lemma \ref{lemma:Psi_L2_map} the estimate
\begin{equation}
  \label{eq:TWN_Psi_sojourn_estimate}
  \begin{split}
      \| \Psi(u) \|_{L^2(\R,\R^d)} &\leq \| \Psi \|_{\bPsi_h} \sqrt{ \vol\left( \set{ x \in \R }{ \inf_{z\in\crit(h)} |u(x) - z| \geq \rho_R } \right) } \\
  &= \| \Psi \|_{\bPsi_h}  \sqrt{ \vol( A_{\rho_R}(u) ) }.
  \end{split}
\end{equation}
This was derived in the context of path spaces.
However, the only point in the proof of Lemma \ref{lemma:Psi_L2_map} where the assumption $u\in \calP(z_-,z_+)$
is used is to ensure the volume of the sojourn set $A_{\rho_R}(u)$ is finite.
Instead, we now assume $u$ is a $\delta$-approximate solution of \eqref{eq:TWN}.
The result then follows by combining estimate \eqref{eq:TWN_Psi_sojourn_estimate} with Theorem \ref{thm:Arho_volume_bound}.
\end{myproof}

\subsection{Isolating trajectory neighbourhoods}
\label{subsec:TWN_isolating_nbhds}

Since equation \eqref{eq:TWN} is not formulated on a phase space, we cannot use the standard definition of isolating neighbourhoods and isolated invariant sets.
Therefore, we introduce here a new notion, that of isolating trajectory neighbourhoods, which isolate bounded solutions without exploiting pointwise phase space considerations.
Results in Section \ref{sec:TWN_Morse_isomorphism} on isolating blocks will show that this is, indeed, an appropriate generalisation of isolating neighbourhoods.


\begin{mydef}[Isolating trajectory neighbourhood]
\label{def:isolating_trajectory_nbhd}
  A subset $E \subset C_{\text{loc}}^1(\R,\R^d)$ is called an isolating trajectory neighbourhood provided that
  \begin{enumerate}
  \item $E$ is translation invariant, and
  \item $E$ is open in the induced topology of $W^{1,\infty}(\R,\R^d)$, and bounded in $L^\infty(\R,\R^d)$, and
  \item the closure of $E$ in $C^1_{\text{loc}}(\R,\R^d)$ coincides with the closure of $E$ in $W^{1,\infty}(\R,\R^d)$, and
  \item there exists $\delta > 0$ such that $\cl E \setmin E$ does not contain any $\delta$-approximate solutions of \eqref{eq:TWN}.
  \end{enumerate}
  We define the invariant set of $E$ by
  \[
  \Inv(E) := \set{ u \in E }{ u' + \Phi(u) + \Psi(u) = 0 }.
  \]
  A set $I$ is called an isolated invariant set if it can be written as $I = \Inv(E)$ for some isolating trajectory neighbourhood $E$.
\end{mydef}

\begin{myremark}
  The definition of an isolating invariant set is rather technical and warrants some further motivation.
The difficulty stems from the fact that the natural space to deal with compactness issues of \eqref{eq:TWN} is $C^1_{\text{loc}}(\R,\R^d)$.
However, open neighbourhoods in $C^1_{\text{loc}}(\R,\R^d)$ are much too large to allow for localisation around bounded solutions;
for this, uniform bounds are much more appropriate.
Thus we need to relate these modes of convergence.
We cannot exploit any pointwise considerations, and so we need to impose a compatibility conditions between the topologies, which is done with property (3).

Sets which do not satisfy property (3) are easily constructed.
Consider, for example, the open annulus 
\[
A = \set{ u \in C^1_{\text{loc}}(\R,\R^d) }{ r_1 < \| u \|_{L^\infty(\R,\R^d)} < r_2 }.
\]
The closure $\cl A$ of $A$ in $W^{1,\infty}(\R,\R^d)$ is the closed annulus, but the $C^1_{\text{loc}}(\R,\R^d)$ closure of $A$ is the closed ball of radius $r_2$.
This is related to a fundamental difference in the characterisation of the boundary components of $\cl A$:
given $u \in \cl A$, to verify whether $\| u \|_{L^\infty(\R,\R^d)} = r_2$ it suffices to prove existence of a single $x_0\in\R$ such that $|u(x_0)| = r_2$;
however, verifying that $\| u \|_{L^\infty(\R,\R^d)} = r_1$ requires inspecting $|u(x)|$ for all values of $x \in \R$.
Property (3) may thus be loosely interpreted as generalising the idea of $\cl E \setmin E$ having a localised, almost pointwise, characterisation.
\end{myremark}

Our first result is a characterisation of isolated invariant sets.
\begin{mylemma}
\label{lemma:TWN_isolated_invset_nbhd}
  Let
  \[
  \Sigma := \set{ u \in C^1_{\text{loc}}(\R,\R^d) }{
    \begin{array}{c}
      u' + \Phi(u) + \Psi(u) = 0 \\
      \| u \|_{L^\infty(\R,\R^d)} < \infty
    \end{array}
    }.
  \]
  Suppose $I \subseteq \Sigma$ is such that:
  \begin{enumerate}
  \item The set $I$ is invariant under translation, $\tau \cdot u \in I$ for any $u \in I$ and $\tau \in \R$.
  \item The set $I$ is isolated and bounded in $\Sigma$ with respect to the topology of $L^\infty(\R,\R^d)$.
  \item The set $I$ is compact in the topology of $C^1_{\text{loc}}(\R,\R^d)$.
  \end{enumerate}
  Then there exists $\epsilon > 0$ so that
  \[
  E := \set{ v \in C^1_{\text{loc}}(\R,\R^d) }{  \exists u \in I \quad \text{such that} \quad \| u - v \|_{L^\infty(\R,\R^d)} < \epsilon }
  \]
  is an isolating trajectory neighbourhood, with $\Inv(E) = I$.
\end{mylemma}
\begin{myproof}
It is clear that properties (1) and (2) in the definition of an isolating trajectory neighbourhood are satisfied.
Isolation of $I$ in $\Sigma$ ensures that with $\epsilon > 0$ sufficiently small, the set $\cl{E} \setmin I$ does not contain any solutions of \eqref{eq:TWN}.
In particular we have $\Inv(E) = I$.
Left to prove is that properties (3) and (4) from Definition \ref{def:isolating_trajectory_nbhd} are satisfied.


To see why property (3) from Definition \ref{def:isolating_trajectory_nbhd} holds, let $(v_n)_n \subset E$ be a sequence converging towards $v_\infty$ 
in the topology of $C^1_{\text{loc}}(\R,\R^d)$.
By definition of $E$ we can find $(u_n)_n \subset I$ such that $\| u_n - v_n \|_{L^\infty(\R,\R^d)} < \epsilon$.
By the compactness assumption on $I$, after choosing a subsequence, we may assume $(u_n)_n$ converges towards $u_\infty \in I$
in the topology of $C^1_{\text{loc}}(\R,\R^d)$.
Thus the sequence $(u_n - v_n)_n$ converges (over a subsequence) towards $u_\infty - v_\infty$, again in the topology of $C^1_{\text{loc}}(\R,\R^d)$.
But then $\| u_\infty - v_\infty \|_{L^\infty(\R,\R^d)} \leq \epsilon$, hence $v_\infty$ belongs to the closure of $E$ in the topology of $W^{1,\infty}(\R,\R^d)$.

Now to prove property (4) from Definition \ref{def:isolating_trajectory_nbhd}, first define
\[
\delta := \inf_{u\in\cl{E} \setmin E} \| v' + \Phi(v) + \Psi(v) \|_{L^\infty(\R,\R^d)},
\]
and let $(v_n)_n \subset \cl{E} \setmin E$ be such that
\begin{equation}
  \label{eq:TWN_isolated_estimate}
  \lim_{n\to\infty} \| v_n' + \Phi(v_n) + \Psi(v_n) \|_{L^\infty(\R,\R^d)} = \delta.
\end{equation}
We will prove that $\delta > 0$, and argue by contradiction, assuming $\delta = 0$.
Then \eqref{eq:TWN_isolated_estimate} tells us that $(v_n)_n$ is a Palais--Smale sequence.
For each $n$, there exists sequences $(x_{n,k})_k \subset \R$ and $(u_{n,k})_k \subset I$ so that 
\[
| u_{n,k}(0) - v_n(x_{n,k}) | \to \epsilon \qquad \text{as} \quad k\to\infty.
\]
Moreover, $| u_{n,k}(x) - v_n(x) | \leq \epsilon$ for any $x \in \R$ and $k \in \N$.
Let $k_n \in \N$ be such that
\[
\epsilon - \frac 1 n \leq | u_{n,k}(0) - v_n(x_{n,k}) | \leq \epsilon \qquad \text{for} \quad k \geq k_n.
\]
Then let $w_n := x_{n,k_n} \cdot v_n$.
Now, by translation invariance of \eqref{eq:TWN_isolated_estimate}, using that $\delta = 0$ per assumption, $(w_n)_n$ is a Palais--Smale sequence as well.
By Theorem \ref{thm:TWN_weak_compactness} there exists a subsequence along which $(w_n)_{n}$ accumulates onto a solution $w_\infty$ of \eqref{eq:TWN}.
Moreover, by compactness of $I$ there exists a subsequence along which $(u_{n,k_n})_n$ accumulates onto $u_\infty \in I$.
By construction of the sequences $(u_{n,k_n})_n$ and $(w_n)_n$, we then have 
\[
| u_\infty(0) - w_\infty(0)| = \epsilon, \qquad \text{and} \quad | u_\infty(x) - w_\infty(x)| \leq \epsilon \quad \text{for all} \quad x \in \R.
\]
Hence we have constructed a solution $w_\infty \in \cl{E} \setmin E$ of \eqref{eq:TWN}.
However, since $I$ is isolated in $\Sigma$, the set $\cl{E} \setmin I$ does not contain any solutions of \eqref{eq:TWN}.
We thus arrived at a contradiction, and must conclude that $\delta > 0$.
\end{myproof}

The downside of this characterisation of isolating trajectory neighbourhoods is that we need prior knowledge about the invariant set $I$.
In applications we usually do not have such knowledge.
In Section \ref{sec:TWN_Morse_isomorphism} we will describe a subclass of isolating trajectory neighbourhoods for which such prior knowledge is not needed.
In the remainder of this section we present two lemmata which demonstrate the relevance of isolating trajectory neighbourhoods in the study of moduli spaces.

\begin{mylemma}
\label{lemma:TWN_invset_cpt}
  Let $E$ be an isolating trajectory neighbourhood.
Then $\Inv(E)$ is compact in $C^1_{\text{loc}}(\R,\R^d)$.
\end{mylemma}
\begin{myproof}
  Let $(u_n)_n \subset \Inv(E)$ be an arbitrary sequence.
  By Theorem \ref{thm:TWN_weak_compactness} the sequence $(u_n)_n$ converges over a subsequence towards $u_\infty$,
  where $u_\infty$ is another solution of \eqref{eq:TWN}.
  We need to verify that $u_\infty \in \Inv(E)$.
  Since $u_\infty \in \cl{E}^{C^1_{\text{loc}}(\R,\R^d)}$, if follows by property (3) in the definition of an isolating trajectory neighbourhood that $u_\infty \in \cl{E}^{W^{1,\infty}(\R,\R^d)}$.
Since by definition $\cl E \setmin E$ does not contain solutions of \eqref{eq:TWN}, it follows that $u_\infty \in E$.
Consequently, $u_\infty \in \Inv(E)$.
\end{myproof}

\begin{mydef}
A sequence of functions $(u_\omega)_{\omega \in \Omega}$ is called an $s$-parameter family if $\Omega \subseteq \R^s$ is connected, $\cl{ \inter{ \Omega } } = \cl{ \Omega }$,
and $\omega \mapsto u_\omega$ is continuous from $\Omega$ into $W^{1,\infty}(\R,\R^d)$.
\end{mydef}

\begin{mylemma}
  Suppose $E$ is an isolating trajectory neighbourhood, and $(u_\omega)_{\omega \in \Omega}$, an $s$-parameter family of solutions of \eqref{eq:TWN}.
  Then if $u_{\omega_0} \in E$ for some $\omega_0 \in \Omega$, it follows that $u_\omega \in E$ for all $\omega \in \Omega$.
  In particular, an isolated invariant set is translation invariant.
\end{mylemma}
\begin{myproof}
  Denote by $p : \Omega \to W^{1,\infty}(\R,\R^d)$ the map $p(\omega) := u_\omega$.
  Since $u_\omega$ is a solution of \eqref{eq:TWN} and $E$ is an isolating trajectory neighbourhood, we have $u_\omega \in E \cup \cl{E}^C$.
  Hence $p(\Omega) \subseteq E \cup \cl{E}^C$.
  But as $\Omega$ is connected and $p$ is continuous, the image $p(\Omega)$ is connected.
  Furthermore, ${p(\Omega) \cap E \neq \emptyset}$ by assumption.
  Hence $p(\Omega) \subseteq E$, completing the proof.
\end{myproof}


\section{Heteroclinic solutions}

\label{sec:TWN_heteroclinics}

In this section we study conditions on $\Psi$ that ensure that solutions \eqref{eq:TWN} which are bounded in $L^\infty$ are heteroclinics.

\subsection{Tame solutions and perturbations}

We will now introduce the notion of tame solutions and perturbations.
These should be interpreted as a notion of smallness of solutions with respect the perturbation $\Psi$,
and by duality a notion of smallness of the perturbation $\Psi$ with respect to solutions,
under which the gradient-like behaviour of \eqref{eq:TWN} persists.

\begin{mydef}[$\Psi$-tame solutions]
A subset $E \subset C^1_{\text{loc}}(\R,\R^d)$ with
\[
R = \sup_{ u \in E } \| u \|_{L^\infty(\R,\R^d)} < \infty,
\]
is called $\Psi$-tame provided
\begin{enumerate}
\item $\Psi$ is $R$-admissible, and
\item there exist $0 \leq \beta < 1$ and $\delta > 0$ such that, for any $u \in E$ which satisfies
\[
\| u' + \Phi(u) + \Psi(u) \|_{L^\infty(\R,\R^d)} \leq \delta
\]
it holds that
\[
\int_\R | g_u( \Psi(u) , u' ) | \d x \leq \beta \calE_{\text{kin}}(u).
\]
\end{enumerate}
A solution $u$ of \eqref{eq:TWN} is called $\Psi$-tame if the singleton set $E = \{u\}$ is $\Psi$-tame.
\end{mydef}

Dual to the notion of solutions being $\Psi$-tame, we have the notion of a perturbation being tame with respect to bounded solutions. 

\begin{mydef}[$E$-tame perturbations]
  Suppose $E \subset C^1_{\text{loc}}(\R,\R^d)$ is such that 
  \[
  \sup_{ u \in E } \| u \|_{L^\infty(\R,\R^d)} < \infty.
  \]
  A perturbation $\Psi \in \bPsi_h$ is called $E$-tame provided that the set $E$ is $\Psi$-tame.
\end{mydef}

In light of Theorem \ref{thm:kinematic_estimate_Psi} these conditions are not unreasonable.
In Section \ref{subsec:E-tame} we will show that for isolating neigbourhoods $E$, 
the $E$-tame perturbations form an open subset of the perturbation space $\bPsi_h$.
In particular, as clearly any bounded sequence of solutions corresponding to the unperturbed equation \eqref{eq:TWN} with $\Psi = 0$ is tame,
this demonstrates that $E$-tame perturbations are in abundance.

\begin{myremark}
  To retain gradient-like behaviour of \eqref{eq:TWN}, it suffices to have the estimate
  \[
\int_\R | g_u( \Psi(u) , u' ) | \d x \leq \calE_{\text{kin}}(u)
\]
for any $u \in E$ which is a solution of \eqref{eq:TWN}.
The reason we require $0 \leq \beta < 1$ and $\delta$-approximate solutions are taken into consideration is to ensure that $E$-tame perturbations are open in $\bPsi_h$,
see Section \ref{subsec:E-tame}.
\end{myremark}

\subsection{Quasi-Lyapunov function}
\label{subsec:TWN_quasi-Lyapunov}
The fundamental structure which will be studied in this section is the quasi-Lyapunov function for \eqref{eq:TWN}.

\begin{mydef}[Quasi-Lyapunov]
\label{def:QL}
  Let $\calL :W^{1,\infty}(\R,\R^d) \to \R$ be given by
  \[
  \calL(u) := \rst{\bigg( \frac 1 2 S(u) \cdot \calN[ S(u) ] + F(u) \bigg)}{0} - \calB(u).
  \]
  Here
  \begin{align*}
        \calB(u) &= \frac 1 2 \int_\R \int_0^y S(u(x-y)) \cdot N(y) DS(u(x)) u'(x) \d x \d y \\
    &\quad + \frac 1 2 \sum_{j \in \Z} \int_0^{\lambda_j} S(u(x-\lambda_j)) \cdot N_j DS(u(x)) u'(x) \d x.
  \end{align*}
  Furthermore, define $\calL_-(u)$, $\calL_+(u)$ by
  \[
  \calL_-(u) := \liminf_{\tau \to -\infty} \calL(\tau \cdot u), \qquad \calL_+(u) := \limsup_{\tau \to \infty} \calL(\tau \cdot u).
  \]
\end{mydef}

In light of Theorem \ref{thm:TWN_weak_compactness} any solution $u$ of \eqref{eq:TWN} with $\| u \|_{L^\infty(\R,\R^d)} < \infty$ is in fact an element of $W^{1,\infty}(\R,\R^d)$.
Hence $\calL$ is defined on bounded solutions of \eqref{eq:TWN}.
We now present an integral identity for \eqref{eq:TWN}, which is a variation of the nonlocal Noether theorem developed in \cite{bakker2018noether}.

\begin{mylemma}
\label{lemma:QL_integral_identity}
  Suppose $u \in C^1_{\text{loc}}(\R,\R^d)$ solves \eqref{eq:TWN} and satisfies $\| u \|_{L^\infty(\R,\R^d)} < \infty$.
  Then the following identity holds:
  \begin{equation*}
    \int_a^b g_u\big( u' + \Psi(u) , u' \big) \d x = \calL(a \cdot u) - \calL(b \cdot u), \qquad -\infty < a \leq b < \infty.
  \end{equation*}
\end{mylemma}
\begin{myproof}
Let
\[
\calA_a^b(u) := \int_a^b \bigg( \frac 1 2 S(u) \cdot \calN[S(u)] + F(u) \bigg) \d x,
\]
where $\calN$ is as in Hypothesis \ref{hyp:calN}, i.e.,
\[
\calN[v](x) = \int_\R N(y) v(x-y) \d y + \sum_{j\in\Z} N_j v(x-\lambda_j).
\]
Since $u$ is a $L^\infty$ bounded solution of \eqref{eq:TWN}, Theorem \ref{thm:TWN_weak_compactness} implies that $\| u \|_{W^{1,\infty}(\R,\R^d)} < \infty$.
We are thereby allowed to differentiate the map $\tau \mapsto \calA_a^b(\tau \cdot u)$ under the integral sign.
Differentiating the identity $\calA_a^b(u) = \calA_{a-\tau}^{b-\tau}(\tau\cdot u)$ with respect to $\tau$ yields
\begin{equation}
  \label{eq:D_action}
  \begin{split}
      0 
      &= - \eval{\bigg( \frac 1 2 S(u) \cdot \calN[S(u)] + F(u) \bigg)}{a}{b} \\
  &\quad + \frac 1 2 \int_a^b DS(u) u' \cdot \calN[S(u)] \d x + \int_a^b DF(u) u' \d x \\
  &\quad + \frac 1 2 \int_a^b S(u) \cdot \calN[ DS(u) u' ] \d x \\
  &= - \eval{\bigg( \frac 1 2 S(u) \cdot \calN[S(u)] + F(u) \bigg)}{a}{b} \\
  &\quad + \int_a^b DS(u) u' \cdot \calN[S(u)] \d x + \int_a^b DF(u) u' \d x \\
  &\quad + \frac 1 2 \int_a^b S(u) \cdot \calN[ DS(u) u' ] \d x - \frac 1 2 \int_a^b D S(u) u' \calN[S(u)] \d x.
  \end{split}
\end{equation}
By definition of the Riemannian metric $g$ and using that $u$ solves \eqref{eq:TWN} we have
\begin{align*}
    \int_a^b DS(u) u' \cdot \calN[S(u)] \d x + \int_a^b DF(u) u' \d x 
  &= \int_a^b g_u\bigg( \nabla_g S(u)^T \calN[S(u)] + \nabla_g F(u) , u' \bigg) \d x \\
  &= - \int_a^b g_u\bigg( u' + \Psi(u), u' \bigg) \d x.
\end{align*}
Substituting this back into \eqref{eq:D_action} and rearranging terms gives
\[
    \int_a^b g_u\big( u' + \Psi(u) , u' \big) \d x = - \eval{\bigg( \frac 1 2 S(u) \cdot \calN[S(u)] + F(u) \bigg)}{a}{b} + \eval{\calI}{a}{b}(u),
\]
where
\begin{align*}
  \eval{\calI}{a}{b}(u) &:= \frac 1 2 \int_a^b \int_\R S(u(x)) \cdot N(y) DS(u(x-y)) u'(x-y) \d y \d x \\
  &\quad - \frac 1 2 \int_a^b \int_\R D S(u(x)) u'(x) \cdot N(y) S(u(x-y)) \d y \d x, \\
  &\quad + \frac 1 2 \int_a^b \sum_{j\in\Z} S(u(x)) \cdot N_j DS(u(x-\lambda_j)) u'(x-\lambda_j) \d x \\
  &\quad - \frac 1 2 \int_a^b \sum_{j\in\Z} DS(u(x)) u'(x) \cdot N_j S(u(x-\lambda_j)) \d x.
\end{align*}
We thus need to verify that
\[
\eval{\calI}{a}{b}(u) = \calB(b\cdot u) - \calB(a\cdot u).
\]
This can be interpreted as a nonlocal analogue of Green's formula, as in \cite{bakker2018noether, gunzburger2010nonlocal}.

The term $\eval{\calI}{a}{b}(u)$ fits within the more general class of functions of the form
\begin{align*}
  \eval{\calI_\mu}{a}{b}(u) &:= \frac 1 2 \int_a^b \int_\R S(u(x)) \cdot K(y) DS(u(x-y)) u'(x-y) \d\mu( y ) \d x \\
  &\quad - \frac 1 2 \int_a^b \int_\R D S(u(x)) u'(x) \cdot K(y) S(u(x-y)) \d\mu( y ) \d x,
\end{align*}
where $\mu$ is a Radon measure on $\R$.
In light of the symmetry and decay imposed by Hypothesis \ref{hyp:calN}, we may assume $\mu(B) = \mu(-B)$ for any Borel set $B \subseteq \R$,
and the map $K : \R \to \Mat_{D\times D}(\R)$ satisfies
\[
K(y) = K(-y) = K(y)^T, \qquad \int_\R ( 1 + |y| ) \| K(y) \| \d\mu( y ) < \infty.
\]
  Note that by Fubini's theorem and the symmetry imposed on $K$ and $\mu$, we have
  \begin{align*}
    \eval{ \calI_\mu }{a}{b}(u) &= \frac 1 2 \int_\R \int_a^b S(u(x)) \cdot K(y) DS(u(x-y)) u'(x-y) \d x \d\mu(y) \\
    &\quad - \frac 1 2 \int_\R \int_a^b D S(u(x)) u'(x) \cdot K(y) S(u(x-y)) \d x \d \mu(y), \\
    &= \frac 1 2 \int_\R \int_a^b S(u(x)) \cdot K(y) DS(u(x+y)) u'(x+y) \d x \d\mu(y) \\
    &\quad - \frac 1 2 \int_\R \int_a^b  S(u(x-y)) \cdot K(y) D S(u(x)) u'(x) \d x \d \mu(y).
   \end{align*}
   A change of variables $x \leadsto x - y$ in the first term in the right hand side yields
   \[
    \eval{ \calI_\mu }{a}{b}(u) = \frac 1 2 \int_\R \left( \int_{a+y}^{b+y} - \int_{a}^{b} \right) S(u(x-y)) \cdot K(y) DS(u(x)) u'(x) \d x \d \mu(y).
  \]
  Note that
  \begin{equation}
    \label{eq:calId_moment}
    \begin{split}
          \left( \int_{a+y}^{b+y} - \int_{a}^{b} \right) &S(u(x-y)) \cdot K(y) DS(u(x)) u(x) \d x \\
      &= \int_b^{b+y} S(u(x-y)) \cdot K(y) DS(u(x)) u(x) \d x \\
      &\quad - \int_a^{a+y} S(u(x-y)) \cdot K(y) DS(u(x)) u(x) \d x.
    \end{split}
  \end{equation}
  Since $\| u \|_{W^{1,\infty}(\R,\R^d)} < \infty$, there exists a constant $C > 0$ so that for any $\tau \in \R$
  \[
    \int_\R \left| \int_\tau^{\tau+y} S(u(x-y)) \cdot K(y) DS(u(x)) u(x) \d x \right| \d \mu(y) \leq C \int_\R | y | \| K(y) \| \d \mu(y) < \infty.
  \]
  Hence the terms in the right hand side of \eqref{eq:calId_moment} are integrable.
  We find that
  \begin{align*}
    \eval{ \calI_\mu }{a}{b}(u) &= \frac 1 2 \int_\R \int_b^{b+y} S(u(x-y)) \cdot K(y) DS(u(x)) u(x) \d x \d \mu(y) \\
    &\quad - \frac 1 2 \int_\R \int_a^{a+y} S(u(x-y)) \cdot K(y) DS(u(x)) u(x) \d x \d \mu(y) \\
    &= \calB_\mu(b\cdot u) - \calB_\mu(a \cdot u),
  \end{align*}
  where
  \[
  \calB_\mu(u) = \frac 1 2 \int_\R \int_0^y S(u(x-y)) \cdot K(y) DS(u(x)) u'(x) \d x \d \mu(y).
  \]
  We then retrieve $\eval{\calI}{a}{b}(u) = \calB(b\cdot u) - \calB(a\cdot u)$ by letting 
  \[
  \mu = \lambda^1 + \sum_{j \in \Z} \delta_{\lambda_j},
  \]
  where $\lambda^1$ be the Lebesgue measure on $\R$, and
  \[
  K(y) = 
  \left\{
    \begin{array}{l l}
      N(y) + N_j & \text{if } y = \lambda_j, \\
      N(y) & \text{otherwise}.
    \end{array}
    \right.
  \]
\end{myproof}

\begin{mylemma}
\label{lemma:QL_inequality}
  Suppose $u$ is a $\Psi$-tame solution of \eqref{eq:TWN} with $\calE_{\text{kin}}(u) < \infty$.
  Then
  \[
  \calE_{\text{kin}}(u) \leq \calL_-(u) - \calL_+(u).
  \]
\end{mylemma}
\begin{myproof}
  By Lemma \ref{lemma:QL_integral_identity} we have 
  \[
    \int_a^b g_u( u' + \Psi(u) , u' ) \d x = \calL(a \cdot u) - \calL(b \cdot u), \qquad -\infty < a \leq b < \infty.
  \]
  Note that
  \[
  \int_a^b g_u( u' , u' ) \d x - \int_a^b | g_u( \Psi(u) , u' ) | \d x \leq \int_a^b g_u( u' + \Psi(u) , u' ) \d x.
  \]
  Combining these identities, we obtain the inequality
  \begin{equation}
    \label{eq:QL_intermediate_inequality}
    \int_a^b g_u( u' , u' ) \d x - \int_a^b | g_u( \Psi(u) , u' ) | \d x \leq \calL(a \cdot u) - \calL(b \cdot u).
  \end{equation}
  The integrals on the left hand side are estimated by
  \[
    \int_a^b g_u( u' , u' ) \d x \leq \int_\R g_u(u',u') \d x = 2 \calE_{\text{kin}}(u) < \infty, 
  \]
  and
  \[
    \int_a^b | g_u( \Psi(u) , u' ) | \d x \leq \int_\R | g_u( \Psi(u) , u' ) | \d x \leq \beta \calE_{\text{kin}}(u) < \infty,
  \]
  where we used the fact that $u$ is $\Psi$-tame and $\calE_{\text{kin}}(u) < \infty$.
  These estimates ensure we can apply the dominated convergence theorem in \eqref{eq:QL_intermediate_inequality} to take the limit inferior as $a \to -\infty$ and $b \to \infty$, to find
  \[
  \int_\R g_u(u',u') \d x - \int_\R | g_u(\Psi(u),u') | \d x \leq \calL_-(u) - \calL_+(u).
  \]
  Finally, since $u$ is $\Psi$-tame we have
  \[
  \calE_{\text{kin}}(u) \leq \int_\R g_u(u',u') \d x - \int_\R | g_u(\Psi(u),u') | \d x.
  \]
  Combining these estimates yields the desired bound.
\end{myproof}


We finish this section with a continuity result for the quasi-Lyapunov function.
\begin{mylemma}
  \label{lemma:QL_continuity}
  Let $(u_n)_n \subset C^1_{\text{loc}}(\R,\R^d)$ be a convergent sequence, with limit point $u_\infty$,
  and suppose that
  \[
  \sup_n \|u_n\|_{W^{1,\infty}(\R,\R^d)} < \infty.
  \]
  Then
  \[
  \lim_{n\to\infty} \calL(u_n) = \calL(u_\infty).
  \]
\end{mylemma}
\begin{myproof}
Convergence of the sequence $\rst{\bigg( \frac 1 2 S(u_n) \cdot \calN[S(u_n)] + F(u_n) \bigg)}{0}$ follows form Lemma \ref{lemma:convolution_weak_continuity}.
To establish the continuity of $\calL$, we thus need to prove convergence of the boundary terms 
\[
  \calB_c(u_n) = \frac 1 2 \int_\R \int_0^y S(u_n(x-y)) \cdot N(y) DS(u_n(x)) u_n'(x) \d x \d y,
\]
  and
\[
  \calB_d(u_n) = \frac 1 2 \sum_{j \in \Z} \int_0^{\lambda_j} S(u_n(x-\lambda_j)) \cdot N_j DS(u_n(x)) u_n'(x) \d x.
\]
Using the uniform bound on $\| u_n \|_{W^{1,\infty}(\R,\R^d)}$, we find there exists $C > 0$ such that the inner integral in $\calB_c(u_n)$ is bounded by
\[
\left| \int_0^y S(u_n(x-y)) \cdot N(y) DS(u_n(x)) u_n'(x) \d x \right| \leq C |y| \|N(y)\|,
\]
and similarly the inner integral in $\calB_d(u_n)$ is bounded by
\[
 \left| \int_0^{\lambda_j} S(u_n(x-\lambda_j)) \cdot N_j DS(u_n(x)) u_n'(x)  \right| \leq C |\lambda_j| \| N_j \|.
\]
By Hypothesis \ref{hyp:calN} these upper bounds $C |y| \| N(y) \|$ and $C |\lambda_j| \|N_j\|$ are integrable and summable, respectively.
This ensures we can apply the dominated convergence theorem to take the limit $n\to\infty$ under the integral signs and sums appearing in $\calB_c(u_n)$ and $\calB_d(u_n)$.
Hence $\calB_c(u_n) \to \calB_c(u_\infty)$ and $\calB_d(u_n) \to \calB_d(u_\infty)$, and consequently $\calL(u_n) \to \calL(u_\infty)$ as $n\to\infty$, as required.
\end{myproof}

\subsection{Gradient-like dichotomy}

We are now prepared to prove the gradient-like dichotomy for \eqref{eq:TWN}.
\begin{mythm}
\label{thm:TWN_gradlike}
  Suppose $u$ is a $\Psi$-tame solution of \eqref{eq:TWN} with $\calE_{\text{kin}}(u) < \infty$.
Then $u$ is either constant, or there exist constant solutions $z_-$, $z_+$ of \eqref{eq:TWN}, with $z_- \neq z_+$, such that
$\big( u(x) , u'(x) \big) \to \big(z_\pm,0\big)$ as $x\to\pm\infty$.
  Furthermore, the kinetic energy satisfies the explicit bound
  \[
        \calE_{\text{kin}}(u) \leq \calL(z_-) - \calL(z_+) = h(z_-) - h(z_+).
  \]
\end{mythm}
\begin{myproof}
  Assume $u$ is nonconstant.
  Since $\| u \|_{L^\infty(\R,\R^d)} < \infty$ we have $| \calL_\pm(u) | < \infty$, so that by Lemma \ref{lemma:QL_inequality} we find that $\calE_{\text{kin}}(u) < \infty$.
  By Lemma \ref{lemma:TWN_limit_sets} we find that $\big( u(x) , u'(x) \big) \to \big( z_\pm , 0 \big)$ as $x \to \pm\infty$,
  where $z_-$, $z_+$ are constant solutions of \eqref{eq:TWN}.
  Left to prove is that $z_- \neq z_+$.
  
  By Theorem \ref{thm:TWN_weak_compactness} we have $\| u \|_{W^{1,\infty}(\R,\R^d)} < \infty$.
  Thus Lemma \ref{lemma:QL_continuity} is applicable to the sequences $( n \cdot u)_{n\in\N}$ and $( (-n) \cdot u)_{n\in\N}$ of translated solutions.
  By the definition of $\calL_\pm$ and Lemma \ref{lemma:QL_continuity} we then have
  \[
  \calL_-(u) = \lim_{n \to \infty} \calL((-n) \cdot u) = \calL(z_-), \qquad \calL_+(u) = \lim_{n\to\infty} \calL(n \cdot u) = \calL(z_+).
  \]
  Combined with Lemma \ref{lemma:QL_inequality}, using that $u$ is nonconstant so that $\calE_{\text{kin}}(u) > 0$, we then find that
  \[
  \calL(z_+) = \calL_+(u) < \calL_-(u) = \calL(z_-).
  \]
  It follows that $z_- \neq z_+$.
  The estimate on $\calE_{\text{kin}}(u)$ now follows from Lemma \ref{lemma:QL_inequality}, combined with the observation that $\calL(z) = h(z)$ for any constant $z\in \R^d$.
\end{myproof}

\begin{myremark}
\label{remark:TWN_generalised_heteroclinics}
  We have done our analysis under the assumption that all constant solutions of \eqref{eq:TWN} are hyperbolic.
  In fact, the construction of the perturbation $\Psi$ only makes sense under this assumption.
  However, in the absence of the perturbation $\Psi$, i.e., assuming $\Psi = 0$, equation \eqref{eq:TWN} makes sense even without assuming hyperbolicity.
  The gradient-like dichotomy extends to this situation to obtain the following result.
  For any bounded solution $u$ the limit sets $\alpha(u)$, $\omega(u)$ are connected and consist of constant solutions of \eqref{eq:TWN}.
  Furthermore, if $u$ is nonconstant, then $\alpha(u) \cap \omega(u) = \emptyset$.
\end{myremark}

\subsection{Moduli spaces}

Let $z_-$, $z_+$ be critical points of $h$.
Denote by
\[
\calM_{\Phi,\Psi}(z_-,z_+) = \set{ u \in C^1_{\text{loc}}(\R,\R^d) }{
  \begin{array}{c}
    u' + \Phi(u) + \Psi(u) = 0 \\
    \big(u(x),u'(x)\big) \to \big( z_\pm , 0 \big) \text{ as } x\to\pm\infty 
  \end{array}
  }
\]
the moduli space of solutions of \eqref{eq:TWN} connecting $z_-$ with $z_+$.
Given an isolating trajectory neighbourhood $E$, we define
\[
\calM_{\Phi,\Psi}(z_-,z_+ ; E) := E \cap \calM_{\Phi,\Psi}(z_-,z_+).
\]
We assume from here on out that $\Psi$ is $R$-admissible, where $R = \sup_{u\in E} \| u \|_{L^\infty(\R,\R^d)}$.
Whenever the choice of $\Phi$ and/or $\Psi$ is clear from the context, we shall suppress them in the notation.

In light of Theorem \ref{thm:TWN_nonlinear_decay_rates} the moduli spaces $\calM(z_-,,z_+)$ is a subspace of the path space $\calP(z_-,z_+)$.
As such, there are now various relevant modes of convergence on $\calM(z_-,z_+)$, namely, the convergence in $C^1_{\text{loc}}(\R,\R^d)$,
uniform convergence, and the convergence in the affine Hilbert space topology of $\calP(z_-,z_+)$.
In this section we further investigate the relations between these modes of convergens, 
and use this knowledge to classify the boundaries of noncompact components of moduli spaces.

\subsubsection{Chains of broken solutions}
We already inspected the breaking of sequences in Section \ref{subsec:sojourn}.
Now, with the gradient-like structure of \eqref{eq:TWN} established, we consider the breaking of solutions once again, this time to characterise noncompactness of moduli spaces.
\begin{mythm}
  \label{thm:TWN_broken_chain}
  Suppose $\Psi$ is $R$-admissible, with $R = \sup_{u\in E} \|u\|_{L^\infty(\R,\R^d)}$.
  Given a sequence $(u_n)_n \subset \calM(z_-,z_+;E)$, there exists $k \in \N \cup \{0\}$ and constant solutions $z_0, z_1, \dots, z_k, z_{k+1}$ of \eqref{eq:TWN},
  with $z_0 = z_-$ and $z_{k+1} = z_+$, such that the following holds:
  \begin{enumerate}
  \item The constant solutions are distinct and satisfy
    \[
    \calL(z_-) > \calL(z_1) > \cdots > \calL(z_k) > \calL(z_+).
    \]
  \item For each $j \in \{1,\dots,k\}$ there exist a sequence $(\tau^0_{j,n})_n \subset \R$,
    such that $(u_n)_n$ breaks along a subsequence of $(\tau^0_{j,n})_n$ towards $z_j$, in the sense of Definition \ref{def:TWN_breaking}.
  \item For each $j \in \{0,\dots,k\}$ there exists a sequence $(\tau^+_{j,n})_n \subset \R$,
    such that $(\tau^+_{j,n}\cdot u_n)_n$ converges over a subsequence in $C^1_{\text{loc}}(\R,\R^d)$ towards a solution $v_j$ of \eqref{eq:TWN}, with
    \[
    \lim_{x\to -\infty} \big( v_j(x), v_j'(x) \big) = \big( z_j, 0 \big), \qquad \lim_{x\to \infty} \big( v_j(x) , v_j'(x) \big) = \big( z_{j+1} , 0 \big).
    \]
  \end{enumerate}
  We call such a $(k+1)$-tuple $(v_0,\dots,v_k)$ a $k$-fold broken chain of solutions.
\end{mythm}
\begin{myproof}
  In light of the exponential decay Theorem \ref{thm:TWN_nonlinear_decay_rates}, one has $\calE_{\text{kin}}(u_n) < \infty$ for each $n$.
  Consequently, the gradient-like dichotomy Theorem \ref{thm:TWN_gradlike} gives the explicit and uniform upper bound
  \[
  \sup_n \calE_{\text{kin}}(u_n) \leq \calL(z_-) - \calL(z_+).
  \]
  Now we can argue as in the proof of Lemma \ref{lemma:Bn_diameter_blowup_breaking}.
  This gives an iterative procedure of constructing a broken chain.
  Theorem \ref{thm:TWN_gradlike} ensures that property (1) is satisfied along any broken chain.
  Since $E$ is an isolating trajectory neighbourhood and $h$ is Morse, 
  there are only finitely many combinations of constant solutions $z_i \in E$ of \eqref{eq:TWN} which satisfy property (1).
  Therefore the iterative procedure of constructing a broken chain must terminate after finitely many steps.
\end{myproof}

\subsubsection{Modes of convergence}
We now investigate the relation between the various modes of convergence on $\calM(z_-,z_+;E)$.

\begin{mythm}
\label{thm:TWN_moduli_convergence}
  Suppose $\Psi$ is $R$-admissible, with $R = \sup_{u\in E} \|u\|_{L^\infty(\R,\R^d)}$.
  Given $(u_n)_n \subset \calM(z_-,z_+;E)$ the following statements are equivalent:
  \begin{enumerate}
  \item $(u_n)_n$ converges to $u_\infty \in \calM(z_-,z_+;E)$ in the topology of $\calP(z_-,z_+)$,
  \item $(u_n)_n$ converges to $u_\infty \in \calM(z_-,z_+;E)$ in the topology of $W^{1,\infty}(\R,\R^d)$,
  \item $(u_n)_n$ converges to $u_\infty \in \calM(z_-,z_+;E)$ in the topology of $C^1_{\text{loc}}(\R,\R^d)$.
  \end{enumerate}
\end{mythm}
\begin{myproof}
We will show (1)$\Longleftrightarrow$(2) and (2)$\Longleftrightarrow$(3).

\paragraph{(1)$\implies$(2).}
Differentiating \eqref{eq:TWN} and rearranging terms, we find
\[
- ( u''_\infty - u''_n ) = D \Phi(u_\infty) u'_\infty - D \Phi(u_n) u_n + D \Psi(u_\infty) u'_\infty - D \Psi(u_n) u'_n.
\]
By Lemma \ref{lemma:Phi_L2_smooth} and Lemma \ref{lemma:Psi_L2_smooth}, and using the assumption that $\| u_\infty - u_n \|_{W^{1,2}(\R,\R^d)} \to 0$ as $n\to\infty$,
we find that $\| u''_\infty - u''_n \|_{L^2(\R,\R^d)} \to 0$ as $n\to\infty$, hence $\| u_\infty - u_n \|_{W^{2,2}(\R,\R^d)} \to 0$ as $n\to\infty$.
  In light of Morrey's inequality this implies $\| u_\infty - u_n \|_{W^{1,\infty}(\R,\R^d)} \to 0$ as $n\to\infty$.

\paragraph{(2)$\implies$(1).}
  By Theorem \ref{thm:TWN_nonlinear_decay_rates}, given $\eta \in \calR_{u_\infty}$ with $\eta > 0$ there exists $C > 0$ and $n_0 \in \N$ such that for $n \geq n_0$
  \begin{align*}
    | u_n(x) - z_- | + | u_n'(x) | \leq C e^{-\eta |x|}, & \qquad x < 0, \\
    |u_n(x) - z_+ | + |u_n'(x) | \leq C e^{-\eta|x|}, & \qquad x \geq 0.
  \end{align*}
  Thus, given $\epsilon > 0$ we may choose $L > 0$ so that
  \[
  \| u_\infty - u_n \|^2_{W^{1,2}(\R \setmin [-L,L] , \R^d)} < \epsilon^2 / 2, \qquad n \geq n_0.
  \]
  Moreover, by uniform convergence combined with the dominated convergence theorem, there exists $n_1 \in \N$ so that
  \[
  \| u_\infty - u_n \|^2_{W^{1,2}([-L,L] , \R^d)} < \epsilon^2 / 2, \qquad n \geq n_1.
  \]
  Hence, for $n \geq \max\{ n_0 , n_1 \}$ we find
  \[
  \| u_\infty - u_n \|_{W^{1,2}(\R,\R^d)} = \sqrt{ \| u_\infty - u_n \|^2_{W^{1,2}(\R \setmin [-L,L] , \R^d)}  +  \| u_\infty - u_n \|^2_{W^{1,2}([-L,L] , \R^d)} } < \epsilon.
  \]

\paragraph{(2)$\implies$(3).}
This implication is trivial.

\paragraph{(3)$\implies$(2).}
  Pick any $0 < \epsilon < | z_+ - z_- |$.
  In light of the exponential decay Theorem \ref{thm:TWN_nonlinear_decay_rates}, one has $\calE_{\text{kin}}(u_n) < \infty$ for each $n$.
  Consequently, the gradient-like dichotomy Theorem \ref{thm:TWN_gradlike} gives the explicit and uniform upper bound
  \[
  \sup_n \calE_{\text{kin}}(u_n) \leq \calL(z_-) - \calL(z_+).
  \]
  Hence that Lemma \ref{lemma:Bn_diameter_blowup_breaking} is applicable with $\rho = \epsilon/2$ and $B_n = A_\rho(u_n)$.
  Since $(u_n)_n$ converges to another element of $\calM(z_-,z_+;E)$, Theorem \ref{thm:TWN_broken_chain} ensures $(u_n)_n$ cannot break along any sequence $(\tau^0_n)_n$.
  It follows from Lemma \ref{lemma:Bn_diameter_blowup_breaking} that there exists $L > 0$ so that
  \[
  \sup_n \diam \set{ x \in \R }{ \inf_{z\in\crit(h)} |u_n(x)-z| \geq \frac \epsilon 2 } \leq L.
  \]
  Hence there exists a sequence $(\tau_n)_n \subset \R$ such that
  \[
  \set{ x \in \R }{ \inf_{z\in\crit(h)} |u_n(x)-z| \geq \frac \epsilon 2 } \subseteq [ \tau_n - L/2 , \tau_n + L/2 ].
  \]

  We claim that $(\tau_n)_n$ is bounded.
  If not, we could find a subsequence $(\tau_{n_k})_k$ which grows unbounded, say, $\tau_{n_k} \to \infty$ as $k \to \infty$.
  But then, for any $x \in \R$,
  \[
  | u_\infty(x) - z_- | \leq | u_\infty(x) - u_n(x) | + | u_n(x) - z_- | < \epsilon
  \]
  for $n$ sufficiently large.
  Hence $\| u_\infty - z_- \|_{L^\infty(\R,\R^d)} < \epsilon < | z_+ - z_- |$.
  But by assumption $u_\infty \in \calM(z_-,z_+)$, so that $u_\infty(x) \to z_+$ as $x\to\infty$, hence $\| u_\infty - z_- \|_{L^\infty(\R,\R^d)} \geq | z_+ - z_- |$.
  A similar contraction is obtained if $\tau_{n_k} \to -\infty$ as $k \to \infty$.
  We conclude that $(\tau_n)_n$ is indeed bounded.

  We thus find that 
  \[
  A_{\epsilon / 2} := \bigcup_n \set{ x \in \R }{ \inf_{z\in\crit(h)} |u_n(x)-z| \geq \frac \epsilon 2 }
  \]
  is a bounded subset of $\R$, say, $A_{\epsilon / 2} \subseteq [-R,R]$, for some $R > 0$.
  Then
  \[
  \| u_\infty - u_n \|_{L^\infty(\R \setmin [-R,R],\R^d)} < \epsilon, \qquad \text{for all} \quad n. 
  \]
  Moreover, by $C^1_{\text{loc}}(\R,\R^d)$ convergence, we can find $n_0 \in \N$ so that
  \[
  \| u_\infty - u_n \|_{L^\infty([-R,R],\R^d)} < \epsilon, \qquad n \geq n_0.
  \]
  Hence
  \[
  \| u_\infty - u_n \|_{L^\infty(\R,\R^d)} < \epsilon, \qquad n \geq n_0.
  \]
  Convergence in $W^{1,\infty}(\R,\R^d)$ then follows from Lemma \ref{lemma:TWN_Phi_uniform_smooth} and Lemma \ref{lemma:TWN_Psi_uniform_smooth} using a bootstrap argument.
\end{myproof}

\subsubsection{Classification of moduli spaces}

We now wish to exploit the fact that $\calM(z_-,z_+)$ is obtained as the zero level set of the map
\[
\partial_x + \Phi(\cdot) + \Psi(\cdot) : \calP(z_-,z_+) \to L^2(\R,\R^d).
\]
To do so, we first introduce the notions of transverse maps and $E$-transverse perturbations.

\begin{mydef}[Transverse maps]
  Equation \eqref{eq:TWN} is said to be transverse on $\calM(z_-,z_+)$ 
  if for any $u \in \calM_{\Phi,\Psi}(z_-,z_+)$ the linearisation 
  \[
  \partial_x + D\Phi(u) + D\Psi(u) : W^{1,2}(\R,\R^d) \to L^2(\R,\R^d)
  \]
  is surjective.
\end{mydef}

\begin{mydef}[$E$-transverse perturbations]
  Let $E$ be an isolating trajectory neighbourhood.
  A perturbation $\Psi$ is said to be $E$-transverse up to index $k$
  provided that for any pair $z_- , z_+ \in E$ of constant solutions to \eqref{eq:TWN} with $m_h(z_-) - m_h(z_+) \leq k$, 
  equation \eqref{eq:TWN} with this perturbation $\Psi$ is transverse on $\calM_{\Phi,\Psi}(z_-,z_+)$. 
\end{mydef}
We will show in Section \ref{sec:TWN_generic} that $E$-transverse perturbations up to index $2$ are generic.

Suppose \eqref{eq:TWN} is transverse on $\calM(z_-,z_+)$.
Surjectivity of the linearisation combined with the Fredholm property allows us to choose a bounded right inverse to
\[
  \partial_x + D\Phi(u) + D\Psi(u) : W^{1,2}(\R,\R^d) \to L^2(\R,\R^d), \qquad u \in \calM(z_-,z_+).
\]
We can therefore apply the implicit function theorem to find that $\calM(z_-,z_+)$ is a $C^3$ smooth finite dimensional manifold of dimension $\dim \calM(z_-,z_+) = m_h(z_-) - m_h(z_+)$.

The classification of $\calM_{\Phi,\Psi}(z_-,z_+;E)$ up to dimension $2$, for $E$-transverse perturbations $\Psi$, then follows a geometric construction which is typical in Morse theory.
Since this argument is well-documented, see e.g.\ \cite{weber2010morse, schwarz1993morse, audin2014morse, salamon1999lectures}, we do not give details here but only summarise the results.
Essential is the proper and free $\R$-action on $\calM_{\Phi,\Psi}(z_-,z_+;E)$ given by $u \mapsto \tau \cdot u$, for $\tau \in \R$.
We then obtain $C^1$ manifolds $\calM_{\Phi,\Psi}(z_-,z_+;E) / \R$ of dimension $m_h(z_-) - m_h(z_+) - 1$.
This reduction on the dimension allows for the exploitation of the classification of $1$-dimensional manifolds without boundary.

By transversality the Morse index $m_h$ acts like a discrete Lyapunov function.
Thus the statement of Theorem \ref{thm:TWN_broken_chain} may be supplemented with the estimate
\[
m_h(z_-) > m_h(z_1) > \cdots > m_h(z_k) > m_h(z_+).
\]
Hence, when $m_h(z_-) - m_h(z_+) = 1$, there cannot be any broken chains in the closure of $\calM_{\Phi,\Psi}(z_-,z_+;E)$.
Combined with Theorem \ref{thm:TWN_moduli_convergence} this then shows the $0$ dimensional quotient manifold $\calM(z_-,z_+;E) / \R$ is compact (hence finite).
If $m_h(z_-) - m_h(z_+) = 2$, then either $(u_n)_n$ converges to another element of $\calM_{\Phi,\Psi}(z_-,z_+;E)$, or it converges to a $1$-fold broken chain of solutions.
Finally, by heteroclinic glueing, 
any $1$-fold broken chain of solution in fact occurs in the boundary of a unique connected component of a $2$ dimensional moduli space $\calM_{\Phi,\Psi}(z_-,z_+;E)$.

This leads to the following classification.
\begin{mythm}
\label{thm:TWN_moduli_closure}
  Let $E$ be an isolating trajectory neighbourhood for \eqref{eq:TWN}.
  Assume that $\Psi$ is $E$-tame and $E$-transverse.
  Let $O \subseteq \calM_{\Phi,\Psi}(z_-,z_+;E)$ be a connected component.
  Then
  \begin{enumerate}
  \item if $m_h(z_-) - m_h(z_+) = 1$, then the set $O / \R$ is finite.
  \item if $m_h(z_-) - m_h(z_+) = 2$, then either $O / \R$ is compact and homeomorphic to $\sphere^1$,
    or $O / \R$ is noncompact and homeomorphic to $(0,1)$.
    The noncompact endpoints $0$ and $1$ correspond to distinct $1$-fold broken chain of solutions.
    Conversely, any $1$-fold broken chain of solutions is a boundary component of a moduli space of dimension $2$.
    In particular, $1$-fold broken chains of solutions always occur in pairs inside $E$.
  \end{enumerate}
\end{mythm}


\section{Generic properties of $\Psi$}
\label{sec:TWN_generic}

In order to apply the results from the previous section, we need to choose $\Psi$ which are both $E$-tame and $E$-transverse.
The aim is to achieve this via small perturbations $\Psi$ of a given perturbation $\Psi_0$.
Then we need to ensure that $E$ remains an isolating trajectory neighbourhood under these perturbations.
This last property we call $E$-stability of $\Psi$.
The strategy is now as follows.
We first show that $E$-stable and $E$-tame perturbations both form open subsets of the perturbation space $\bPsi_h$.
Then, if our initial perturbation $\Psi_0$ is both $E$-stable and $E$-tame (e.g. $\Psi = 0$, for which $E$-admissibility is trivial), 
we obtain an open neighbourhood of $\Psi_0$ on which both these desirable properties hold.
We can then apply a Sard--Smale transversality construction on this open neighbourhood, to obtain a residual subset of perturbations
which are $E$-stable, $E$-tame, and $E$-transverse.

\subsection{$E$-stable perturbations}

\begin{mydef}[$E$-stable perturbations]
  Let $\Psi_0 \in \bPsi_h$ and suppose $E$ is an isolating trajectory neighbourhood for \eqref{eq:TWN} with $\Psi = \Psi_0$.
  Then a perturbation $\Psi_1 \in \bPsi_h$ said to be an $E$-stable perturbation of $\Psi_0$ if $E$ is an isolating trajectory neighbourhood for \eqref{eq:TWN} with $\Psi=\Psi_1$.
\end{mydef}

We will now prove that any small perturbation of $\Psi_0$ is in fact $E$-stable.

\begin{mythm}
\label{thm:E-stable}
  Let $\Psi_0 \in \bPsi_h$, and suppose $E$ is an isolating trajectory neighbourhood for \eqref{eq:TWN} with $\Psi = \Psi_0$.
  Then there exists $\epsilon > 0$ such that any $\Psi \in \bPsi_h$ which satisfies $\| \Psi - \Psi_0 \|_{\bPsi_h} < \epsilon$ is $E$-stable.
\end{mythm}
\begin{myproof}
Since $E$ is an isolating trajectory neighbourhood for \eqref{eq:TWN} with $\Psi = \Psi_0$,
there exists $\delta > 0$ such that
\[
\| u' + \Phi(u) + \Psi_0(u) \|_{L^\infty(\R,\R^d)} > \delta, \qquad \text{for any} \quad u \in \cl{E} \setmin E.
\]
Now invoke Lemma \ref{lemma:TWN_Psi_uniform_bound} to obtain a constant $C > 0$ such that
\[
\sup_{ u \in E } \| \Psi(u) \|_{L^\infty(\R,\R^d)} \leq C \| \Psi \|_{\bPsi_h}, \qquad \text{for all} \quad \Psi \in \bPsi_h.
\]
Then set $\epsilon = \delta / (2 C)$ and suppose $\Psi \in \bPsi_h$ is such that $\| \Psi - \Psi_0 \|_{\bPsi_h} < \epsilon$.
Then for any $u \in \cl{E} \setmin E$ one has
\begin{align*}
  \| u' + \Phi(u) + \Psi(u) \|_{L^\infty(\R,\R^d)} &\geq \| u' + \Phi(u) + \Psi_0(u) \|_{L^\infty(\R,\R^d)} - \| \Psi(u) - \Psi_0(u) \|_{L^\infty(\R,\R^d)} \\
  &> \delta - C \| \Psi - \Psi_0 \|_{\bPsi_h} \\
  &> \delta / 2.
\end{align*}
Thus $\cl{E} \setmin E$ does not contain any $\delta/2$-approximate solutions of \eqref{eq:TWN} with perturbation $\Psi$, whenever $\| \Psi - \Psi_0 \|_{\bPsi} < \epsilon$.
This shows that $E$ is an isolating trajectory neighbourhood for any such $\Psi$.
\end{myproof}

\subsection{$E$-tame perturbations}
\label{subsec:E-tame}

\begin{mythm}
\label{thm:E-tame}
  Let $\Psi_0 \in \bPsi_h$, and suppose $\Psi_0$ is an $E$-tame perturbation.
  Then there exists $\epsilon > 0$ such that any $\Psi \in \bPsi_h$ with $\| \Psi - \Psi_0 \|_{\bPsi_h} < \epsilon$ is an $E$-tame perturbation.
\end{mythm}
\begin{myproof}
Since $\Psi_0$ is $E$-tame, there exist $\delta_0 > 0$ and $0 \leq \beta_0 < 1$ such that for any $u \in E$ which satisfies
\begin{equation}
  \label{eq:delta0_approximate_solution}
  \| u' + \Phi(u) + \Psi_0(u) \|_{L^\infty(\R,\R^d)} \leq \delta_0,
\end{equation}
it follows that
  \[
  \int_\R \big| g_u\big( \Psi_0(u) , u' \big) \big| \d x \leq \beta_0 \calE_{\text{kin}}(u).
  \]
By Theorem \ref{thm:kinematic_estimate_Psi}, after decreasing $\delta_0$ if necessary, there exists $C_{R,\Psi_0} > 0$ such that
for any $u \in E$ which satisfies \eqref{eq:delta0_approximate_solution}, and for any $\Psi \in \bPsi_h$, we have
\[
\| \Psi(u) - \Psi_0(u) \|_{L^2(\R,\R^d)} \leq C_{R,\Psi_0} \| \Psi - \Psi_0 \|_{\bPsi_h} \sqrt{ \calE_{\text{kin}}(u) }.
\]
Consequently, by the Cauchy--Schwarz inequality, and using the definition of the Riemannian metric $g$, there exists $C > 0$ such that
\begin{equation}
  \label{eq:Psi_pert_kin_bound}
  \int_\R | g_{u}( \Psi(u) - \Psi_0(u) , u' ) | \d x \leq C \| \Psi - \Psi_0 \|_{\bPsi_h} \calE_{\text{kin}}(u).
\end{equation}

Let $\beta_0 < \beta < 1$ and choose $\epsilon < \min\{ \delta_0 / 2 , (\beta-\beta_0)/C \}$.
In light of Lemma \ref{lemma:R-admissible_open} we can choose $\epsilon$ sufficiently small so that any $\Psi \in \bPsi_h$ with $\| \Psi - \Psi_0 \|_{\bPsi_h} < \epsilon$ 
is $R$-admissible, where $R = \sup_{u\in E} \| u \|_{L^\infty(\R,\R^d)}$.
Now choose any $0 < \delta < \delta_0 / 2$ and $\Psi \in \bPsi_h$ with $\| \Psi - \Psi_0 \|_{\bPsi_h} < \epsilon$, and suppose that $u \in E$ is such that
\begin{equation}
  \label{eq:delta_approximate_solution}
  \| u' + \Phi(u) + \Psi(u) \|_{L^\infty(\R,\R^d)} \leq \delta.
\end{equation}
Then $u$ satisfies \eqref{eq:delta0_approximate_solution} and consequently, in light of \eqref{eq:Psi_pert_kin_bound}, we have
\[
\int_\R | g_{u}( \Psi(u) - \Psi_0(u) , u' ) | \d x \leq  (\beta-\beta_0) \calE_{\text{kin}}(u).
\]
Consequently, for any $u \in E$ satisfying \eqref{eq:delta_approximate_solution} we have
\begin{align*}
  \int_\R | g_u( \Psi(u) , u' ) | \d x &\leq \int_R | g_u( \Psi_0(u) , u' ) | \d x + \int_\R | g_{u}( \Psi(u) - \Psi_0(u) , u' ) | \d x \\
  &\leq \beta \calE_{\text{kin}}(u).
\end{align*}
This shows that $\Psi$ is $E$-tame.
\end{myproof}

\subsection{$E$-transverse perturbations}

Let $\Psi_0 \in \bPsi_h$, and suppose that $\Psi_0$ is $E$-tame and $E$ is an isolating trajectory neighbourhood for \eqref{eq:TWN} with $\Psi = \Psi_0$.
We now denote by $\bPsi_h(E,\Phi+\Psi_0)$ the open connected subset of $\bPsi_h$ containing $\Psi_0$, consisting of those perturbations $\Psi$
which are $E$-stable and $E$-tame.
In light of Theorem \ref{thm:E-stable} and Theorem \ref{thm:E-tame} this set is nonempty.
In this section we will construct a residual subset $\bPsi_h^\pitchfork(E,\Phi+\Psi_0)$ of $\bPsi_h(E,\Phi+\Psi_0)$ consisting of perturbations $\Psi$ which are $E$-transverse.

Given $Z^2 \subseteq ( E \cap \crit(h) )^2$, let
\[
\Theta(Z^2) := \set{ \theta = (\ell , (z_- , z_+)) \in \N \times Z^2 }{ z_- \neq z_+ },
\]
and define a linear subspace $\bPsi_h(Z^2)$ of $\bPsi_h$ by
\[
\bPsi_h(Z^2) := \set{ \Psi \in \bPsi_h }{ 
  \begin{array}{c}
    \rst{ \Psi(u) }{x} = \sum_{\theta\in \Theta(Z^2)} \rst{ \alpha_{\theta} * \psi_{\theta}(u) }{(x,x)} \\
    \alpha_\theta \in W^{1,2}_{\eta_0}(\R^2,\R^d), \qquad \theta \in \Theta(Z^2)
  \end{array}
  }.
\]
We then define a subspace $\bPsi_h(Z^2 , E,\Phi+\Psi_0)$ of $\bPsi_h(E,\Phi+\Psi_0)$ by
\[
\bPsi_h(Z^2 , E,\Phi+\Psi_0) := \bPsi_h(E,\Phi+\Psi_0) \cap \bPsi_h(Z^2).
\]
Note that $\bPsi_h(Z^2, E,\Phi+\Psi_0)$ is an open subset of $\bPsi_h(Z^2)$.


\begin{mylemma}
\label{lemma:TWN_perturbation_surjective}
    Let $z_-, z_+ \in E$ be constant solutions of \eqref{eq:TWN}, with $z_- \neq z_+$,
    and pick any path $u \in \calP(z_-,z_+)$.
    Let $Z^2 := \{ (z_-,z_+) \}$.
    Then the map
    \[
    \bPsi_h(Z^2) \to L^2(\R,\R^d), \qquad \Psi \mapsto \Psi(u)
    \]
    has dense image.
\end{mylemma}
\begin{myproof}
Recall from Section \ref{sec:TWN_setup} that, with $\theta = ( \ell , z_- , z_+ )$, where $\ell \in \N$, we have
\[
\rst{ \psi_\theta(u) }{ \xi } := \rst{ \sigma_{\ell,z_-}(u) }{ \xi_1 } \rst{ \sigma_{\ell,z_+}(u) }{ \xi_2 }.
\]
Hence it suffices to verify that the image of the map
$P : W^{1,2}_{\eta_0}(\R^2,\R^d) \times \N \to L^2(\R,\R^d)$, given by
\[
P(\alpha,\ell) := \iint_{\R^2} \alpha\big(\xi_1,\xi_2\big) \rst{ \sigma_{\ell,z_-}(u) }{x-\xi_1} \rst{ \sigma_{\ell,z_+}(u) }{ x-\xi_2 } \d \xi_1 \d \xi_2,
\]
has dense image.

The main property of the maps $\sigma_{\ell,z_-}(u)$ and $\sigma_{\ell,z_+}(u)$ which we will exploit is that we have the pointwise limit
\begin{equation}
  \label{eq:TWN_pointwise_sigma_localisation}
  \lim_{\ell \to \infty} \rst{ \sigma_{\ell,z_\pm}(u) }{x} = \1_{D_{z_\pm}(u)}(x), 
\end{equation}
where
\[
D_{z_\pm}(u) := \set{ x \in \R }{ \rho(z_\pm) < |u(x)-z_\pm| < 2 \rho(z_\pm) },
\] 
where we recall $\rho(z_-)$ and $\rho(z_+)$ are given by
\[
\rho(z_\pm) = \frac 1 3 \inf_{ \substack{ z \in\crit(h) \\ z \neq z_\pm } } | z - z_\pm |.
\]
In particular, $| z_+ - z_- | > 2 \max\{ \rho(z_+) , \rho(z_-) \}$.
Since $u : \R \to \R^d$ defines a smooth path connecting $z_-$ with $z_+$, it follows that $D_{z_-}(u)$ and $D_{z_+}(u)$ are bounded and nonempty open subsets of $\R$.

We now observe the following.
As $D_{z_-}(u)$ and $D_{z_+}(u)$ are bounded, nonempty, and open, they have positive Lebesgue measures,
and thus
\[
x_L := \essinf D_{z_-}(u), \qquad x_R := \esssup D_{z_+}(u)
\]
define numbers $-\infty < x_L < x_R < \infty$.
Let $\epsilon_0 > 0$ be such that 
\[
(x_L , x_L + \epsilon_0) \subset D_{z_-}(u), \qquad (x_R - \epsilon_0, x_R) \subset D_{z_+}(u).
\]
Then, for any $x_0 \in \R$ and $0 < \epsilon < \epsilon_0$, define 
\[
\xi_1(x_0,\epsilon) = x_0 - x_L - \epsilon, \qquad \xi_2(x_0,\epsilon) = x_0 - x_R + \epsilon.
\]
It follows that
\begin{equation}
  \label{eq:TWN_spatial_localiser_limit}
  \1_{D_{z_-}(u)}( x - \xi_1(x_0,\epsilon) )\1_{D_{z_+}(u)}( x - \xi_2(x_0,\epsilon) ) = \1_{( x_0 - \epsilon , x_0 + \epsilon)}(x).
\end{equation}
We also refer back to Figure \ref{fig:TWN_perturbation_intuition} from Section \ref{sec:TWN_setup}.

In order to verify density of the image of $P$ in $L^2(\R,\R^d)$, it suffices to see that $\img( P )^{\perp} = \{0\}$.
Thus, assuming $v \in L^2(\R,\R^d)$ is such that $\langle P(\alpha,\ell) , v \rangle_{L^2(\R,\R^d)} = 0$ for all $(\alpha,\ell) \in W^{1,2}_{\eta_0}(\R^2,\R^d) \times \N$,
we will show that $v = 0$.

First, for given $x_0 \in \R$, and $w \in \R^d$, and sufficiently small $\epsilon > 0$, we construct
a sequence $((\alpha_n , \ell_n))_n \subset  W^{1,2}_{\eta_0}(\R^2,\R^d) \times \N$ for which we have the weak convergence in $L^2(\R,\R^d)$
\begin{equation}
  \label{eq:TWN_dual_localiser_weak_limit}
  P(\alpha_n,\ell_n) \rightharpoonup w \; \1_{(x_0-\epsilon, x_0+\epsilon)}, \qquad \text{as} \quad n \to \infty.
\end{equation}
First, by letting $\ell \to \infty$ using the dominated convergence theorem, we find using \eqref{eq:TWN_pointwise_sigma_localisation} that for any $\alpha \in W^{1,2}_{\eta_0}(\R^2,\R^d)$
\begin{equation}
  \label{eq:TWN_dual_localiser_l_limit}
  P(\alpha,\ell) \rightharpoonup \iint_{\R^2} \alpha(\xi_1,\xi_2) \1_{D_{z_-}(u)}( x - \xi_1 )\1_{D_{z_+}(u)}( x - \xi_2 ) \d \xi_1 \d \xi_2, \qquad \ell \to \infty.
\end{equation}
Choose a sequence $(\beta_n)_n \subset W^{1,2}_{\eta_0}(\R^2,\R)$ of such that we have convergence of distributions
\[
\beta_n \rightharpoonup \delta_{\big(\xi_1(x_0,\epsilon),\xi_2(x_0,\epsilon)\big)}, \qquad \text{as} \quad n\to\infty,
\]
where $\delta_\xi$ denotes the Dirac delta in $\xi \in \R^2$.
Then define $\alpha_n := w \beta_n$.
Using \eqref{eq:TWN_dual_localiser_l_limit} and a diagonal argument, we obtain a sequence $(\alpha_n , \ell_n)_n \subset W^{1,2}_{\eta_0}(\R^2,\R^d) \times\N$ such that
\[
P(\alpha_n,\ell_n) \rightharpoonup w \; \1_{D_{z_-}(u)}( x - \xi_1(x_0,\epsilon) )\1_{D_{z_+}(u)}( x - \xi_2(x_0,\epsilon) ), \qquad \text{as} \quad n \to \infty.
\]
Combined with \eqref{eq:TWN_spatial_localiser_limit} we now recover \eqref{eq:TWN_dual_localiser_weak_limit}, as desired.

Now, if $v \in L^2(\R,\R^d)$ is such that $\langle P(\alpha,\ell) , v \rangle_{L^2(\R,\R^d)} = 0$ holds for all elements $(\alpha,\ell) \in W^{1,2}_{\eta_0}(\R^2,\R^d) \times \N$,
we obtain
\[
\int_{x_0-\epsilon}^{x_0+\epsilon} w \cdot v(x) \d x = \lim_{n\to\infty} \langle P(\alpha_n,\ell_n) , v \rangle_{L^2(\R,\R^d)} = 0, \qquad x_0 \in \R,\; \epsilon > 0.
\]
By the Lebesgue differentiation theorem it follows that $w \cdot v(x_0) = 0$ for $x_0 \in \R \setmin A(w)$, where $A(w)$ is of measure zero.
Letting $\{ w_1,\dots,w_d \}$ be a basis for $\R^d$, we find that $v(x_0) = 0$ for $x_0 \in \R \setmin A$, where $A = \bigcup_{i=1}^d A(w_i)$ is of measure zero.
Thus $v = 0 \in L^2(\R,\R^d)$.
\end{myproof}

We are now prepared to prove the genericity of $E$-transverse perturbations.

\begin{mythm}
\label{thm:TWN_Psi_generic_transverse}
%
Given a subset $Z^2 \subseteq ( E \cap \crit(h) )^2$, there exists a dense subset
\[
\bPsi_h^\pitchfork(Z^2, E , \Phi+\Psi_0) \subset \bPsi_h(Z^2 , E , \Phi+\Psi_0),
\]
residual in the sense of Baire spaces, so that for any pair $(z_- , z_+) \in Z^2$ whose Morse indices satisfy $m_h(z_-) - m_h(z_+) \leq 2$,
and any $\Psi \in \bPsi_h^\pitchfork(Z^2, E , \Phi+\Psi_0)$, equation \eqref{eq:TWN} is transverse on $\calM_{\Phi,\Psi}(z_-,z_+)$.

In particular, we obtain a dense subset
\[
\bPsi_h^\pitchfork(E,\Phi+\Psi_0) \subset \bPsi_h(E,\Phi+\Psi_0),
\]
residual in the sense of Baire spaces, such that any $\Psi \in \bPsi_h^\pitchfork(E,\Phi+\Psi_0)$ is $E$-transverse up to index $2$.
\end{mythm}
\begin{myproof}
  With the density of the perturbations established by Lemma \ref{lemma:TWN_perturbation_surjective}, 
  the remainder of the proof now follows a standard transversality argument, see for example \cite{schwarz1993morse}.
  Fix $(z_- , z_+) \in Z^2$ and define the universal moduli space
  \[
  \calM_{\text{univ}}(z_-,z_+) := \set{ (u,\Psi) \in \calP(z_-,z_+) \times \bPsi_h(Z^2, E,\Phi+\Psi_0) }{ (u,\Psi) \text{ solves \eqref{eq:TWN}} }.
  \]
  We will denote by 
  \[
  \calM_{\Phi,\Psi}(z_-,z_+) := \set{ u }{ (u,\Psi) \in \calM_{\text{univ}}(z_-,z_+) }
  \]
  the moduli space of heteroclinics from $z_-$ to $z_+$ corresponding to a fixed perturbation $\Psi$.
  Finally, let
  \[
  \calM_{\Phi,\Psi}(z_-,z_+;E) := E \cap \calM_{\Phi,\Psi}(z_-,z_+)
  \]
  denote the moduli space of heteroclinics in $E$ from $z_-$ to $z_+$ corresponding to a fixed perturbation $\Psi$.
  
  First note that, in light of the gradient-like dichotomy Theorem \ref{thm:TWN_gradlike}, and since the perturbations $\Psi \in \bPsi_h(Z^2 , E,\Phi+\Psi_0)$ are $E$-tame,
  we have 
  \[
  \calM_{\Phi,\Psi}(z_-,z_+;E) = \{ z \} \qquad \text{whenever} \quad z_- = z_+ = z.
  \]
  Consequently, by Lemma \ref{lemma:hyperbolicity_implies_invertible} the moduli space $\calM_{\Phi,\Psi}(z_-,z_+;E)$ 
  satisfies the transversality condition whenever $z_- = z_+$ and $\Psi \in \bPsi_h(Z^2 , E,\Phi+\Psi_0)$.
  
  Now assume $z_- \neq z_+$.
  We note that $\calM_{\text{univ}}(z_-,z_+)$ is the zero set of the map $\Lambda(u,\Psi) = u' + \Phi(u) + \Psi(u)$.
  For $(u,\Psi) \in \calM_{\text{univ}}(z_-,z_+)$, the linearisation of $\Lambda$ around $(u,\Psi)$ takes the form
  \[
  D\Lambda(u,\Psi)[\delta u,\delta\Psi] = A[\delta u] + B[\delta \Psi],
  \]
  where
  \begin{equation*}
    \begin{array}{l l}
    A : W^{1,2}(\R,\R^d) \to L^2(\R,\R^d), \qquad & A[\delta u] = [\delta u]' + D \Phi(u)[\delta u] + D \Psi(u)[\delta u], \\
    B : \bPsi_h(Z^2) \to L^2(\R,\R^d), \qquad & B[\delta \Psi] = [\delta\Psi](u).
    \end{array}
  \end{equation*}
  By Theorem \ref{thm:TWN_Fredholm} the map $A$ is Fredholm of index $\ind(A) = m_h(z_-) - m_h(z_+)$,
  and by Lemma \ref{lemma:TWN_perturbation_surjective} the map $B$ has dense image (in fact, it suffices to restrict $B$ to the subspace $\bPsi_h(\{ z_-,z_+ \}) \subseteq \bPsi_h(Z^2)$).
  Consequently the map $D\Lambda(u,\Psi)$ has a bounded right inverse depending $C^3$ smoothly on the basepoint $(u,\Psi) \in \calM_{\text{univ}}(z_-,z_+)$,
  and we can therefore apply the implicit function theorem to derive that $\calM_{\text{univ}}(z_-,z_+)$ is a $C^3$ smooth Hilbert manifold.

  Now consider the projection
  \[
  \pi : \calM_{\text{univ}}(z_-,z_+) \to \bPsi_h(Z^2, E,\Phi+\Psi_0), \qquad (u,\Psi) \mapsto \Psi.
  \]
  By a standard transversality argument \cite{kato2012perturbation} this defines a $C^3$ smooth Fredholm map of index $\ind(\pi) = \ind(A) = m_h(z_-) - m_h(z_+)$.
  By the Sard--Smale theorem \cite{smale1965infinite} the regular values of $\pi$ are residual (i.e.\ comeager) whenever $\ind(\pi)$ is lower than the degree of smoothness of $\pi$.
  Thus, we obtain a residual subset $\bPsi_h^{z_-,z_+} \subseteq \bPsi_h(Z^2, E,\Phi+\Psi_0)$ of regular values whenever $m_h(z_-) - m_h(z_+) \leq 2$.
  For any $\Psi \in \bPsi_h^{z_-,z_+}$, it follows that $\calM_{\Phi,\Psi}(z_-,z_+) = \pi^{-1}(\Psi)$ is a $C^3$ smooth manifold.
  Since $\calM_{\Phi,\Psi}(z_-,z_+;E)$ is an isolated subset of $\calM_{\Phi,\Psi}(z_-,z_+)$,
  it follows that the transversality condition is satisfied on $\calM_{\Phi,\Psi}(z_-,z_+;E)$ whenever $\Psi \in \bPsi_h^{z_-,z_+}$.

  Now define
  \[
  \bPsi_h^\pitchfork(Z^2 , E,\Phi+\Psi_0) := \bigcap_{(z_- , z_+) \in Z^2} \bPsi_h^{z_-,z_+}.
  \]
  This being the finite intersection of residual subsets of $\bPsi_h(Z^2, E,\Phi+\Psi_0)$, we find that the set $\bPsi_h^\pitchfork(Z^2, E,\Phi+\Psi_0)$ is residual.
  By construction, for any $(z_- , z_+) \in Z^2$ with $m_h(z_-) - m_h(z_+) \leq 2$, and any $\Psi \in \bPsi_h^\pitchfork(Z^2 , E,\Phi+\Psi_0)$,
  equation \eqref{eq:TWN} is transverse on $\calM(z_-,z_+;E)$.
  In particular, if we define
  \[
  \bPhi_h^\pitchfork(E,\Phi+\Phi_0) := \bPsi_h^\pitchfork(Z^2_{\text{max}} , E,\Phi+\Phi_0), \qquad Z^2_{\text{max}} = (E \cap \crit(h))^2,
  \]
  it follows that any $\Psi \in \bPsi_h^\pitchfork(E,\Phi+\Psi_0)$ is $E$-transverse up to index $2$.
\end{myproof}


\section{Conley--Floer homology}
\label{sec:TWN_homology}

We are now prepared to construct the Conley--Floer homology groups.
We first outline the construction of the homology based on generic perturbations $\Psi \in \bPsi_h^\pitchfork(E,\Phi)$.
Then we show that, up to canonical isomophisms, the resulting homology groups are in fact independent of the particular choice of $\Psi$.
This results in the definition of Conley--Floer homology groups for the unperturbed equation \eqref{eq:TWN}.
Finally, towards the end of this section we give conditions under which the homology groups depend continuously on $\Phi$.

\subsection{Homology groups for generic data}

Assume $E$ is an isolating trajectory neighbourhood for \eqref{eq:TWN} with $\Psi = 0$.
Clearly, $\Psi=0$ is an $E$-tame perturbation, hence $\bPsi_h(E,\Phi)$ is nonempty.
We now choose a $\Psi \in \bPsi_h^\pitchfork(E,\Phi)$, which now gives us an $E$-stable, $E$-tame, $E$-transverse perturbation.

We now go though the usual steps in the construction of Floer homology.
We refer to \cite{weber2010morse, audin2014morse, schwarz1993morse} for details.
Let
\[
C_n(E,\Phi) = \set{ \sum_j \alpha_j z_j }{ \alpha_j \in \Z_2,\; z_j \in \crit(h) \cap E, \; m_h(z) = n }
\]
denote the vector space with $\Z_2$ coefficients generated by $z \in \crit(h) \cap E$ with $m_h(z) = n$.
Then introduce a linear map 
\[
\partial_n : C_n(E,\Phi) \to C_{n-1}(E,\Phi),
\]
defined on generators $z \in \crit(h)$ of $C_n(E,\Phi)$ by
\begin{equation}
  \label{eq:TWN_Floer_boundary}
  \partial_n z = \sum_{\substack{ z' \in \crit(h) \cap E \\ m_h(z')=n-1}} n(z,z') z',
\end{equation}
where $n(z,z')$ is the binary count of the number of elements of $\calM_{\Phi,\Psi}(z,z';E) / \R$.
By Theorem \ref{thm:TWN_moduli_closure} this map is well defined.
Additionally, this theorem shows that $\partial_{n}  \partial_{n+1} = 0$, since this composition encodes a binary count of $1$-fold broken chains of solutions,
which by Theorem \ref{thm:TWN_moduli_closure} always occur in pairs trapped by the isolating trajectory neighbourhood $E$.
We arrive at the following theorem.
\begin{mythm}
\label{thm:TWN_generic_Floer_homology}
  For any $\Psi \in \bPsi_h^\pitchfork(E,\Phi)$, the boundary operator $\partial_n$, defined by \eqref{eq:TWN_Floer_boundary}, satisfies the fundamental relation
  \[
  \partial_n \partial_{n+1} = 0.
  \]
  Consequently, the Conley--Floer homology groups for \eqref{eq:TWN} with generic perturbations $\Psi$
\[
\HF_n(E,\Phi,\Psi;\Z_2) := \frac{ \ker{ \partial_n } }{ \img{ \partial_{n+1} } }
\]
are well-defined.
\end{mythm}

\subsection{Continuation results}

In this section we collect various results on how the homology groups depend on the data $E$, $\Phi$, and $\Psi$.
We do so, for the most part, by studying maps between homology groups induced by homotopies of $\Phi$.
Since we have not specified a topology on the space of nonlinearities $\Phi$, we need to specify what we mean by a homotopy.

\begin{mydef}[Homotopy of $\Phi$]
\label{def:TWN_Phi_homotopy}
  Let $E \subset C^1_{\text{loc}}(\R,\R^d)$. 
  Suppose that, for each $\beta \in [0,1]$, we have a map $\Phi_\beta$ of the form
  \[
  \Phi_\beta(u) = \nabla_{g_\beta} S_\beta(u)^T \calN_\beta[S_\beta(u)] + \nabla_{g_\beta} F_\beta(u),
  \]
  satisfying the hypotheses of Section \ref{sec:TWN_setup} (with the exception of the hyperbolicity assumption in Definition \ref{def:TWN_hyperbolic}).
  Then $\Phi_\beta$ is said to be a homotopy over $E$ provided that, for any $\beta_0 \in [0,1]$, we have
   \[
  \lim_{\beta \to \beta_0} \sup_{u\in \cl E} \| \Phi_{\beta_0}(u) - \Phi_\beta(u) \|_{L^\infty(\R,\R^d)} = 0,
  \]
  and furthermore, the Riemannian metric $g_\beta$ varies continuously with $\beta$.
\end{mydef}
Note that the individual maps $\calN_\beta$, $S_\beta$, and $F_\beta$ are not assumed to depend continuously on $\beta$.
In particular, the ``intermediate dimension'' $D$ used in the definition of $\calN$ and $S$ may vary with $\beta$; this is an idea we will exploit in Section \ref{subsec:TWN_local_continuation}.
The quasi-Lyapunov functions may therefore also depend discontinuously on $\beta$.
However, the continuity assumption on $g$ ensures continuity of the kinetic energy functional,
and in particular uniformity in the estimates of the gradient-like dichotomy of Theorem \ref{thm:TWN_gradlike}.


\begin{mydef}[Relation by continuation]
  We say $\Phi_0$ and $\Phi_1$ are related by continuation over $E \subset W^{1,\infty}(\R,\R^d)$,
if there exists a homotopy $\Phi_\beta$ between $\Phi_0$ and $\Phi_1$, such that for each $\beta \in [0,1]$ the set $E$ is an isolating trajectory neighbourhood for 
\[
u' + \Phi_\beta(u) = 0.
\]
\end{mydef}

We will show that, for maps which are related by continuation, the induced homology groups are isomorphic.
To do so, we first present a local continuation result in Section \ref{subsec:TWN_local_continuation}, which shows that the Conley--Floer homology is invariant under small perturbations
of $\Phi$ and $\Psi$.
We do so by studying appropriately constructed ``continuation systems'', which are essentially higher-dimensional version of \eqref{eq:TWN} which encode the homotopy in the dynamics.
This is where the bulk of the technical work is located.
In Section \ref{subsec:TWN_HF_arbitrary} we use these results to define the homology groups for arbitrary $\Phi$, without assuming the Morse property, and with $\Psi \to 0$.
We then prove a global continuation result in Section \ref{subsec:TWN_global_continuation}.
Finally, in Section \ref{subsec:TWN_independence_E} we consider the effect of perturbing the isolating trajectory neighbourhood $E$.

\subsubsection{Local continuation} 
\label{subsec:TWN_local_continuation}

In this section, we study the effect of small perturbations $\Phi_\epsilon$ of $\Phi$.
We first specify what we mean by a ``small perturbation''.
\begin{mydef}[$\epsilon$-perturbation]
  We say $\Phi_\epsilon$ is an $\epsilon$-perturbation of $\Phi$ over $E$, provided that
  \[
  \sup_{u\in\cl E} \| \Phi_\epsilon(u) - \Phi(u) \|_{L^\infty(\R,\R^d)} < \epsilon,
  \]
  and
  \[
  \sup_{u\in\cl E} \; \sup_{\substack{ v_1 , v_2 \in \R^d \\ |v_1| = |v_2| = 1 }} \| g_{\epsilon,u}(v_1,v_2) - g_u(v_1,v_2) \|_{L^\infty(\R,\R)} < \epsilon.
  \]
\end{mydef}

Now suppose $\Phi_0$ and $\Phi_1$, given by
\[
\Phi_i(u) = \nabla_{g_i} S_i(u)^T \calN_i[S_i(u)] + \nabla_{g_i} F_i(u), \qquad i \in \{0,1\},
\]
are both $\epsilon$-perturbations of some fixed $\Phi$.
Note that the ``indermediate dimension'' $D$ used in the definition of $\calN_i$ and $S_i$ may be different for $i=0$ and $i=1$.
To study the effect of these perturbations on the homology groups, we modify the geometric approach outlined in \cite{weber2010morse}.
The idea is to construct a ``continuation system'' of the type \eqref{eq:TWN} which contains the chain complexes corresponding to $\Phi_0$ and $\Phi_1$.
The relation $\partial_n \partial_{n+1} = 0$ on the continuation system then induces a map between the Conley--Floer homologies, 
which turns out to be an isomorphism for $\epsilon > 0$ sufficiently small.

We start out with specifying an appropriate homotopy between $\Phi_0$ and $\Phi_1$.
First, for given $0 \leq \beta \leq 1$, we define 
\[
\bm{\calN} =
\begin{pmatrix}
 \calN_{0} & 0 \\
 0 & \calN_{1}
\end{pmatrix},
\qquad
\bS_\beta(u) =
\begin{pmatrix}
  \sqrt{ 1-\beta } S_0(u) \\
  \sqrt{ \beta } S_1(u)
\end{pmatrix}
\]
and
\[
F_\beta = (1-\beta) F_0 + \beta F_1, \qquad
g_\beta = (1-\beta) g_0 + \beta g_1.
\]
Then let 
\begin{equation}
  \label{eq:TWN_model_homotopy}
      \Phi_\beta(u) := \nabla_{g_\beta} \bS_\beta(u)^T \bm{\calN}[ \bS_\beta(u) ] + \nabla_{g_\beta} F_{\beta}(u).
\end{equation}
Now $\Phi_\beta$ defines a homotopy between $\Phi_0$ and $\Phi_1$, and our first result is that, for $\epsilon > 0$ sufficiently small,
this homotopy induces a relation by continuation.

\begin{mylemma}
\label{lemma:TWN_local_continuation_neighbourhood}
  Let $\Phi$ be as specified in Section \ref{sec:TWN_setup}, and let $E$ be an isolating trajectory neighbourhood for \eqref{eq:TWN}. 
  There exists $\epsilon > 0$ such that if $\Phi_0$ and $\Phi_1$ are both $\epsilon$-perturbations of $\Phi$, 
  then $\Phi_0$ and $\Phi_1$ are related by continuation over $E$, via the homotopy given by \eqref{eq:TWN_model_homotopy}.
\end{mylemma}
\begin{myproof}
  Note that \eqref{eq:TWN_model_homotopy} can be written as
  \begin{align*}
  \Phi_\beta(u) &=  (1-\beta) \Phi_0(u) + \beta \Phi_1(u) \\
   &\quad + (1-\beta) \bigg( \big(\nabla_{g_\beta} - \nabla_{g_0} \big) S_0(u) \calN_0[S_0(u)] + \big( \nabla_{g_\beta} - \nabla_{g_0} \big) F_0(u) \bigg) \\
  &\quad + \beta \bigg( \big( \nabla_{g_\beta} - \nabla_{g_1} \big) S_1(u) \calN_1[S_1(u)] + \big( \nabla_{g_\beta} - \nabla_{g_1} \big) F_1(u) \bigg).
\end{align*}
  Hence, when $g_0 = g_1$, the homotopy \eqref{eq:TWN_model_homotopy} is convex, and it follows immediately that $\| \Phi_\beta(u) - \Phi(u) \|_{L^\infty(\R,\R^d)} < \epsilon$,
  uniformly over $u \in \cl E$ and $\beta \in [0,1]$.
  When $g_0 \neq g_1$ we get additional correction terms, which are of order $\epsilon$ by the smallness assumption on the Riemannian metrics combined with the uniform boundedness of $E$.
  Thus
  \[
  \sup_{\beta \in [0,1]} \; \sup_{u \in \cl E} \| \Phi_\beta(u) - \Phi(u) \|_{L^\infty(\R,\R^d)} = O(|\epsilon|), \qquad \epsilon \to 0,
  \]
  for arbitrary $\epsilon$-perturbations $\Phi_0$ and $\Phi_1$ of $\Phi$.
  Let $\delta_0 := \inf_{ u \in \cl E \setmin E } \| u' + \Phi(u) \|_{L^\infty(\R,\R^d)}$, which is strictly positive by our assumptions.
  Choosing $\epsilon > 0$ sufficiently small, we may assume that
  \[
  \sup_{\beta \in [0,1]} \; \sup_{u \in \cl E} \| \Phi_\beta(u) - \Phi(u) \|_{L^\infty(\R,\R^d)} < \delta_0/2.
  \]
  It follows that
    \begin{align*}
    \inf_{ u \in \cl E \setmin E } \| u' + \Phi_\beta(u) \|_{L^\infty(\R,\R^d)} &\geq \inf_{ u \in \cl E \setmin E } \| u' + \Phi(u) \|_{L^\infty(\R,\R^d)} \\
      &\quad - \sup_{u \in \cl E} \| \Phi_\beta(u) - \Phi(u) \|_{L^\infty(\R,\R^d)} \geq \delta_0/2,
  \end{align*}
  uniformly in $\beta \in [0,1]$.
  This proves that $\Phi_0$ and $\Phi_1$ are related by continuation over $E$, via the homotopy given by \eqref{eq:TWN_model_homotopy}.
\end{myproof}

We now select $\epsilon$-perturbations as in Lemma \ref{lemma:TWN_local_continuation_neighbourhood}, and construct the continuation system, first without any perturbation $\Psi$.
Let $P : \R \to \R$ be the cubic polynomial given by 
\[
P(\lambda) := \lambda ( \lambda^2 - 3 ).
\]
Choose a smooth function $\beta : \R \to [0,1]$ with $\beta(\lambda) = 0$ for $\lambda \leq -1/2$ and $\beta(\lambda) = 1$ for $\lambda \geq 1/2$,
with sufficiently rapid convergence so that the composition maps $\sqrt{ \beta } : \R \to [0,1]$ and $\sqrt{ 1 - \beta } : \R \to [0,1]$ are smooth.
Now, given $\mu > 0$ (which later on will be chosen sufficiently small), consider the system of equations
\begin{equation}
  \label{eq:TWN_continuation}
  \left\{
  \begin{split}
    - u' &= \Phi_{\beta(\lambda)}(u), \\
    - \lambda' &= \mu \bigg( \frac 1 2 S_1(u) \cdot \calN_1[S_1(u)] - \frac 1 2 S_0(u) \cdot \calN_0[S_0(u)] \bigg) \beta'(\lambda) \\
    &\quad + \mu \bigg( F_1(u) - F_0(u) \bigg) \beta'(\lambda) + P'(\lambda).
  \end{split}
\right.
\end{equation}
Remark that this is, in fact, again an equation of the type \eqref{eq:TWN}, namely
\[
  \bu' + \bPhi(\bu) = 0,
\]
where
\[
\bu = 
\begin{pmatrix}
  u \\ \lambda
\end{pmatrix},
\qquad
\bPhi(\bu) = \nabla_{\bg} \bS(\bu)^T \bm{\calN}[ \bS(\bu) ] + \nabla_{\bg} \F(\bu).
\]
Here
\[
\bm{\calN} =
\begin{pmatrix}
 \calN_{0} & 0 \\
 0 & \calN_{1}
\end{pmatrix},
\qquad
\bS(\bu) = \bS_{\beta(\lambda)}(u) =
\begin{pmatrix}
  \sqrt{ 1-\beta(\lambda) } S_0(u) \\
  \sqrt{ \beta(\lambda) } S_1(u)
\end{pmatrix},
\]
and
\[
\F(\bu) = F_{\beta(\lambda)}(u) + \mu^{-1} P(\lambda), \qquad \bg = g_{\beta(\lambda)} + \mu^{-1} \d \lambda \otimes \d \lambda.
\]
Let
\[
\E := \set{ \bu = (u , \lambda) \in W^{1,\infty}(\R,\R^{d+1}) }{ 
  \begin{array}{c}
    u \in E \\
    \lambda \in W^{1,\infty}(\R,\R) \text{ with } \| \lambda\|_{L^\infty(\R,\R)} < 2
  \end{array}
}.
\]
\begin{mylemma}
\label{lemma:TWN_continuation_nbhd}
  For sufficiently small $\mu > 0$ the set $\E$ is an isolating trajectory neighbourhood for \eqref{eq:TWN_continuation}.
\end{mylemma}
\begin{myproof}
  Note that $\E = E \times B$, where
  \[
  B = \set{ \lambda \in W^{1,\infty}(\R,\R) }{ \sup_{x\in\R} | \lambda(x) | < 2}.
  \]
  As such, it is clear that $\E$ is open in $W^{1,\infty}(\R,\R^d)$ and bounded in $L^\infty(\R,\R^{d+1})$.
  Furthermore, the closure of $B$ in $C^1_{\text{loc}}(\R,\R)$ is given by
  \[
  \cl B = \set{ \lambda \in W^{1,\infty}(\R,\R) }{ \sup_{x\in\R} | \lambda(x) | \leq 2 },
  \]
  hence coincides with the closure in $W^{1,\infty}(\R,\R)$.
  Consequently, the closure of $\E$ in $C^1_{\text{loc}}(\R,\R^{d+1})$ coincides with the closure in $W^{1,\infty}(\R,\R^{d+1})$.
  Left to prove is that, for sufficiently small $\delta > 0$, the boundary
  \begin{equation}
    \label{eq:TWN_continuation_E_boundary}
    \cl{\E} \setmin \E =  \bigg( \big( \cl{E} \setmin E \big) \times \cl{B} \bigg) \cup \bigg( E \times \big( \cl{B} \setmin B \big) \bigg)
  \end{equation}
  does not contain $\delta$-approximate solutions of \eqref{eq:TWN_continuation}, with $\mu > 0$ sufficiently small.
  
  As we established in the proof of Lemma \ref{lemma:TWN_local_continuation_neighbourhood}, there exists $\delta > 0$, independent of $\beta \in [0,1]$,
 such that $\cl{E}\setmin E$ does not contain any $\delta$-approximate solutions of \eqref{eq:TWN} with $\Phi = \Phi_\beta$.
  With this choice of $\delta$ fixed, will now show that $\cl \E \setmin \E$ does not contain any $\delta$-approximate solutions for sufficiently small $\mu$.
  Pick any $\bu = (u,\lambda) \in \cl{\E} \setmin \E$.
  Suppose first that $u \in \cl{E} \setmin E$.
  Then
  \[
  \| \bu' + \bPhi(\bu) \|_{L^\infty(\R,\R^d)} \geq \inf_{0 \leq \beta \leq 1} \| u' + \Phi_\beta(u) \|_{L^\infty(\R,\R^d)} > \delta.
  \]
  This estimate is uniform in $\mu > 0$.
  Now suppose that $u \in E$.
  In light of \eqref{eq:TWN_continuation_E_boundary} we then have $\lambda \in \cl{B} \setmin B$.
  Hence there exists a sequence $(x_n)_n \subset \R$ such that 
  \[
  \lim_{n\to\infty} |\lambda(x_n)| = 2, \qquad \lim_{n\to\infty} |\lambda'(x_n)| = 0.
  \]
  Then
   \[
   \| \lambda' + P'(\lambda) \|_{L^\infty(\R,\R^d)} \geq \lim_{n\to\infty} | \lambda'(x_n) + P'(\lambda(x_n)) | = | P'(\pm 2) |.
   \]
  Now
  \[
  \| \bu' + \bPhi(\bu) \|_{L^\infty(\R,\R^d)} \geq \| \lambda' + P'(\lambda) \|_{L^\infty(\R,\R^d)} - \mu C_\E \geq | P'(\pm 2) | - \mu C_\E,
  \]
  where
  \begin{align*}
      C_\E &= \sup_{(u,\lambda) \in \cl \E} \left\| \bigg( \frac 1 2 S_1(u) \cdot \calN_1[S_1(u)] - \frac 1 2 S_0(u) \cdot \calN_0[S_0(u)] \bigg) \beta'(\lambda) \right\|_{L^\infty(\R,\R^d)} \\
    &\quad + \sup_{(u,\lambda) \in \cl \E} \left\|\bigg( F_1(u) - F_0(u) \bigg) \beta'(\lambda) \right\|_{L^\infty(\R,\R^d)}.
   \end{align*}
   Hence, whenever $\mu < ( |P'(\pm 2)| - \delta ) / C_\E$ we find that $\cl \E \setmin \E$ does not contain any $\delta$-approximate solutions of \eqref{eq:TWN}.
\end{myproof}

For each $0 \leq \beta \leq 1$, we have a function $h_\beta : \R^d \to \R$ (not necessarily Morse) associated with $\Phi_\beta$, given by
\[
h_\beta(z) = (1-\beta) h_0(z) + \beta h_1(z) = \frac 1 2 \bS_\beta(z) \cdot \bm{\calN}[ \bS_\beta(z) ] + F_\beta(z).
\]
For $\mu > 0$ we define $\bh : \R^{d+1} \to \R$ by
\[
\bh(\bz) = h_{\beta(\lambda)}(z) + \mu^{-1} P(\lambda), \qquad \bz = (z,\lambda) \in \R^d \times \R.
\]
\begin{mylemma}
\label{lemma:TWN_continuation_splitting}
  For sufficiently small $\mu > 0$, any $\bz = (z,\lambda) \in \crit(\bh)$ with $|z| \leq R$, where $R = \sup_{u \in \cl E} \| u \|_{L^\infty(\R,\R^d)}$, satisfies one of the following: 
  \begin{enumerate}
  \item $\bz \in \crit(h_0) \times \{-1\}$, and $m_{\bh}(\bz) = m_{h_0}(z) + 1$, or
  \item $\bz \in \crit(h_1) \times \{1\}$, and $m_{\bh}(\bz) = m_{h_1}(z)$.
  \end{enumerate}
  In particular, $\bh$ is a Morse function on $B_R(0) \times \R$, and 
  \[
  C_n( \E, \bPhi ) \iso C_{n-1}( E , \Phi_0 ) \oplus C_n( E , \Phi_1 ).
  \]
\end{mylemma}
\begin{myproof}
  Evaluating $\partial_\lambda \bh(z,\lambda) = 0$, we obtain the identity
  \begin{equation}
    \label{eq:Morse_mu_balance}
    \begin{split}
    - \mu^{-1} P'(\lambda) &= \bigg( \frac 1 2 S_1(z) \cdot \calN_1[S_1(z)] - \frac 1 2 S_0(z) \cdot \calN_0[S_0(z)] \bigg) \beta'(\lambda) \\
    &\quad + \bigg( F_1(z) - F_0(z) \bigg) \beta'(\lambda).
    \end{split}
   \end{equation}
   For $\lambda \not\in [-1/2,1/2]$ we have $\beta'(\lambda) = 0$, hence the right hand side vanished for those values of $\lambda$.
   However, the left hand side is nonzero unless $\lambda = \pm 1$.
   Now define $\mu_0(R) > 0$ by
   \begin{align*}
        \mu_0(R)^{-1} &:= \max_{\substack{ |z| \leq R \\ |\lambda| \leq 1/2}} \left| \frac{ \beta'(\lambda) }{ P'(\lambda) } \bigg( \frac 1 2 S_1(z) \cdot \calN_1[S_1(z)] - \frac 1 2 S_0(z) \cdot \calN_0[S_0(z)] \bigg) \right| \\
     &\quad + \max_{\substack{ |z| \leq R \\ |\lambda| \leq 1/2}} \left|  \frac{ \beta'(\lambda) }{ P'(\lambda) } \bigg( F_1(z) - F_0(z) \bigg) \right|.
   \end{align*}
   Then for any $\mu <\mu_0(R)$, equation \eqref{eq:Morse_mu_balance} does not admit any solution with $|z| \leq R$ and $\lambda \not\in \{-1,1\}$.

   The critical points of $\bh$ at $\lambda = -1$ correspond to critical points of $h_0$.
   Inspecting \eqref{eq:TWN_continuation} we find there is one additional unstable direction, corresponding to $\partial_\lambda$,
   which results in a shift in the Morse index.
   At $\lambda = 1$ the critical points correspond to the critical points of $h_1$, and there are no additional unstable directions, hence the Morse index remains unchanged.
\end{myproof}

Thus far we have constructed the continuation system without perturbation $\Psi$.
Since the space $\bPsi_{h_\beta}$ is only defined when $h_\beta$ is Morse, there is no clear-cut way to incorporate a homotopy of $\Psi$ in \eqref{eq:TWN_model_homotopy}.
Nevertheless, now that a continuation system without perturbations is defined, we can construct perturbations of the continuation system which have the properties required by the continuation theorem.
\begin{mylemma}
  \label{lemma:TWN_continuation_perturbation}
  For sufficiently small $\mu > 0$, there exists $\delta > 0$ such that the following holds.
  Suppose $\Psi_i \in \bPsi_{h_i}^\pitchfork(E,\Phi_i)$, with $\| \Psi_i \|_{\bPsi_{h_i}} < \delta$, where $i \in \{0,1\}$.
Then there exists $\bPsi \in \bPsi_{\bh}^\pitchfork(\E,\bPhi)$ such that, concerning the perturbed continuation system
  \begin{equation}
    \label{eq:TWN_continuation_perturbed}
    \bu' + \bPhi(\bu) + \bPsi(\bu) = 0,
  \end{equation}
  the following holds.
  \begin{enumerate}
  \item If $\bu = (u,\lambda) \in \E$ is a heteroclinic solution of \eqref{eq:TWN_continuation_perturbed}, then either $\lambda$ is contant (and equal to $-1$ or $1$),
    or $\lambda(x) \to \pm 1$ as $x\to\pm\infty$.
  \item There is a 1-to-1 correspondence between solutions $\bu = (u,\lambda) \in \E$ of \eqref{eq:TWN_continuation_perturbed} with $\lambda = - 1$, or $\lambda = 1$,
   and solutions $u\in E$ of
   \[
   u' + \Phi_0(u) + \Psi_0(u) = 0, \qquad \text{or} \quad u' + \Phi_1(u) + \Psi_1(u) = 0, 
   \]
   respectively.
  \end{enumerate}
\end{mylemma}
\begin{myproof}
  We will first construct a perturbation $\bPsi \in \bPsi_{\bh}(\E,\bPhi)$.
  We then verify that properties (1) and (2) hold.
  Finally, we verify that the constructed perturbation $\bPsi$ is indeed $\E$-transverse up to index $2$.
  
  The perturbations $\Psi_0$ and $\Psi_1$ are of the form
  \[
    \rst{ \Psi_i(u) }{x} = \sum_{\theta \in \Theta_{h_i}} \rst{ \alpha_\theta^i * \psi_\theta(u) }{ (x,x) }, \qquad i \in \{0,1\},
  \]
  where $\psi_{\bm\theta}$ is as defined in Section \ref{subsec:TWN_perturbations}, and $\alpha_\theta^i : \R^2 \to \R^d$.
  Now define $\bPsi_0 \in \bPsi_{\bh}$ by 
  \[
  \rst{ \bPsi_0(\bu) }{x} = \sum_{\bm\theta \in \bm\Theta_{\bh}} \rst{ \bm\alpha_{\bm\theta} * \bm\psi_{\bm\theta}(\bu) }{(x,x)},
  \]
  where
  \[
   \bm\psi_{\bm\theta}(\bu) := \psi_{(\ell,z_1,z_2)}(u), \qquad \text{when} \quad \bm\theta = (\ell,(z_1,\lambda_1),(z_2,\lambda_2)), \quad \bu = (u,\lambda),
  \]
  and $\bm\alpha_{\bm\theta} : \R^2 \to \R^{d+1}$ is given by
  \[
  \bm\alpha_{\bm\theta} := \left\{
    \begin{array}{l l}
      \big( \alpha^0_{(\ell , z_1 , z_2)} , 0 \big) \quad &  \text{when} \quad \bm\theta = ( \ell , (z_1,-1) , (z_2,-1) ), \\ 
      \vspace{-3ex} & \\
      \big( \alpha^1_{(\ell , z_1 , z_2)} , 0 \big) &  \text{when} \quad \bm\theta = ( \ell , (z_1,1) , (z_2,1) ), \\ 
      0 & \text{otherwise}.
    \end{array}
    \right.
  \]
  Note that $\| \bPsi_0 \|_{\bPsi_{\bh}} = \| \Psi_0 \|_{\bPsi_{h_0}} + | \Psi_1 \|_{\bPsi_{h_1}} < 2 \delta$.
  Hence, for sufficiently small $\delta$, it follows from Theorem \ref{thm:E-stable} and Theorem \ref{thm:E-tame} that $\bPsi_0 \in \bPsi_{\bh}(\E,\bPhi)$, 
  i.e.\  $\bPsi_0$ is $\E$-stable and $\E$-tame.
  Let
  \[
  Z^2 := \set{ (\bz_- , \bz_+) \in ( \E \cap \crit(\bh))^2 }{ \bz_- = ( z_- , -1), \quad \bz_+ = ( z_+ , 1 ) },
  \]
  and invoke Theorem \ref{thm:TWN_Psi_generic_transverse} to obtain a perturbation $\bPsi' \in \bPsi_\bh^\pitchfork( Z^2 , E , \bPhi + \bPsi_0)$, which we may choose arbitrarily small.
  Explicitly, $\bPsi'$ is of the form
  \[
      \rst{ \bPsi'(\bu) }{x} = \sum_{\bm\theta \in \bm\Theta_{\bh}} \rst{ \bm\beta_{\bm\theta} * \bm\psi_{\bm\theta}(\bu) }{(x,x)},
  \]
  where $\psi_{\bm\theta}$ is an arbitrary basic perturbation as defined in Section \ref{subsec:TWN_perturbations}, 
  and $\bm\beta_{\bm\theta} : \R^2 \to \R^{d+1}$ satisfies
  \[
  \bm\beta_{\bm\theta} = 0 \qquad \text{when} \quad \bm\theta = (\ell , (z_- , -1) , (z_+, -1) ) \quad \text{or} \quad \bm\theta = (\ell , (z_- , 1) , (z_+,1) ).
  \]
  Hence for any $\bm\theta \in \bm\Theta_{\bh}$, we have either $\bm\alpha_{\bm\theta} = 0$, or $\bm\beta_{\bm\theta} = 0$ (and possibly $\bm\alpha_{\bm\theta} = \bm\beta_{\bm\theta} = 0$).
  Now define $\bPsi := \bPsi_0 + \bPsi'$.
  
  We still need to specify $\mu > 0$, for which we require some auxiliary definitions.
  Define $V : L^\infty(\R, \R^d \times \R) \to L^\infty(\R,\R)$ by
  \begin{align*}
    V(u,\lambda) &:= \bigg( \frac 1 2 S_1(u) \cdot \calN_1[S_1(u)] - \frac 1 2 S_0(u) \cdot \calN_0[S_0(u)] \bigg) \beta'(\lambda) \\
    &\quad + \bigg( F_1(u) - F_0(u) \bigg) \beta'(\lambda),
  \end{align*}
  so that for solutions $\bu = (u,\lambda) \in \E$ of the unperturbed continuation system \eqref{eq:TWN_continuation}, we have $-\lambda' = \mu V(u,\lambda) + P'(\lambda)$.
  Since $V$ is bounded on $\E$, and since $P'(s) < 0$ for $-1/2 \leq s \leq 1/2$,
  there exists $L > 0$ so that for $\mu > 0$ sufficiently small we have
  \[
  \mu V(u,\lambda) + P'(s) \leq L, \qquad (u,\lambda) \in \E, \; s \in [-1 / 2 , 1 / 2].
  \]
  This established the choice of $\mu > 0$.
  Note that we may assume that $\| \bPsi' \|_{\bPsi_{\bh}} < L / 2$.

  We now prove that property (1) holds.
  Suppose $\bu = (u,\lambda) \in \E$ is a heteroclinic solution of \eqref{eq:TWN_continuation_perturbed}.
  We note here that $\pi_{d+1} \bPsi_0 = 0$, where $\pi_i : \R^{d+1} \to \R$ denotes the orthogonal projection onto the $i$-th component,
  so that
  \[
  - \lambda' = \mu V(u,\lambda) + P'(\lambda) + \pi_{d+1} \bPsi'(u,\lambda).
  \]
  Suppose $x_0 \in \R$ is such that $| \lambda(x_0) | < 1 / 2$.
  Then our choices of $\mu$ and $\bPsi'$ imply that
  \[
  - \lambda'(x_0) = \mu V(u,\lambda)(x_0) + P'(\lambda(x_0)) + \pi_{d+1} \bPsi'(u,\lambda)(x_0) \leq L / 2 < 0.
  \]
  Since $\bu$ is assumed to be a heteroclinic, it follows from Lemma \ref{lemma:TWN_continuation_splitting} that $\lambda(x)$ must accumulate onto $\{ -1 , 1\}$
  as $x\to\pm\infty$.
  Hence, if $| \lambda(x_0) | < 1 / 2$ for some $x_0 \in \R$, it must be that $\lambda(x) \to \pm 1$ as $x\to\pm\infty$.
  
  Now assume that $|\lambda(x)| \geq 1 / 2$ for all $x \in \R$.
  Then either $\lambda(x) \to -1$ as $x\to\pm\infty$, or $\lambda(x) \to 1$ as $x\to\pm\infty$, i.e., 
  $\lambda$ is constant or homoclinic to either $-1$ or $1$.
  Note that $V(u,\lambda) = 0$, since by construction $\beta(\lambda)$ is constant for $|\lambda| \geq 1/2$.
  Inspecting the construction of the maps
\[
\rst{ \bm\psi_{\bm\theta}(\bu) }{(\xi_1,\xi_2)} = \rst{ \bm\sigma_{\ell,\bz_1}(\bu) }{\xi_1} \rst{ \bm\sigma_{\ell,\bz_2}(\bu) }{\xi_2}, \qquad \bu = (u,\lambda)
\]
in Section \ref{sec:TWN_setup}, we find that
\begin{align*}
  \sigma_{\ell,\bz_1}(u,\lambda) &= 0, \qquad \text{if} \quad \inf_{x\in\R} \lambda(x) > 1/2, \\
  \sigma_{\ell,\bz_2}(u,\lambda) &= 0, \qquad \text{if} \quad \sup_{x\in\R} \lambda(x) < -1/2,
\end{align*}
whenever $\bz_1 = (z_1,-1)$ and $\bz_2 = (z_2,1)$, with $z_1 \in \crit(h_0)$ and $z_2 \in \crit(h_1)$.
  We thus find that, whenever $\| \lambda \|_{L^\infty(\R,\R)} \geq 1 / 2$ and $\mu > 0$ is sufficiently small, we have
  \[
  \bm\psi_{\bm\theta}(u,\lambda) = 0 \qquad \text{when} \quad \bm\theta \neq (\ell , (z_- , \pm 1) , (z_+,\pm 1) ).
  \]
  From this we find that $\bPsi'(u,\lambda) = 0$ for such $\lambda$.
  Hence $-\lambda' = P'(\lambda)$, which excludes the possibility of $\lambda$ being homoclinic to either $-1$ or $1$.
  We conclude that $\lambda = -1$ or $\lambda = 1$, which verifies property (1).
  
  We now verify property (2) holds.
  Again inspecting the construction of the maps $\psi_\theta$ and $\bm\psi_{\bm\theta}$ in Section \ref{sec:TWN_setup}, 
  we see that, with $\mu > 0$ sufficiently small, for any $u \in E$ we have
  \[
  \bm\psi_{\bm\theta}(u,\pm 1) = 
  \left\{
    \begin{array}{l l}
      \psi_{(\ell , z_1, z_2)}(u) \quad &  \text{when} \quad \bm\theta = ( \ell , (z_1,\pm 1) , (z_2,\pm 1) ), \\
      0 & \text{otherwise}.
    \end{array}
    \right.
  \]
  This results in
  \[
  \bPsi(u,-1) = ( \Psi_0(u) , 0 ), \qquad \bPsi(u,1) = ( \Psi_1(u) , 0 ).
  \]
  From this, property (2) follows.

  Left to verify is that this perturbation $\bPsi$ is $\E$-transverse up to index $2$.
  Suppose that $\bu = (u,\lambda)$ is a solution of \eqref{eq:TWN_continuation_perturbed} with $\big( \bu(x) , \bu'(x) \big) \to \big( \bz_\pm , 0 \big)$ as $x \to \pm\infty$.
  When $\bz_- = ( z_- , -1 )$ and $\bz_+ = ( z_+ , 1)$, with $m_{\bh}(\bz_-) - m_{\bh}(\bz_+) \leq 2$,
  it follows by the construction of $\bPsi$ (and particularly, the choice of $\bPsi'$) that \eqref{eq:TWN} is transverse on $\calM_{\bPhi,\bPsi}(\bz_-,\bz_+;\E)$.
  Now consider the case where $\bz_\pm = ( z_\pm , -1 )$, with $z_\pm \in \crit( h_0 )$.
  As we have just seen, it now follows that $\lambda = -1$, and $\bPsi(\bu) = ( \Psi_0(u) , 0 )$.
  The linearisation of \eqref{eq:TWN_continuation_perturbed} around $\bu \in \calP(\bz_-,\bz_+)$ now attains the form
  \[
  \partial_x +
  \begin{pmatrix}
    D \Phi_0(u) + D \Psi_0(u)  & 0 \\
    0 & P''(-1)
  \end{pmatrix}
  : W^{1,2}(\R,\R^d \oplus \R) \to L^2(\R,\R^d \oplus \R),
  \]
  which is surjective whenever $m_{\bh}(\bz_-) - m_{\bh}(\bz_+) = m_{h_0}(z_-) - m_{h_0}(z_+) \leq 2$,
  since we assume $\Psi_0 \in \bPsi_{h_0}^\pitchfork(E,\Phi_0)$ and $-1$ is a nondegenerate critical point of $P$.
  Likewise, when $\bz_\pm = ( z_\pm ,  1 )$ with $z_\pm \in \crit(h_1)$, 
  surjectivity of the linearisation follows from $E$-transversality of $\Psi_1$ and nondegeneracy of $1 \in \crit(P)$.
  We conclude that $\bPsi \in \bPsi_{\bh}^\pitchfork(\E , \bPhi)$.
  This completes the proof.
\end{myproof}


This establishes the continuation system.
We are now prepared to prove the local continuation property of $\HF_*$.
\begin{mythm}
\label{thm:TWN_local_continuation}
  Let $\Phi$ be as specified in Section \ref{sec:TWN_setup}, and let $E$ be an isolating trajectory neighbourhood for \eqref{eq:TWN}.
  Then there exist $\delta > 0$ and $\epsilon > 0$ such that the following holds.
  Suppose $\Phi_0$ and $\Phi_1$ are $\epsilon$-perturbations of $\Phi$.
  Choose arbitrary $\Psi_i \in \bPsi_{h_i}^\pitchfork(E,\Phi_i)$ with $\| \Psi_i \|_{\bPsi_{h_i}} < \delta$, for $i \in \{0,1\}$.
  Then there exist canonical isomorphisms
  \[
  \Theta^{1,0}_n : \HF_n(E,\Phi_0,\Psi_0;\Z_2) \to \HF_n(E,\Phi_1,\Psi_1;\Z_2), \qquad n \in \N.
  \]
  Moreover, suppose $\Phi_0$, $\Phi_1$, and $\Phi_2$ are all $\epsilon$-perturbations of $\Phi$, and select perturbations $\Psi_i \in \bPsi_{h_i}^\pitchfork(E,\Phi_i)$ with $\| \Psi_i \|_{\bPsi_{h_i}} < \delta$,
  for $i \in \{0,1,2\}$.
  Then the induced maps on homology satisfy the transitive relation
  \[
  \Theta^{2,0}_n = \Theta^{2,1}_n \Theta^{1,0}_n.
  \]
\end{mythm}
\begin{myproof}
  Now that the construction of a continuation system for \eqref{eq:TWN} is established, we can follow the geometric idea outlined in \cite{weber2010morse}.
  Let $\mu > 0$ be sufficiently small, select $\epsilon > 0$ as specified by Lemma \ref{lemma:TWN_local_continuation_neighbourhood}, 
  and choose $\delta > 0$ as specified by Lemma  \ref{lemma:TWN_continuation_perturbation}.
  Then $\Phi_0$ and $\Phi_1$ are related by continuation over $E$, so that we obtain a corresponding continuation system.
  Choose a perturbation $\bPsi \in \bPsi_h^\pitchfork(\E,\bPhi)$ as specified by Lemma \ref{lemma:TWN_continuation_perturbation}, so that the perturbed continuation system
  \[
    \bu' + \bPhi(\bu) + \bPsi(\bu) = 0
  \]
  is $\E$-transverse up to index $2$.
  Hence the Conley--Floer homology groups are well-defined for the perturbed continuation system.
  In particular, we now obtain boundary operators $\Delta^{1,0}_n : C_n(\E,\Phi) \to C_{n-1}(\E,\Phi)$ which satisfy the fundamental relation $\Delta^{1,0}_n  \Delta^{1,0}_{n+1} = 0$.

  With regards to the splitting of the critical groups obtained in Lemma \ref{lemma:TWN_continuation_splitting}, the boundary operator $\Delta_n^{1,0}$ is now given by
  \[
  \Delta^{1,0}_n =
  \begin{pmatrix}
    \partial^0_{n-1} & 0 \\
    \Theta^{1,0}_{n-1} & \partial^1_n
  \end{pmatrix}.
  \]
  By Lemma \ref{lemma:TWN_continuation_perturbation} there are no heteroclinic connections from $(z,1)$ to $(z',-1)$,
  which is why the upper right corner of the matrix $\Delta^{1,0}_n$ is $0$.
  The map $\partial^0_n$ counts heteroclinic connections from $(z,-1)$ to $(z',-1)$ with $m_{h_0}(z) = m_{h_0}(z') + 1$.
  Likewise, $\partial^1_n$ counts heteroclinic connections from $(z,1)$ to $(z',1)$ with $m_{h_1}(z) = m_{h_1}(z') + 1$.
  In light of Lemma \ref{lemma:TWN_continuation_perturbation}, these maps $\partial^i_n$ are nothing but the boundary operators corresponding to
  \[
  u' + \Phi_i(u) + \Psi_i(u) = 0, \qquad i \in \{0,1\}.
  \]
  Finally, the map $\Theta^{1,0}_*$ counts isolated heteroclinic connections from $(z,-1)$ to $(z',1)$ with $m_{h_0}(z) = m_{h_1}(z')$.
  From the fundamental relation $\Delta^{1,0}_n  \Delta^{1,0}_{n+1} = 0$ we derive that $\Theta^{1,0}_*$ is a chain homomorphism, i.e.,
  \[
  \Theta^{1,0}_{n-1}  \partial^0_n = \partial^1_n  \Theta^{1,0}_n.
  \]
  Hence $\Theta^{1,0}_*$ induces a homomorphism between the homology groups, which we again denote by the same symbol,
  \[
  \Theta^{1,0}_n : \HF_n(E,\Phi_0,\Psi_0;\Z_2) \to \HF_n(E,\Phi_1,\Psi_1;\Z_2).
  \]
  As we already know that the Conley--Floer homology groups are independent of the chosen perturbations $\Psi$,
  it now suffices to check that $\Theta^{1,0}_*$ are in fact isomorphisms.

  To see that the maps $\Theta^{1,0}_*$ are in fact isomorphisms we now study their dependence on the maps $\Phi_i$.
  Suppose $\Phi_0$, $\Phi_1$, $\Phi_2$ and $\Phi_3$ are all $\epsilon$-perturbations of $\Phi$.  
  For any $j, k \in \{0,\dots,3\}$, let us denote by 
  \[
  \bu' + \bPhi^{k,j}(\bu) + \bPsi^{k,j}(\bu) = 0
  \]
  the perturbed continuation system induced by the homotopy \eqref{eq:TWN_model_homotopy} from $\Phi_j$ to $\Phi_k$.
  Repeating the preceding construction, for any $j,k \in \{0,\dots,3\}$ we obtain a map
  \[
  \Theta^{k,j}_n : C_n(E,\Phi_j) \to C_n(E,\Phi_k).
  \]
  Furthermore, repeating the estimates in Lemma \ref{lemma:TWN_local_continuation_neighbourhood}, decreasing $\epsilon > 0$ and $\mu > 0$ if needed,
  we find that $\bPhi^{1,0}$ and $\bPhi^{3,2}$ are related by continuation over $\E$.
  Correspondingly, we obtain a higher dimensional continuation system, for which we again obtain a boundary operator $\bm{\Delta}_n : \C_n \to \C_{n-1}$.
  With respect to the splitting
  \[
    \C_n \iso C_{n-2}(E,\Phi_0) \oplus C_{n-1}(E,\Phi_1) \oplus C_{n-1}(E,\Phi_2) \oplus C_n(E,\Phi_3)
  \]
  the boundary operator $\bm{\Delta}_n$ is represented by
  \[
  \bm{\Delta}_n =
  \begin{pmatrix}
        \partial^0_{n-2} & 0 & 0 & 0 \\
    \Theta^{1,0}_{n-2} & \partial^1_{n-1} & 0 & 0 \\
    \Theta^{2,0}_{n-2} & 0 & \partial^2_{n-1} & 0 \\
    \Lambda^{3,0}_{n-2} & \Theta^{3,1}_{n-1} & \Theta^{3,2}_{n-1} & \partial^3_n
  \end{pmatrix},
  \]
  where $\Lambda^{3,0}_{n-2}$ counts isolated heteroclinic connections from $(z,-1,-1)$ to $(z',1,1)$, with $m_{h_0}(z) = m_{h_3}(z') - 1$.
  Evaluating the relation $\bm{\Delta}_n \bm{\Delta}_{n+1} = 0$ results in the chain homotopy
  \[
  \Lambda^{3,0}_{n-2} \partial^0_{n-1} + \partial^3_{n} \Lambda^{3,0}_{n-1} = \Theta^{3,2}_{n-1} \Theta^{2,0}_{n-1} - \Theta^{3,1}_{n-1} \Theta^{1,0}_{n-1}.
  \]
  This implies that on the level of homology we have 
  \begin{equation}
    \label{eq:TWN_Theta_functorial}
      \Theta^{3,2}_* \Theta^{2,0}_* = \Theta^{3,1}_* \Theta^{1,0}_*.
  \end{equation}
  If $\Phi_2 = \Phi_3$ this yields the functorial relation
  \[
  \Theta^{2,0}_* = \Theta^{2,1}_* \Theta^{1,0}_*,
  \]
  where we used that $\Theta^{2,2}_* = \id$ on the level of chain groups.
  Letting $\Phi_2 = \Phi_0$ we then see that $\Theta^{1,0}_*$ is an isomorphism on the level of homology, with inverse given by $\Theta^{0,1}_*$.

  Finally, remark that the construction of $\Theta^{1,0}_*$ required the choice of a $\beta : \R \to [0,1]$ and $\mu > 0$,
  and strictly speaking we should denote this map by $\Theta^{1,0}_*(\beta,\mu)$.
  For different choices $\widetilde \beta : \R \to [0,1]$ and $\widetilde \mu > 0$, we obtain a map $\Theta^{1,0}_*(\widetilde \beta , \widetilde \mu)$, 
  which on the level of chain groups may be distinct from $\Theta^{1,0}_*(\beta,\mu)$.
  From \eqref{eq:TWN_Theta_functorial}, with $\Theta^{1,0}_* = \Theta^{1,0}_*(\beta,\mu)$, and $\Theta^{2,0}_* = \Theta^{1,0}_*(\widetilde \beta , \widetilde \mu)$,
  and $\Phi_1 = \Phi_2 = \Phi_3$, we find that on the level of homology groups
  \[
  \Theta^{1,0}_*(\beta,\mu) = \Theta^{1,0}_*(\widetilde \beta , \widetilde \mu).
  \]
  Hence on the level of homology groups, the map $\Theta^{1,0}_*$ is canonically defined by the ``endpoints'' $( \Phi_0 , \Psi_0 )$ and $( \Phi_1 , \Psi_2 )$.
\end{myproof}


\subsubsection{Homology groups for arbitrary data}
\label{subsec:TWN_HF_arbitrary}

We may now define the Conley--Floer homology groups for nonlinearities $\Phi$ for which the constant solutions of \eqref{eq:TWN} are degenerate (i.e.\ not hyperbolic).
For this we use the following stability result for isolating trajectory neigbhourhoods.
\begin{mylemma}
  \label{lemma:TWN_isolating_traj_nbhd_Phi_stable}
  Let $\Phi$ be as described in Section \ref{sec:TWN_setup}, with the exception of the Morse property imposed in Section \ref{subsec:Morse_function_h}.
  Suppose $E$ is an isolating trajectory neighbourhood for \eqref{eq:TWN}.
  For sufficiently small $\epsilon > 0$ there exists an $\epsilon$-perturbation $\Phi_\epsilon$ of $\Phi$, for which the following properties hold:
  \begin{enumerate}
  \item The set $E$ is an isolating trajectory neighbourhood for \eqref{eq:TWN} with $\Phi = \Phi_\epsilon$.
  \item All constant solutions $z \in E$ of \eqref{eq:TWN} with $\Phi = \Phi_\epsilon$ are hyperbolic.
  \end{enumerate}
  A perturbation which satisfies these properties is said to be stable and regular with respect to $E$.
\end{mylemma}
\begin{myproof}
  Remark that, by Thom transversality, given $\epsilon > 0$, the local potential $F$ appearing in 
\[
\Phi(u) = \nabla_g S(u)^T \calN[S(u)] + \nabla_g F(u)
\]
may be perturbed slightly into $F_\epsilon$, with 
\[
\| F_\epsilon - F \|_{L^\infty(\R^d,\R)} + \| \nabla_g F_\epsilon - \nabla_g F \|_{L^\infty(\R^d,\R^d)} < \epsilon,
\]
to ensure all constant solutions of \eqref{eq:TWN} which lie in $E$ are hyperbolic.
This results in a perturbed nonlinearity $\Phi_\epsilon$, defined by
\[
\Phi_\epsilon(u) = \nabla_g S(u)^T \calN[S(u)] + \nabla_g F_\epsilon(u).
\]
Clearly this is an $\epsilon$-perturbation of $\Phi$.
Property (2) follow immediately from the definition of $\Phi_\epsilon$.
Since $E$ is an isolating trajectory neighbourhood for the unperturbed equation \eqref{eq:TWN}, there exists $\delta > 0$ so that $\cl{E} \setmin E$
does not contain any $\delta$-approximate solutions of the unperturbed equation \eqref{eq:TWN}.
If we let $0 < \epsilon < \delta / 2$, it follows from property (1) that $\cl{E} \setmin E$ does not contain any $\delta/2$-approximate solutions of \eqref{eq:TWN} with $\Phi = \Phi_\epsilon$.
Hence for sufficiently small $\epsilon > 0$, property (1) is also satisfied.
\end{myproof}

As a consequence of the preceding lemma, for a given stable and regular $\epsilon$-perturbation $\Phi_\epsilon$, and $\Psi_\epsilon \in \bPsi_{h_\epsilon}^\pitchfork(E,\Phi_\epsilon)$,
the homology $\HF_*(E,\Phi_\epsilon,\Psi_\epsilon;\Z_2)$ is well-defined for sufficiently small $\epsilon > 0$.
The local continuation property of $\HF_*$ ensures that, up to canonical isomorphism, this homology is independent of the choice of $\Phi_\epsilon$ and $\Psi_\epsilon$.
We then define $\HF_*(E,\Phi;\Z_2)$ as the isomorphism type of such perturbations.
\begin{mydef}[Conley--Floer homology with arbitrary data]
  Suppose $\Phi$ satisfies the Hypotheses from Section \ref{sec:TWN_setup}, but $E$ contains constant solutions of \eqref{eq:TWN} which are not hyperbolic.
  Given an isolating trajectory neighbourhood $E$, the Conley--Floer homology groups for the pair $(E,\Phi)$ are defined by
  \[
  \HF_*(E,\Phi;\Z_2) := \varprojlim \HF_*(E,\Phi_\epsilon,\Psi_\epsilon;\Z_2),
  \]
  where the inverse limit is taken over all perturbations $\Phi_\epsilon$ of $\Phi$ satisfying properties (1)--(3) from Lemma \ref{lemma:TWN_isolating_traj_nbhd_Phi_stable},
  and all $\Psi_\epsilon \in \bPsi_{h_\epsilon}^{\pitchfork}(E,\Phi)$.
\end{mydef}

\subsubsection{Global continuation}
\label{subsec:TWN_global_continuation}

We now arrive at a global continuation theorem for the Conley--Floer homology.
\begin{mythm}
  \label{thm:TWN_HF_homotopy_invariance}
  Suppose $\Phi_\beta$ with $0 \leq \beta \leq 1$ is a homotopy over $E$,
  which defines a relation by continuation over $E$.
  Then the homotopy induces an isomorphism
  \[
  \HF_*( E , \Phi_0 ; \Z_2 ) \iso \HF_*( E , \Phi_1 ; \Z_2 ).
  \]
\end{mythm}
\begin{myproof}
  Given $\beta \in [0,1]$, denote by $\delta(\Phi_\beta) > 0$ and $\epsilon(\Phi_\beta) > 0$ the values of $\delta > 0$ and $\epsilon > 0$ 
prescribed by Theorem \ref{thm:TWN_local_continuation} with $\Phi = \Phi_\beta$.
  By compactness we may select $\beta_1 , \dots , \beta_n \in [0,1]$, with $\beta_1 = 0$ and $\beta_n = 1$,
  so that $\Phi_{\beta_{k+1}}$ is an $\epsilon(\Phi_{\beta_k})$-perturbation of $\Phi_{\beta_k}$, for $k \in \{1,\dots,n-1\}$, and vice versa.
Analogous to Lemma \ref{lemma:TWN_isolating_traj_nbhd_Phi_stable}, we may slightly perturb the homotopy
so that the maps $\Phi_{\beta_1} , \dots , \Phi_{\beta_n}$ are regular over $E$, while preserving the isolating property of $E$.
Since, up to canonical isomorphisms, the homologies are independent of the particular choice of such a perturbation, 
we again denote this perturbed homotopy by $\Phi_\beta$.
Let $\delta := \min\{ \delta(\Phi_{\beta_1}) , \dots , \delta(\Phi_{\beta_n}) \}$, and for each $k \in \{1,\dots,n\}$, choose a perturbation 
$\Psi_k \in \bPsi_{h_{\beta_k}}^\pitchfork(E,\Phi_{\beta_k})$ with $\| \Psi_k \|_{\bPsi_{h_{\beta_k}}} < \delta$.
We can now invoke Theorem \ref{thm:TWN_local_continuation} to find that
\[
\HF_*(E,\Phi_{\beta_k},\Psi_k;\Z_2) \iso \HF_*(E,\Phi_{\beta_{k+1}},\Psi_{k+1};\Z_2), \qquad k \in \{1,\dots,n-1\}.
\]
Consequently, we find that
\begin{align*}
  \HF_*(E,\Phi_0;\Z_2) &\iso \HF_*(E,\Phi_{\beta_1} , \Psi_{1};\Z_2) \iso \HF_*(E,\Phi_{\beta_2} , \Psi_2;\Z_2) \iso \cdots \\
  &\iso \HF_*(E,\Phi_{\beta_{n-1}},\Psi_{n-1};\Z_2) \iso \HF_*(E,\Phi_{\beta_n},\Psi_{n};\Z_2) \\
  &\iso \HF_*(E,\Phi_1;\Z_2).
\end{align*}
The transitive relation of the isomorphisms obtained in Theorem \ref{thm:TWN_local_continuation} ensures that this isomorphism between $\HF_*(E,\Phi_0;\Z_2)$ and $\HF_*(E,\Phi_1;\Z_2)$
is independent of the particular choice of intermediate points $\beta_1,\dots,\beta_n \in [0,1]$ and the corresponding generic perturbations.
\end{myproof}

\subsubsection{Independence of $E$}
\label{subsec:TWN_independence_E}

Clearly, the construction of the homology for generic data only depends on the invariant set $\Inv(E)$, but not on the particular choice of $E$.
This remains to when considering the homology for nongeneric data.
\begin{mylemma}
\label{lemma:TWN_independence_isolating_nbhd}
  Suppose that $E_1$ and $E_2$ are isolating trajectory neighbourhoods, with $\Inv(E_1) = \Inv(E_2)$.
  Then there is a canonical isomorphism
  \[
  \HF_*(E_1 ,\Phi;\Z_2) \iso \HF_*(E_2,\Phi;\Z_2).
  \]
\end{mylemma}
\begin{myproof}
  We first choose an appropriate perturbed nonlinearity $\Phi_\epsilon$, stable and regular with respect to both $E_1$ and $E_2$,
  as in Lemma \ref{lemma:TWN_isolating_traj_nbhd_Phi_stable}.
  We note that, as $\Inv(E_1) = \Inv(E_2)$, the intersection $\bPsi_h^\pitchfork(E_1,\Phi_\epsilon) \cap \bPsi_h^\pitchfork(E_2,\Phi_\epsilon)$ is nonempty.
  Select a generic perturbation $\Psi_\epsilon \in \bPsi_h^\pitchfork(E_1,\Phi_\epsilon) \cap \bPsi_h^\pitchfork(E_2,\Phi_\epsilon)$ with $\| \Psi_\epsilon \|_{\bPsi_h} < \epsilon$.
  Now, by definition of the Conley--Floer homology for nongeneric data, we have canonical isomorphisms
  \[
  \HF_*(E,\Phi;\Z_2) \iso \HF_*(E,\Phi_\epsilon,\Psi_\epsilon;\Z_2), \qquad \text{for} \quad E=E_1, \; E_2.
  \]
  Moreover, we have
  \[
  \HF_*(E_1,\Phi_\epsilon,\Psi_\epsilon;\Z_2) = \HF_*(E_2,\Phi_\epsilon,\Psi_\epsilon;\Z_2),
  \]
  where the equality holds on the level of chain complexes.
  As this holds for any such generic perturbations $\Phi_\epsilon$ and $\Psi_\epsilon$, we obtain the canonical identification of the homology groups for nongeneric data.
\end{myproof}

\subsection{Direct sum formula}

We now state an algebraic relation between homology groups induced by set relations on isolating trajectory neighbourhoods.

\begin{mylemma}
\label{lemma:TWN_HF_direct_sum}
  Suppose that $E_1$ and $E_2$ are disjoint isolating trajectory neighbourhoods, $E_1 \cap E_2 = \emptyset$.
  Then $E_1 \cup E_2$ is an isolating trajectory neighbourhood, and
  \[
  \HF_*(E_1 \cup E_2,\Phi;\Z_2) \iso \HF_*(E_1,\Phi;\Z_2) \oplus \HF_*(E_2,\Phi;\Z_2).
  \]
\end{mylemma}
\begin{myproof}
  The fact that $E_1 \cup E_2$ is another isolating trajectory neighbourhood follows directly from the disjointness of $E_1$ and $E_2$.
  To prove the splitting of the homology, we first choose an appropriate perturbed nonlinearity $\Phi_\epsilon$, stable and regular with respect to $E_1\cup E_2$,
  as in Lemma \ref{lemma:TWN_isolating_traj_nbhd_Phi_stable}.
  Then select a generic perturbation $\Psi_\epsilon \in \bPsi_h^\pitchfork(E_1 \cup E_2,\Phi_\epsilon)$ with $\| \Psi_\epsilon \|_{\bPsi_h} < \epsilon$.
  Then, by disjointness of $E_1$ and $E_2$, we also have $\Psi_\epsilon \in \bPsi_h^\pitchfork(E_1,\Phi_\epsilon) \cap \bPsi_h^\pitchfork(E_2,\Phi_\epsilon)$.
  Now, by definition of the Conley--Floer homology for nongeneric data, we have
  \[
  \HF_*(E,\Phi;\Z_2) \iso \HF_*(E,\Phi_\epsilon,\Psi_\epsilon;\Z_2), \qquad \text{for} \quad E=E_1, \; E_2, \; E_1 \cup E_2.
  \]

  On the level of chain complexes, the splitting of chain groups into the direct sum
  \[
  C_n( E_1 \cup E_2 , \Phi_\epsilon ) \iso C_n( E_1 , \Phi_\epsilon ) \oplus C_n( E_2 , \Phi_\epsilon ),
  \]
  and the boundary operator $\partial_n^\epsilon$ associated with the generic perturbations $\Phi_\epsilon$, $\Psi_\epsilon$ factorises through this splitting. 
  Consequently
  \[
  \HF_n(E_1 \cup E_2,\Phi_\epsilon,\Psi_\epsilon;\Z_2) \iso \HF_n(E_1,\Phi_\epsilon,\Psi_\epsilon;\Z_2) \oplus \HF_n(E_2,\Phi_\epsilon,\Psi_\epsilon;\Z_2).
  \]
  As this holds for any such generic perturbations $\Phi_\epsilon$, $\Psi_\epsilon$, we obtain the splitting of the homology for nongeneric data, as well.
\end{myproof}



\section{A Morse isomorphism}
\label{sec:TWN_Morse_isomorphism}

In this section we formulate a version of a Morse isomorphism for the Conley--Floer homology developed in the previous section.
This isomorphism relates the Conley--Floer homology groups with the singular homology of certain subsets of $\R^d$.
Before discussing this theorem, we first introduce the notion of isolating blocks.

\subsection{Isolating blocks}

In this section we introduce a class of isolating trajectory neighbourhoods, built up from submanifolds of $\R^d$.
Isolating blocks where introduced in \cite{conley1971isolated} in the context of smooth flows.
In \cite{wilson1973lyapunov}, a subclass of isolating blocks, so-called isolating blocks-with-corners, was introduced, which provide slightly more regularity compared to isolating blocks.
Here, we generalise the notion of isolating blocks-with-corners, which for convenience we will refer to as isolating blocks.
To do so, we formally interpret $\Phi$ as a multivalued vector field on $\R^d$, with values attainable in a point $u_0 \in \R^d$ given by $\Phi(u)(0)$,
where $u : \R \to \R^d$ is an arbitrary curve with $u(0) = u_0$.
We point out the obvious similarities to Conley theory for discrete multivalued maps, introduced in \cite{kaczynski1995conley}.

We first introduce a few notions about continuous vector fields.
For a given manifold $M$, a vector field along a continuous curve $\gamma : \R \to M$ is a continuous map $V : \R \to T M$ such that $V(x) \in T_{\gamma(x)} M$.
Given a codimension $1$ submanifold without boundary $A \subset M$, a vector field $V$ along a curve $\gamma$ is said to be transverse to $A$ if,
whenever $\gamma(x) \in A$, the vector $V(x) \in T_{\gamma(x)} M$ is transverse to $T_{\gamma(x)} A$.

Given $u : \R \to \R^d$, we exploit the canonical identification between $\R^d$ and $T_{u(x)} \R^d$ to identify $\Phi(u)$ with a vector field $\phi_u$ along $u$.

\begin{figure}
  \centering
  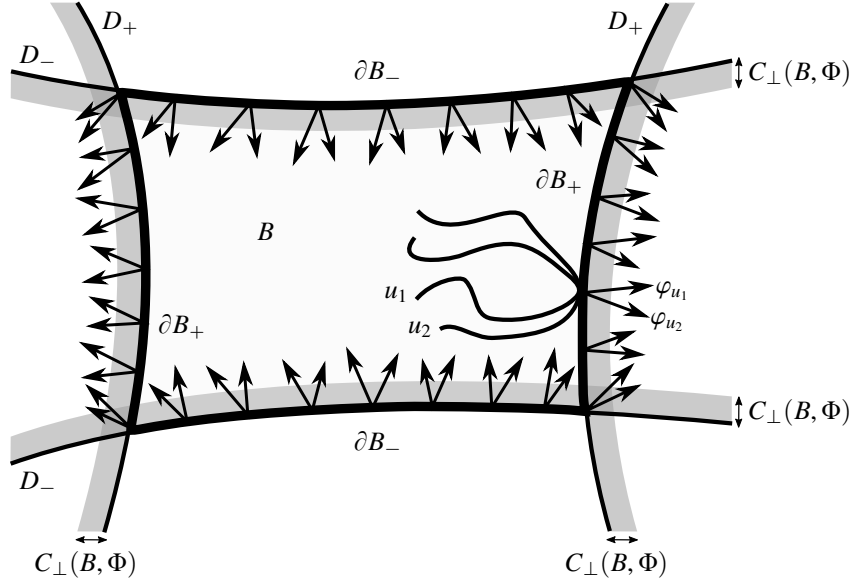
  \caption{An isolating block, where $\Phi$ is interpreted as a multivalued vector field on $\bdy B$.}
  \label{fig:TWN_isolating_block}
\end{figure}

\begin{mydef}[Isolating block]
\label{def:TWN_isolating_block}
  A $d$-dimensional compact manifold with corners $B \subset \R^d$ is said to be an isolating block for \eqref{eq:TWN} provided that the following is true:
  \begin{enumerate}
  \item The boundary of $B$ decomposes as 
    \[
    \bdy B = \bdy B_+ \cup \bdy B_- \cup \bdy B_0,
    \]
    where $\bdy B_\pm$ (when nomemtpy) are compact $(d-1)$-dimensional submanifolds.
   There exist open $(d-1)$-dimensional submanifolds $D_\pm \subset \R^d$ such that $\bdy B_\pm \subset D_\pm$, and the transverse intersection $D_+ \cap D_- = \bdy B_0$ (when nonempty)
    is of dimension $d-2$.
  \item For any given point $u_0 \in \bdy B_0$, any line segment $\ell(u_0)$ passing through $u_0$ can not be contained in $B$, i.e., $\ell(u_0) \not\subseteq B$.
  \item Given any continuous curve $u : \R \to B$,
    the vector field $\phi_u$ along $u$ is transverse to both $D_+$ and $D_-$, inward pointing on $\bdy B_- \setmin \bdy B_0$ and outward pointing on $\bdy B_+ \setmin \bdy B_0$.
    \item There exists an optimal constant $C_\perp(B,\Phi) \geq 0$ so that for any $u : \R \to B$ we have
      \[
      | \phi_u(x) \cdot \nu_\pm(u(x)) | \geq C_\perp(B,\Phi), \qquad \text{whenever} \quad u(x) \in D_\pm.
      \]
      Here $\nu_+$ and $\nu_-$ denote unit normal vector fields (with respect to the Euclidean metric on $\R^d$) along $D_+$ and $D_-$, respectively.
  \end{enumerate}
\end{mydef}
Note that condition (2) implies the isolating block $B$ does not admit any internal tangencies at the corners.
Any solution $u$ of \eqref{eq:TWN} which hits the boundary $\bdy B$ must therefore enter $B$ by passing through $\bdy B_+$ and exit $B$ via $\bdy B_-$.
The sets $\bdy B_+$ and $\bdy B_-$ are called sets of ingress and egress, respectively.
A schematic depiction of an isolating block for \eqref{eq:TWN} is given in Figure \ref{fig:TWN_isolating_block}.

Note that isolating blocks are preserved under small perturbations of $\Phi$, where the smallness is expressed in terms of the constant $C_\perp(B,\Phi)$.
We will exploit this property in Section \ref{subsec:TWN_Morse_isomorphism} and Section \ref{subsec:TWN_global_homology}.

  We now consider the set
  \[
  E(B) := \set{ u \in C^1_{\text{loc}}(\R,\R^d) }{ \cl{ u(\R) } \subset \inter{B} }.
  \]

\begin{mylemma}
\label{lemma:TWN_isolating_block_nbhd}
  Suppose $B$ is an isolating block for \eqref{eq:TWN}.
  Then, for any $u \in \cl{E(B)} \setmin E(B)$, we have
  \[
  \| u' + \Phi(u) \|_{L^\infty(\R,\R^d)} \geq C_\perp(B,\Phi).
  \]
  In particular, $E(B)$ is an isolating trajectory neighbourhood.
\end{mylemma}
\begin{myproof}
  Choose the unit normal vector fields $\nu_+$ and $\nu_-$ along $D_+$ and $D_-$ in such a manner that
  their restrictions to $\bdy B_\pm \setmin \bdy B_0$ are outward pointing with respect to $B$.
  Then
  \[
    \pm \big( \phi_u(x) \cdot \nu_\pm(u(x)) \big) > C_\perp(B,\Phi), \qquad \text{whenever} \quad u(x) \in D_\pm.
  \]
  Now let $u \in \cl{ E(B) } \setmin E(B)$.
  Note that, without loss of generality, we may assume that
  \begin{equation}
    \label{eq:isolating_block_delta_bound}
     \| u' + \Phi(u) \|_{L^\infty(\R,\R^d)} < \infty,
  \end{equation}
  for otherwise the desired estimate is trivially satisfied.

  
  Since $u \in \cl{ E(B) } \setmin E(B)$, there exists a sequence $(x_n)_n \subset \R$ such that $( u(x_n) )_n$ accumulates onto $\bdy B$.
  Let $u_n := x_n \cdot u$.
  In light of \eqref{eq:isolating_block_delta_bound} and the translational invariance of the $L^\infty$ norm, 
  we see that $(u_n)_n$ is a sequence of $\delta$-approximate solutions (for sufficiently large $\delta > 0$).
  Furthermore, as $B$ is compact, we have $\sup_n \| u_n \|_{L^\infty(\R,\R^d)} < \infty$.
  Consequently, using Theorem \ref{thm:TWN_weak_compactness}
  we can extract a subsequence, again denoted by $(u_n)_n$, converging in $C^1_{\text{loc}}(\R,\R^d)$ to $u_\infty \in \cl{ E(B) }$.
  By construction we have $u_\infty(0) \in \bdy B$.

  Let us assume $u_\infty(0) \in D_+$; the argument for $u_\infty(0) \in D_-$ is identical.
  Then
  \[
  \phi_u(x_n) \cdot \nu_+(u(x_n)) = \phi_{u_n}(0) \cdot \nu_+(u_n(0)) \to \phi_{u_\infty}(0) \cdot \nu_+(u_\infty(0)) \geq C_\perp(B,\Phi).
  \]
  Hence, given any $0 < L < C_\perp(B,\Phi)$, for sufficiently large $n$ we have
  \[
  \Phi(u)(x_n) \cdot \nu_+(u_\infty(0)) > L,
  \]
  where we used the canonical identification between $\R^d$ and $T_{u_\infty(0)}$ to identify $\Phi(u)(x_n) \in \R^d$ with a vector in $T_{u_\infty(0)} \R^d$.
  Suppose first that $u_\infty(0) \in \bdy B_+$.
  Then $u_\infty$ must have a tangency to $D_+$ at the point $u_\infty(0)$.
  Hence $u'(x_n) \cdot \nu_+(u_\infty(0)) \to u_\infty'(0) \cdot \nu_+(u_\infty(0)) = 0$.
  On the other hand, if $u_\infty(0) \in \bdy B_0$, then $u_\infty'(0) = 0$ since $B$ does not admit any internal tangencies at the corners.
  So again, we find that $u'(x_n) \cdot \nu_+(u_\infty(0)) \to 0$.
  Therefore, we have
  \[
  \| u' + \Phi(u) \|_{L^\infty(\R,\R^d)} \geq \limsup_{n\to\infty}  \big( u'(x_n) + \Phi(u)(x_n) \big) \cdot \nu_+(u_\infty(0)) > L. 
  \]
  Since this holds true for any $0 < L < C_\perp(B,\Phi)$, we conclude that
  \[
  \| u' + \Phi(u) \|_{L^\infty(\R,\R^d)} \geq C_\perp(B,\Phi).
  \]
\end{myproof}

\subsection{A Morse isomorphism}
\label{subsec:TWN_Morse_isomorphism}

Compared to a general isolating trajectory neighbourhood, 
the big advantage of isolating blocks is that it is relatively straightforward to check whether a given set is an isolating block,
without having prior knowledge about the isolated invariant set $\Inv( E(B) )$.
After all, it only involves estimates on $\phi_u$ along the boundary of a compact set.
The other big advantage is that it allows us to relate the Conley--Floer homology groups to the singular homology of the isolating block $B$ relative to the set of egress $\bdy B_-$.
This Morse isomorphism indicates the similarity between Conley--Floer homology for nonlocally coupled systems on the one hand, 
and on the other hand the homological Conley index for flows.

Before stating and proving a Morse isomorphism for Conley--Floer homology, we note that, in applications, one could use the estimate
\[
\sup_{u \in \cl{E(B)}\setmin E(B)} \| \nabla_g S(u)^T \calN[ S(u) ] \|_{L^\infty(\R,\R^d)} \leq C  \sup_{z_1 , z_2 \in B} \| \nabla_g S(z_1) \| | S(z_2) |,
\]
where the constant $C$ depends on $\calN$, only, so that verification of the hypotheses required in the Morse isomorphism reduce to pointwise considerations.

\begin{mythm}
\label{thm:TWN_Morse_isomorphism}
Consider nonlinearities
\[
\Phi_0(u) = \nabla_g F(u), \qquad \Phi_1(u) = \nabla_g S(u)^T \calN[S(u)] + \nabla_g F(u).
\]
  Suppose $B$ is an isolating block for \eqref{eq:TWN} with $\Phi = \Phi_1$.
  Assume furthermore that
  \[
  \sup_{u \in \cl{E(B)}\setmin E(B)} \| \nabla_g S(u)^T \calN[ S(u) ] \|_{L^\infty(\R,\R^d)} < C_\perp(B,\Phi_0).
  \]
  Then $\Phi_0$ and $\Phi_1$ are related by continuation over $E(B)$, and
  \[
  \HF_*(E(B),\Phi_1;\Z_2) \iso \HF_*(E(B),\Phi_0;\Z_2) \iso \homology_*( B , \bdy B_- ; \Z_2 ),
  \]
  where $\homology_*$ denotes the relative singular homology with $\Z_2$ coefficients.
\end{mythm}
\begin{myproof}
  We will first verify that $\Phi_0$ and $\Phi_1$ are related by continuation over $E(B)$.
  For $0 \leq \beta \leq 1$, let
  \[
  \Phi_\beta(u) := \beta \nabla_g S(u)^T \calN[S(u)] + \nabla_g F(u).
  \]
  By the hypothesis of the theorem there exists $0 < \epsilon \leq 1$ so that for any $u \in \cl{ E(B) }$ we have
  \[
  \| \Phi_\beta(u) - \Phi_0(u) \|_{L^\infty(\R,\R^d)} = \beta \| \nabla_g S(u)^T \calN[ S(u) ] \|_{L^\infty(\R,\R^d)} \leq (1-\epsilon) C_\perp(B,\Phi_0).
  \]
   Consequently, using Lemma \ref{lemma:TWN_isolating_block_nbhd} we recover
  \begin{align*}
    \| u' + \Phi_\beta(u) \|_{L^\infty(\R,\R^d)} &\geq \| u' + \Phi_0(u) \|_{L^\infty(\R,\R^d)} - \| \Phi_\beta(u) - \Phi_0(u) \|_{L^\infty(\R,\R^d)} \\
    &\geq C_\perp(B,\Phi_0) - \| \Phi_\beta(u) - \Phi_0(u) \|_{L^\infty(\R,\R^d)} \\
    &\geq \epsilon C_\perp(B,\Phi_0).
  \end{align*}
  This shows that the maps $\Phi_0$ and $\Phi_1$ are related by continuation over $E(B)$.
  Consequently, by Theorem \ref{thm:TWN_HF_homotopy_invariance},
  \[
  \HF_*(E(B),\Phi_0 ; \Z_2) \iso \HF_*(E(B),\Phi_1 ; \Z_2).
  \]
  
  We will now identify the Conley--Floer homology of $( E(B) , \Phi_0 )$ with the Morse homology of $(B,F,g)$.
  At this stage we perturb $F$ slightly, resulting in a perturbation $\widetilde \Phi_0(u) := \nabla_{g} \widetilde F(u)$ of $\Phi_0$, where 
  $(\widetilde F , g)$ is Morse--Smale on $B$, i.e., $\rst{ F }{ B } : B \to \R$ is a Morse function, 
  and for any $z_0 , z_1 \in \crit(\widetilde F) \cap B$, the stable and unstable manifolds $W^s(z_0)$ and $W^u(z_1)$ of the negative gradient flow
  \begin{equation}
    \label{eq:TWN_local_gradflow}
    z' = - \nabla_{g} \widetilde F(z)
  \end{equation}
  intersect transversely.
  Following results in \cite{austin1995morse}, we may choose such perturbations as small as we wish, 
  thus ensuring that $B$ is also an isolating block for \eqref{eq:TWN} with $\Phi = \widetilde \Phi_0$,
  and
  \[
  \HF_*(E(B),\Phi_0 ; \Z_2) \iso \HF_*(E(B), \widetilde \Phi_0 ; \Z_2).
  \]
  Now recall that $\HF_*(E(B),\widetilde \Phi_0 ; \Z_2)$ is defined as the isomorphism type of 
  \[
  \HF_*(E(B), \widetilde \Phi_0, \Psi ; \Z_2),
  \]
  for generic perturbations $\Psi \in \bPsi_h^\pitchfork(E(B),\widetilde \Phi_0)$.
  By the Morse--Smale property, we have $0 \in \bPsi_h^\pitchfork(E(B),\widetilde \Phi_0)$.
  Hence the boundary operator associated with $\HF_*(E(B),\widetilde \Phi_0 ; \Z_2)$ counts the flow lines of the gradient flow \eqref{eq:TWN_local_gradflow}.
  We thus recover the Morse homology of $(B, F , g)$, i.e.,
  \[
  \HF_*(E(B),\Phi_0 ; \Z_2) \iso \HM_*(B,F,g;\Z_2).
  \]

  Combining our observations up to this point, we have the identification
  \[
  \HF_*(E(B),\Phi_1 ; \Z_2) \iso \HM_*(B,F,g;\Z_2).
  \]
  The fact that the Morse homology computes the singular homology of $B$ relative to the ``exit set'' $\bdy{B_- }$ is a standard result in Morse theory \cite{schwarz1993morse}.
\end{myproof}

\subsection{Global Conley--Floer homology}
\label{subsec:TWN_global_homology}

We now introduce a notion of stability of isolating blocks.
\begin{mydef}[Stabilising isolating blocks]
  A sequence $( B^R )_{R \geq R_0}$ of is said to be a family of stabilising isolating blocks for \eqref{eq:TWN}, provided that the following holds.
  \begin{enumerate}
  \item For each $R \geq R_0$, the set $B^R$ is an isolating block for \eqref{eq:TWN}.
  \item The isolating blocks exhaust $\R^d$, in the sense that
    \[
    B^{R_1} \subseteq B^{R_2} \qquad \text{whenever} \quad R_0 \leq R_1 \leq R_2, \qquad \text{and} \quad \bigcup_{R \geq R_0} B^R = \R^d.
    \]
  \item The invariant set captured by the isolating blocks is stable, i.e.,
    \[
    \Inv(E(B^{R_1})) = \Inv(E(B^{R_2})), \qquad \text{for any} \quad R_1 , R_2 \geq R_0.
    \]
  \end{enumerate}
\end{mydef}
Note that, in particular, by Lemma \ref{lemma:TWN_independence_isolating_nbhd} the homology $\HF_*(E(B^R),\Phi;\Z_2)$ is up to canonical isomorphism independent of $R \geq R_0$.
This leads us to define the global Conley--Floer homology as the isomorphism type of $\HF_*(E(B^{R_0'},\Phi;\Z_2)$.
\begin{mydef}[Global Conley--Floer homology]
  Suppose \eqref{eq:TWN} admits a family of stabilising isolating blocks $(B^R)_{R \geq R_0}$.
  We then define the global Conley--Floer homology of \eqref{eq:TWN} by
\[
\HF_*^{\text{global}}(\Phi;\Z_2) := \varprojlim \HF_*(E(B^R),\Phi_\epsilon;\Z_2).
\]
  where the inverse limit is taken over $R \geq R_0$ and the associated canonical isomorphisms of the homology groups.
\end{mydef}

After Definition \ref{def:TWN_isolating_block} we already mentioned that isolating blocks are stable with respect to small perturbations, 
where smallness is expressed in terms of the constant $C_\perp(B,\Phi)$.
Given a family $(B^R)_{R\geq R_0}$ of stabilising isolating blocks, the perturbation then only has to be small in terms of $C_\perp(B^R,\Phi)$, for $R$ sufficiently large.
Stated differently, families of isolating blocks are stable with respect to perturbations which are ``small at infinity''.
This is the content of the next lemma.

\begin{mylemma}
\label{lemma:TWN_stabalising_blocks}
Consider nonlinearities $\Phi_0$ and $\Phi_1$ as defined in Section \ref{sec:TWN_setup}.
  Suppose $(B^R)_{R\geq R_0}$ is a family of stabilising isolating blocks for \eqref{eq:TWN} with $\Phi = \Phi_0$.
  Assume furthermore that  
  \[
  \limsup_{R \to \infty} \frac{ \sup_{u \in \cl{E(B^R)}\setmin E(B^R)} \| \Phi_1(u) - \Phi_0(u) \|_{L^\infty(\R,\R^d)} }{ C_\perp(B^R,\Phi_0) } < 1.
  \]
  Then there exists $R_0' \geq R_0$ so that $(B^R)_{R \geq R_0'}$ is a family of stabilising isolating blocks for \eqref{eq:TWN} with $\Phi = \Phi_1$.
\end{mylemma}
\begin{myproof}
Given $R \geq R_0$, let $u : \R \to B^R$ be a continuous curve, and identify $\Phi_0(u)$ and $\Phi_1(u)$ with vector fields $\phi_u^0$ and $\phi_u^1$ along $u$, respectively.
  Using the hypothesis of the lemma, we find that there exists $0 \leq \epsilon < 1$ and $R_0' \geq R_0$, such that for any $R \geq R_0'$ and any continuous curve $u \in \cl{E(B^R)} \setmin E(B^R)$ we have
  \[
  \| \phi_u^1  - \phi_u^0 \|_{L^\infty(\R,\R^d)}  = \| \Phi_1(u) - \Phi_0(u) \|_{L^\infty(\R,\R^d)} \leq \epsilon C_\perp(B^R,\Phi_0).
  \]
  Fix $R \geq R_0'$ and denote by $D_{R,\pm}$ the $(d-1)$-dimensional manifolds from the definition of isolating blocks, 
  and choose unit normal vector fields $\nu_{R,\pm}$ on $D_{R,\pm}$ (the choice of orientation is irrelevant).
  Let $u : \R \to B^R$ be a continuous curve and $x \in \R$ be such that $u(x) \in D_{R,\pm}$.
  Then
  \begin{align*}
    | \phi_u^1(x) \cdot \nu_{R,\pm}(u(x)) | &\geq | \phi_u^0(x) \cdot \nu_{R,\pm}(u(x)) | - | ( \phi_u^1(x) - \phi_u^0(x) ) \cdot \nu_{R,\pm}(u(x)) | \\
    &\geq C_\perp(B^R,\Phi_0) - \| \phi_u^1  - \phi_u^0 \|_{L^\infty(\R,\R^d)} \\
    &\geq (1-\epsilon) C_\perp(B^R,\Phi_0).
  \end{align*}
  Hence, according to Definition \ref{def:TWN_isolating_block}, $B^R$ is an isolating block for \eqref{eq:TWN} with $\Phi = \Phi_1$.
\end{myproof}

We are now prepared to formulate and prove a global version of the Morse isomorphism for the Conley--Floer homology groups of \eqref{eq:TWN}.
As with the local Morse isomorphism, we remark that verification of the hypotheses required in the global Morse isomorphism reduce to pointwise considerations,
by using the estimate
\[
\sup_{u \in \cl{E(B^R)}\setmin E(B^R)} \| \nabla_g S(u)^T \calN[ S(u) ] \|_{L^\infty(\R,\R^d)} \leq C  \sup_{z_1 , z_2 \in B^R} \| \nabla_g S(z_1) \| | S(z_2) |,
\]
where the constant $C$ depends on $\calN$, only.

\begin{mythm}
\label{thm:TWN_global_Morse_isomorphism}
Consider nonlinearities
\[
\Phi_0(u) = \nabla_g F(u), \qquad \Phi_1(u) = \nabla_g S(u)^T \calN[S(u)] + \nabla_g F(u).
\]
  Suppose $(B^R)_{R\geq R_0}$ is a family of stabilising isolating blocks for \eqref{eq:TWN} with $\Phi = \Phi_0$.
  Assume furthermore that  
  \[
  \limsup_{R \to \infty} \frac{ \sup_{u \in \cl{E(B^R)}\setmin E(B^R)} \| \nabla_g S(u)^T \calN[ S(u) ] \|_{L^\infty(\R,\R^d)} }{ C_\perp(B^R,\Phi_0) } < 1.
  \]
  Then there exists $R_0' \geq R_0$ so that $(B^R)_{R \geq R_0'}$ is a family of stabilising isolating blocks for \eqref{eq:TWN} with $\Phi = \Phi_1$.
  Furthermore,
  \[
  \HF_*^{\text{global}}( \Phi_1 ; \Z_2 ) \iso \homology_*( B^R , \bdy B^R_- ; \Z_2 ),
  \]
  where $\homology_*$ denotes the relative singular homology with $\Z_2$ coefficients. 
\end{mythm}
\begin{myproof}
  By Lemma \ref{lemma:TWN_stabalising_blocks} there exists $R_0' \geq R_0$ so that $(B^R)_{R\geq R_0'}$ is a family of stabilising isolating blocks for \eqref{eq:TWN} with $\Phi = \Phi_1$.
  After updating the value of $R_0'$ if necessary, we may assume
  \[
   \sup_{u \in \cl{E(B^{R_0'})}\setmin E(B^{R_0'})} \| \nabla_g S(u)^T \calN[ S(u) ] \|_{L^\infty(\R,\R^d)} < C_\perp(B^{R_0'},\Phi_0).
  \]
  It now follows from Theorem \ref{thm:TWN_Morse_isomorphism} that
  \[
  \HF_*( E(B^{R_0'}) , \Phi_1 ; \Z_2 ) \iso \homology_*( B^{R_0'} , \bdy B^{R_0'}_- ; \Z_2 ).
  \]
  This, combined with the definition of the global Conley--Floer homology, concludes the proof.
\end{myproof}


\section{Application to global bifurcation problems}
\label{sec:TWN_application}

In this section we formulate a ``forcing theorem'', deriving lower bounds on the existence and multiplicity of travelling waves in nonlocal equations.
We finish with giving examples of how this theory can be applied to study travelling waves in neural field equations and nonlocal diffusion equations.

\subsection{Detection of heteroclinics}

Let $E$ be an isolating trajectory neighbourhood for \eqref{eq:TWN} with $\Psi = 0$,
and let $Z \subset E$ consist of hyperbolic constant solutions of \eqref{eq:TWN}.
We do not assume that $Z$ contains all the hyperbolic constant solutions of \eqref{eq:TWN}.
Hence such sets can be constructed in applications, either by hand or through rigorous computational means.

Recall from Section \ref{subsec:TWN_asymptotics} the definition of $\alpha$- and $\omega$-limit sets, and define
\begin{align*}
  W^-_E(Z) &:= \set{ \text{ solutions $u$ of \eqref{eq:TWN} in $E$ } }{ \alpha(u) \in Z } / \R, \\
  W^+_E(Z) &:= \set{ \text{ solutions $u$ of \eqref{eq:TWN} in $E$ } }{ \omega(u) \in Z } / \R,
\end{align*}
where the $\R$-action is translation $\tau \mapsto \tau \cdot u$.
We denote the equivalence class of $u$ in $W^-_E(Z)$ or $W^+_E(Z)$ by $[u]$.

Elements $[u]$ in the intersection $W^-_E(Z) \cap W^+_E(Z)$ have the property that $u$ is a heteroclinic solution of \eqref{eq:TWN}, connecting distinct hyperbolic endpoints selected from $Z$.
The elements $[u]$ of the symmetric difference $W^-_E(Z) \triangle W^+_E(Z)$ have either $\alpha(u) \in Z$, or $\omega(u) \in Z$, but not both.
The function $u$ is a generalised heteroclinic solutions of \eqref{eq:TWN}, conform Remark \ref{remark:TWN_generalised_heteroclinics}.
More precisely, the limit sets $\alpha(u)$ and $\omega(u)$ are nonempty, disjoint, and consist of critical points of $h$.
However, those critical points might be degenerate.

We can now exploit the properties of the homology to obtain a ``forcing theorem'', 
giving lower bounds on the number of heteroclinics contained in an isolating trajectory neighbourhood,
in terms of prior information on constant solutions and the homology groups.

\begin{mythm}
\label{thm:TWN_Morse_inequality}
Let $E$ be an isolating trajectory neighbourhood for \eqref{eq:TWN} with $\Psi = 0$.
Let $Z \subset E$ consist of hyperbolic constant solutions of \eqref{eq:TWN} in $E$.
  We then have the following Morse inequality:
  \[
  2 \cdot \# [ W^-_E(Z) \cap W^+_E(Z) ] + \# [ W^-_E(Z) \triangle W^+_E(Z) ] \geq \# Z - \rank \HF_*(E,\Phi;\Z_2).
  \]
  In particular
  \begin{multline*}
      2 \cdot \# \big\{ \text{ generalised heteroclinic solutions of \eqref{eq:TWN} in $E$, modulo translation } \big\} \\
      \geq \# \big\{ \text{ hyperbolic constant solutions of \eqref{eq:TWN} in $E$ } \big\} - \rank \HF_*(E,\Phi;\Z_2).
  \end{multline*}
\end{mythm}
\begin{myproof}
  Let us write $r(E) := \rank \HF_*(E,\Phi;\Z_2)$.
  Suppose $C \subseteq Z$ is a subset of $Z$ such that
  \begin{equation}
    \label{eq:C_rank_estimate}
    \# C - r(E) > 0.
  \end{equation}
  We claim that this inequality ensures that $C$ is not isolated in $\Inv(E)$ with respect to the topology of $L^\infty(\R,\R^d)$.
  We argue by contradiction, assuming that $C$ is isolated in $\Inv(E)$.
  Let us now enumerate $C = \{ z_1 , \dots , z_m \}$.
  By Lemma \ref{lemma:TWN_isolated_invset_nbhd}, the isolation property ensures that we may find mutually disjoint subsets 
  $W , V_1 , \dots , V_m \subset E$, each set being an isolating trajectory neighbourhood of \eqref{eq:TWN},
  with 
  \[
  \Inv(V_i) = \{ z_i \}, \qquad i \in \{1,\dots,m\},
  \]
  and
  \[
  \Inv(E) = \Inv(W) \cup \bigcup_{i=1}^{m} \Inv(V_i).
  \]
  By Lemma \ref{lemma:TWN_HF_direct_sum} we have
  \[
  \HF_n(E,\Phi;\Z_2) \iso \HF_n(W,\Phi;\Z_2) \oplus \bigoplus_{i=1}^m \HF_n(V_i,\Phi;\Z_2).
  \]
  But now the rank of the homology on the right hand side is at least $m = \# C \geq r(E) + 1$, whilst by definition the homology on the left hand side is of rank $r(E)$.
  Since we arrived at a contradiction, we conclude that $C$ is not isolated in $\Inv(E)$.
  More precisely, we find that given any subset $C \subseteq Z$ satisfying \eqref{eq:C_rank_estimate},
  there exists a solution $u \in E$ of \eqref{eq:TWN} for which $\alpha(u) \subseteq C$ or $\omega(u) \subseteq C$, or possibly both.

  We will use this observation to detect the heteroclinic solutions using an iterative procedure, 
  where after $k$ steps we denote by $n^{\cap}_k$ and $n^{\triangle}_k$ the number of detected elements of $W^-_E(Z) \cap W^+_E(Z)$ and $W^-_E(Z) \triangle W^+_E(Z)$, respectively.
  We initialise these variables by setting $n^{\cap}_0 = n^{\triangle}_0 = 0$.
  We start the iterative construction by setting $C_1 = Z$.
  If $\# Z - r(E) = 0$ the conclusion of the theorem is trivial, so that henceforth without loss of generality we may assume that
  \[
  \# C_1 - r(E) = \# Z - r(E) > 0.
  \]
  As per the observation in the previous paragraph (with $C = C_1$), 
  we conclude that there exists a solution $u_1 \in E$ of \eqref{eq:TWN} for which $\alpha(u_1) \subseteq C_1$ or $\omega(u_1) \subseteq C_1$, or possibly both.
  Now set
  \begin{equation}
    \label{eq:Ck_iteration}
    C_{k+1} := C_k \setmin \big( C_k \cap ( \alpha( u_k ) \cup \omega( u_k ) ) \big),
  \end{equation}
  first for $k = 1$, and then iteratively for all consecutive  values of $k \geq 1$.
  If $\# C_{k+1} - r(E) > 0$ we may repeat the preceding argument (with $C = C_{k+1}$) to obtain a solution $u_{k+1} \in E$ of \eqref{eq:TWN} 
  for which $\alpha(u_{k+1}) \subseteq C_{k+1}$ or $\omega(u_{k+1}) \subseteq C_{k+1}$, or possibly both.
  Observe that, in light of the definition \eqref{eq:Ck_iteration} of the spaces $C_1,\dots,C_{k+1}$, the detected solutions $[u_1],\dots,[u_{k+1}]$ must be distinct. 
  
  To keep track of the number and type of heteroclinic solutions detected, we let
  \[
  \begin{pmatrix}
    n_k^{\cap} \\ n_k^{\triangle} 
  \end{pmatrix}
  =
  \left\{
    \begin{array}{l l}
      \begin{pmatrix}
        n_{k-1}^{\cap} + 1 \\
        n_{k-1}^{\triangle}
      \end{pmatrix}
   & \text{if} \quad [u_k] \in W^-_E(Z) \cap W^+_E(Z), \\ \\
      \begin{pmatrix}
        n_{k-1}^{\cap} \\
        n_{k-1}^{\triangle} + 1
      \end{pmatrix}
   & \text{if} \quad [u_k] \in W^-_E(Z) \triangle W^+_E(Z).
    \end{array}
  \right.
  \]
  Hence after $k$ iterations we have detected $n^{\cap}_k$ distinct elements of $W^-_E(Z) \cap W^+_E(Z)$ and $n^{\triangle}_k$ distinct elements of $W^-_E(Z) \triangle W^+_E(Z)$.
  Recalling \eqref{eq:Ck_iteration}, we see that when $ n^{\triangle}_k - n^{\triangle}_{k-1} = 1$, 
  i.e. when $[u_k] \in W^-_E(Z) \triangle W^+_E(Z)$, the space $C_{k+1}$ is obtained from $C_k$ by removing a single element.
  On the other hand, when $n^{\cap}_k - n^{\cap}_{k-1} = 1$, i.e. when $[u_k] \in W^-_E(Z) \cap W^+_E(Z)$, then $C_{k+1}$ is obtained by removing either $1$ or $2$ elements from $C_k$,
  depending on whether $[u_k] \not\in W^-_E(C_k) \cap W^+_E(C_k)$ or $[u_k] \in W^-_E(C_k) \cap W^+_E(C_k)$, respectively.
  As such, we obtain the bound
  \[
  \# C_{k} \setmin C_{k+1} \leq 2 ( n^{\cap}_k - n^{\cap}_{k-1} ) + n^{\triangle}_k - n^{\triangle}_{k-1},
  \]
  where the inequality is strict if $[u_k] \in W^-_E(Z) \cap W^+_E(Z)$ but $[u_k] \not\in W^-_E(C_k) \cap W^+_E(C_k)$.
  Hence
  \[
    \# C_1 = \# C_{k+1} + \# C_k \setmin C_{k+1} + \cdots + \# C_2 \setmin C_1 \leq \# C_{k+1} + 2 n^{\cap}_k + n^{\triangle}_k,
  \]
  and using the fact that $Z = C_1$ we now find that
  \begin{equation}
    \label{eq:TWN_forcing_bound}
    \# C_{k+1} - r(E)  + 2 n^{\cap}_k + n^{\triangle}_k \geq \# Z - r(E).
  \end{equation}
  As long as $\# C_{k+1} - r(E) > 0$ we can proceed with another iteration.
  Since the cardinality of $C_k$ decreases with each step, the iteration must terminate after finitely many steps $k = k_*$, at which point $\# C_{k_* + 1} - r(E) \leq 0$.
  With the aid of \eqref{eq:TWN_forcing_bound} we then find that
  \begin{align*}
    2 \cdot \# [ W^-_E(Z) \cap W^+_E(Z) ] + \# [ W^-_E(Z) \triangle W^+_E(Z) ] &\geq 2 n^{\cap}_{k_*} + n^{\triangle}_{k_*} \geq \# Z - r(E) - ( \# C_{k_*+1} - r(E) ) \\
    &\geq \# Z - r(E).
  \end{align*}
  This concludes the proof.
\end{myproof}

Of course, when more detailed information on $\Phi$ is available, the above forcing theorem can be combined with bounds on Morse indices and energy estimates
to obtain more precise information about the detected heteroclinics.

\subsection{Global structure of travelling fronts}
\label{subsec:TWN_wave_application}

We now return to Examples \ref{ex:TWN_NFE} and \ref{ex:TWN_nonlocal_diffusion} from Section \ref{sec:TWN_intro}, and compute the associated global Conley--Floer homologies.

\subsubsection{Neural field equations}
Travelling fronts, propagating with wave speed $c \neq 0$, in neural field equations
\[
  u_t = -u + N*\sigma(u), \qquad u(t,x) \in \R^d,
\]
satisfy
\begin{equation}
  \label{eq:TWN_NFE_TW}
  c u' + N*\sigma(u) - u = 0.
\end{equation}
Here $\sigma(u) = (\sigma_1(u_1),\dots,\sigma_d(u_d))$, where $\sigma_i : \R \to \R$ are assumed to be bounded and strictly increasing,
and the matrix-valued interaction kernel $N$ is assumed to be symmetric.
We first cast this equation into the form \eqref{eq:TWN}.
Let
\[
S(u) =
\begin{pmatrix}
  c^{-1} \sigma_1(u_1) \\ \vdots \\ c^{-1} \sigma_d(u_d)
\end{pmatrix}, \qquad u = (u_1,\dots,u_d) \in \R^d.
\]
Define a Riemannian metric on $\R^d$ by setting
\[
g_u( v_1 , v_2 ) = \langle v_1 , DS(u) v_2 \rangle = c^{-1} \sum_{i=1}^d \sigma_i'(u_i) v_{1,i} v_{2,i}, \qquad v_1 , v_2 \in T_u \R^d = \R^d,
\]
so that $\nabla_g S(u) = \id$ for any $u \in \R^d$.
Finally, let $F:\R^d \to \R$ be given by
\[
F(u) := - c^{-2} \sum_{i=1}^d \int_0^{u_i} \sigma_i'(s) s \d s, \qquad u = (u_1,\dots,u_d) \in \R^d.
\]
Then $\nabla_g F(u) = - c^{-1} u$, so that
\[
\Phi_1(u) := \nabla_g S(u)^T \calN[S(u)] + \nabla_g F(u) = c^{-1} \big( N*\sigma(u) - u \big),
\]
thus casting \eqref{eq:TWN_NFE_TW} into the form \eqref{eq:TWN} with $\Phi = \Phi_1$.

Given $R > 0$, let $B^R$ denote the closed ball of radius $R$ in $\R^d$.
Then $\nabla_g F$ is transverse and inward pointing on $\bdy B^R$.
In particular, $(B^R)_{R > 0}$ is a family of stabilising isolating neighbourhoods for \eqref{eq:TWN} with $\Phi_0(u) = \nabla_g F(u)$,
and, with the notation of Section \ref{sec:TWN_Morse_isomorphism}, we have $C_\perp(B^R,\Phi_0) = c^{-1} R$.
Since by assumption $S$ uniformly bounded, we find that there exists $C > 0$ such that
\[
\frac{ \sup_{u \in \cl{E(B^R)}\setmin E(B^R)} \| \nabla_g S(u)^T \calN[ S(u) ] \|_{L^\infty(\R,\R^d)} }{ C_\perp(B^R,\Phi_0) } \leq \frac{ C }{ R } \to 0, \qquad \text{as} \quad R \to \infty.
\]
Hence Theorem \ref{thm:TWN_global_Morse_isomorphism} is applicable.
We conclude there exists $R_0' > 0$ so that $(B^R)_{R\geq R_0'}$ is a family of isolating trajectory neighbourhoods for \eqref{eq:TWN} with $\Phi = \Phi_1$, and
  \[
  \HF_*^{\text{global}}( \Phi_1 ; \Z_2 ) \iso \homology_*( B^R , \bdy B^R_- ; \Z_2 ), \qquad \text{where} \quad R \geq R_0'.
  \]
The relative homology on the right hand side computes the reduced singular homology of a $d$-sphere, hence
\[
\HF_n^{\text{global}}( \Phi_1 ; \Z_2 ) \iso \left\{
  \begin{array}{l l}
    \Z_2 & \text{if } n = d, \\
    0 & \text{otherwise}.
  \end{array}
  \right.
\]
Now Theorem \ref{thm:TWN_Morse_inequality} proves existence (and multiplicity) of travelling fronts for any wave speed $c \neq 0$, in the neural field equation with symmetric interaction kernel,
provided the neural field equation supports at least $2$ spatially homogeneous hyperbolic stationary states.

\subsubsection{Nonlocal reaction-diffusion}
Travelling fronts, propagating with wave speed $c \neq 0$, in the nonlocal reaction-diffusion equations of gradient form
\[
u_t = N*u + \nabla F(u), \qquad u(t,x) \in \R^d,
\]
where the interaction kernel $N$ is symmetric and may be continuous, discrete, or a combination of these two, satisfy
\begin{equation}
  \label{eq:TWN_NRD_TW}
    c u' + N*u + \nabla F(u) = 0.
\end{equation}
After rescaling we may assume $c = 1$, and it becomes obvious that this equation is of the type \eqref{eq:TWN}.

Assuming $\nabla F$ grows superlinear, the isomorphism type of the global Conley--Floer homology now depends on the asymptotics of $F$.
In the scalar case, $d = 1$, the global homology is always of rank $0$ or $1$.
The situation drastically changes when considering higher dimensional systems.
Even when $d = 2$, the global homology may be of arbitrary rank.
To see this, consider, for given $k \in \N$ and $n > 2$, a potential $F$ of the form 
\begin{equation}
  \label{eq:TWN_NRD_example_potential}
  F(u) = F_{k,n}(u_1,u_2) + P_{n-1}(u_1,u_2), \qquad u = (u_1,u_2) \in \R^2,
\end{equation}
where $P_{n-1}$ is an arbitrary polynomial of degree at most $n-1$,
and $F_{k,n}$ is given in polar coordinates $(\rho,\theta)$ by
\[
F_{k,n}(\rho,\theta) = \rho^n \sin( k \theta ), \qquad \rho \geq 0, \; \theta \in [0,2\pi).
\]

We will now construct families of stabilising isolating blocks.
First suppose $k \geq 2$ .
Let $B^R \subset \R^2$ be the filled polygon with $2k$ vertices $V^R_0, \dots , V^R_{2k-1}$, given in polar coordinates by
\[
V^R_i = \big( R , i \pi / k \big), \qquad 0 \leq i \leq 2k - 1.
\]
Given $0 \leq i \leq 2k - 2$, denote by $Q^R_i$ the edge connecting $V^R_i$ with $V^R_{i+1}$;
denote by $V^R_{2k-1}$ the edge connecting $V^R_{2k-1}$ with $V^R_0$.
Then $\nabla F_{k,n}$ is outward pointing on $Q^R_i$ for $i$ even, and inward pointing for $i$ odd.
Since $\nabla P_{n-1}$ is a lower order perturbation of $\nabla F_{k,n}$, we find that for sufficiently large $R$ the set $B^R$ is an isolating block for \eqref{eq:TWN} with $\Phi = \Phi_0$,
where $\Phi_0(u) = \nabla F(u)$, with sets of ingress, egress, and corners of $\bdy B^R$ given by
\[
\bdy B^R_+ = \bigcup_{0\leq j \leq k-1} Q^R_i, \qquad \bdy B^R_- = \bigcup_{0\leq j \leq k-1} Q^R_{i+1}, \qquad \bdy B^R_0 = \bigcup_{0\leq j \leq 2k-1} V^R_i.
\]

When $k = 1$, the family of isolating blocks needs to be defined slightly different.
We let $B^R \subset \R^2$ be given by
\[
B^R = \set{ (x,y) \in \R^2 }{ - R \leq x \leq R, \quad | y | \leq | x^2 - R^2 | }.
\]
Then $\nabla F_{k,n}$ is outward pointing on $\bdy B^R_+ := \bdy B^R \cap \R \times [0,\infty)$, and inward
pointing on $\bdy B^R_- := \bdy B^R \cap \R \times (-\infty,0]$.
Since $\nabla P_{n-1}$ is a lower order perturbation of $\nabla F_{k,n}$, we find that for sufficiently large $R$ the set $B^R$ is an isolating block for \eqref{eq:TWN} with $\Phi = \Phi_0$.

We now note that, regardless of the value of $k \geq 1$, we have
\[
\liminf_{R\to\infty} \frac{ C_\perp(B^R,\Phi_0) }{ R^{n-1} } > 0.
\]
By assumption $n > 2$, so that we obtain the estimate
\[
\frac{ \sup_{u \in \cl{E(B^R)}\setmin E(B^R)} \| \calN[ u ] \|_{L^\infty(\R,\R^d)} }{ C_\perp(B^R,\Phi_0) } \leq \frac{ C R }{ C_\perp(B^R,\Phi_0) } \to 0, \qquad \text{as} \quad R \to \infty.
\]
Hence by Theorem \ref{thm:TWN_global_Morse_isomorphism} there exists $R_0' > 0$ so that $(B^R)_{R\geq R_0'}$ is a family of isolating trajectory neighbourhoods for \eqref{eq:TWN} 
with $\Phi = \Phi_1$, where $\Phi_1(u) = \calN[u] + \nabla_g F(u)$, and
  \[
  \HF_*^{\text{global}}( \Phi_1 ; \Z_2 ) \iso \homology_*( B^R , \bdy B^R_- ; \Z_2 ), \qquad \text{where} \quad R \geq R_0'.
  \]
  The homotopy type of $(N , \bdy N_-)$ is that of a marked $2$-sphere with $k \geq 1$ punctures, and we find that
\[
\HF_n^{\text{global}}( \Phi_1 ; \Z_2 ) \iso \left\{
  \begin{array}{l l}
    \underbrace{ \Z_2 \oplus \cdots \oplus \Z_2 }_{ k-1 \text{ times} } & \text{if } n = 1, \\
    0 & \text{otherwise}.
  \end{array}
  \right.
\]
Now by Theorem \eqref{thm:TWN_Morse_inequality}, there exist travelling fronts for any wave speed $c \neq 0$ in the nonlocal reaction-diffusion equation
with potential $F$ given by \eqref{eq:TWN_NRD_example_potential}, provided that the reaction-diffusion equation has at least $k$ spatially homogeneous hyperbolic steady states.





\bibliography{biblio}
\bibliographystyle{abbrv}

\end{document}